%% file: main.tex
\documentclass[preprint]{elsarticle}

\usepackage[margin=3cm]{geometry}

\usepackage{lineno,hyperref}
\modulolinenumbers[5]
\usepackage{amssymb}

\usepackage{subfigure}

\usepackage{amsmath} 
\usepackage{color} 
\usepackage{mathtools} 

\usepackage[]{algorithm2e} 

\makeatletter
\def\ps@pprintTitle{%
 \let\@oddhead\@empty
 \let\@evenhead\@empty
 \def\@oddfoot{}%
 \let\@evenfoot\@oddfoot}
\makeatother

\begin{document}

\begin{frontmatter}

\title{Simulation and validation of surfactant-laden drops in two-dimensional Stokes flow}
\author[add1]{Sara P\aa lsson \corref{cor1}}
\ead{sarapal@kth.se}
\author[add2]{Michael Siegel}
\author[add1]{Anna-Karin Tornberg}

\cortext[cor1]{Corresponding author}
\address[add1]{KTH Mathematics, Linn\'e Flow Centre, 100 44 Stockholm Sweden}
\address[add2]{Department of Mathematical Sciences, New Jersey Institute of Technology, Newark, NJ 07102}

\begin{abstract}
Performing highly accurate simulations of droplet systems is a challenging problem. This is primarily due to the interface dynamics which is complicated further by the addition of surfactants. This paper presents a boundary integral method for computing the evolution of surfactant-covered droplets in 2D Stokes flow. The method has spectral accuracy in space and the adaptive time-stepping scheme allows for control of the temporal errors. Previously available semi-analytical solutions (based on conformal-mapping techniques) are extended to include surfactants, and a set of algorithms is introduced to detail their evaluation.  These semi-analytical solutions are used to validate and assess the accuracy of the boundary integral method, and it is demonstrated that the presented method maintains its high accuracy even when droplets are in close proximity.
\end{abstract}


\end{frontmatter}

\input{srcfiles/intro5}
\input{srcfiles/formulation}
\input{srcfiles/bie}
\input{srcfiles/meth_bie}

\input{srcfiles/meth_est}
\input{srcfiles/meth_spec}
\input{srcfiles/nummethods}
\input{srcfiles/validation4}

\input{srcfiles/results}
\input{srcfiles/conclusions2}
\input{srcfiles/acknowledgement}

\appendix
\input{srcfiles/appendix_est}
\input{srcfiles/appendix_spec}

\newpage

\section*{References}
\bibliographystyle{plainnat}
\bibliography{paperI_ref.bib}

\end{document}

%% file: srcfiles/intro5.tex

\section{Introduction}
\noindent Microfluidics is the study of fluids at the microscopic scale. It is a field of study with a range of applications \cite{droplet_review}, such as: drug delivery, diagnostic chips and microreactors. Droplets in microfluidic systems have a large surface-area to volume ratio, which makes interfacial forces important; a good review is provided by \citeauthor{Rallison1984} \cite{Rallison1984}. Surfactants locally influence the surface tension of a droplet interface. This may create non-uniformity in surface tension which in turn generates a stress opposed to the flow, known as the Marangoni stress. Thus, the addition of surfactants can strongly influence the behaviour of the system \cite{anna2016,droplet_reviewTan}.

\citeauthor{Stone1990}
\cite{Stone1990} studied how surfactants change the drop deformation and breakup using a numerical boundary integral method. They found that for the same strain rate, surfactant-covered drops would reach a more deformed steady state than their clean counterparts. The effect of viscosity-ratio on deformation was investigated for surfactant-covered droplets in an extension of this work \cite{Milliken1993a}. It was shown that for cases with near-zero diffusion of surfactants along the interface of the droplet, the steady state deformation was independent of viscosity-ratio.

A microfluidic system is typically one in which viscous forces are dominant. Microfluidic multiphase flows can therefore be accurately modelled using Stokes equations. With the addition of surfactants, it is also necessary to couple Stokes equations with a convection-diffusion equation on the droplet interfaces. A discontinuity in normal stress across the droplet interface is created by the surface tension forces. This represents a major challenge for the simulation of these flow types. One way to avoid this difficulty is to ensure that the computational grid is always coincident with the interfaces. However, this becomes very expensive when droplets deform or advect as it requires the domain to be remeshed at every time-step. The most common approach for multiphase flow simulation is instead to utilise some form of interface tracking. This involves solving Stokes equations on a fixed computational grid in the entire domain. The locations of droplet interfaces are computed separately using either an explicit boundary discretisation or an implicit representation; such as a level-set function or volume-of-fluid method. Some examples of this kind of methods which also deal with surfactants include those by \citeauthor{khatri_soluble} \cite{Khatri2011,khatri_soluble} and \citeauthor{Muradoglu2014} \cite{Muradoglu2014}, who considered the Navier-Stokes equations and surfactants using a front-tracking scheme and finite differences in 2D and 3D respectively. Another approach is to use a diffuse interface method as
\citeauthor{Teigen2011} \cite{Teigen2011}.

Accurate treatment of the discontinuity in normal stress remains a significant challenge for all interface tracking methods when the interfaces do not align with the computational grid. The most common approach to handle this has been to regularise the surface tension forces, such as introduced in e.g. \citeauthor{Brackbill1992} \cite{Brackbill1992} or in the immersed boundary method by \citeauthor{Peskin1977} \cite{Peskin1977}. However, regularisation limits the accuracy of the method to first order near the interfaces. Another approach is to impose the jump conditions in the normal stress directly. \citeauthor{Leveque1994} \cite{Leveque1994} obtained second order accuracy using an immersed interface method based on finite differences for clean drops, but required drops and bulk to have the same viscosity. This constraint still remained as \citeauthor{xu2006level} \cite{xu2006level} extended the method in \cite{Leveque1994} to simulate surfactant-covered droplets in Stokes flow. Care has to be taken when imposing the jump condition to avoid  a time step constraint that is dependent on how the interfaces cut the underlying grid. This is thoroughly discussed for finite element methods in  \citeauthor{Hansbo2014} \cite{Hansbo2014}, where they manage to design a second order method and avoid such a constraint by adding suitable stabilising terms.

Reformulating Stokes equations in integral form on the droplet interfaces avoids the aforementioned issues. This reduces the dimension of the problem as only the interfaces have to be discretised and gives an explicit representation of the interface which does not have to be coupled to an underlying grid. Furthermore, it naturally handles the discontinuity in normal stress. A review of boundary integral methods for Stokes flow is given in \cite{Pozrikidis2001}. For drops, boundary integral methods including surfactants in 3D are among others \cite{Bazhlekov2006numerical,Pozrikidis2018,Sorgentone2018}. In two dimensions, boundary integral methods have been used for e.g. vesicles as by \citeauthor{Marple2016} \cite{Marple2016} who simulated vesicle suspensions in confined flows and \citeauthor{Quaife2016}  \cite{Quaife2016} who focused on developing an adaptive time-stepping scheme using an implicit spectral deferred correction method. This was extended by \citeauthor{Bystricky} \cite{Bystricky} for rigid body suspensions. \citeauthor{kropinskilushi2011} \cite{kropinskilushi2011} simulated surfactant-covered bubbles with a boundary integral method and used a spectral method to compute the surfactant concentration. \citeauthor{Booty2013} \cite{Booty2013} extended this work to include a model for soluble surfactants in the no-diffusion limit in the bulk. However, the numerical solution of the boundary integral equations has its own challenges. Primarily, the near-singular behaviour of the discretised integrals when evaluating interfaces in close proximity. In 2D, \citeauthor{Ojala2015} \cite{Ojala2015} used a special quadrature scheme to ensure high accuracy also for close to touching clean drops. Other options for closely interacting drops includes quadrature by expansion (QBX) by \citeauthor{Klockner2013} \cite{Klockner2013}, which can be extended to 3D, and was done so efficiently for the special case of solid spheroids in \cite{KlintebergQBX}.

This paper proposes an efficient and accurate boundary integral method for the simulation of deforming bubbles and droplets in two dimensional Stokes flow; for both clean and surfactant-covered droplets. The method is an extension of that proposed by \citeauthor{kropinskilushi2011} \cite{kropinskilushi2011} and uses the special quadrature method in \cite{Ojala2015} to enable accurate simulations of close drop-drop interactions. To match the spectral accuracy in space, an adaptive time-integration scheme is employed and also described herein.

This paper also aims to validate the proposed numerical method using conformal mapping theory inspired by the approaches in \cite{crowdy2005pair,Tanveer1995}. For single surfactant-covered bubbles in an extensional flow at steady state, exact solutions are computed  by \citeauthor{siegelvalid} \cite{siegelvalid}. Furthermore, \citeauthor{crowdy2005pair} \cite{crowdy2005pair} used conformal mapping theory to semi-analytically compute the deformation of a pair of clean bubbles in an extensional flow. With this approach, a conformal mapping of the interface is considered and the time-dependent problem of deforming bubbles is rewritten into a system of ODEs for the conformal mapping parameters, whilst the flowfield can be described analytically. In this paper, this approach has been extended to surfactant-covered bubbles. Together, these solutions enable a precise and quantitative validation of the numerical method proposed by this paper. A comprehensive evaluation of the method is undertaken, examining temporal development, droplets in close proximity and the distribution of surfactants.

Additionally, the errors introduced when evaluating the near-singular integrals numerically using standard Gauss-Legendre quadrature are estimated. This work extends that of \citeauthor{AfKlinteberg2017} \cite{AfKlinteberg2017,klintebergadapt} where estimates based on contour integrals where derived for Laplace's and Helmholtz equations. The estimates are in excellent agreement with the observed numerical errors. Details of the estimates and their derivation can be found in \ref{sec:app_est}.

The paper is organised as follows: in \S\ref{sec:probform} the model, governing equations and boundary integral formulation are introduced. \S\ref{sec:methodBIE} explains the numerical solution of the boundary integral equation, the special quadrature method and provides results for the error estimates. The complete numerical method for the coupled problem is explained in \S\ref{sec:nummethod} and \S\ref{sec:valid} details the methods used for validation. In \S\ref{sec:results} numerical results for drop simulation and validation is presented.

%% file: srcfiles/formulation.tex

\section{Problem formulation}
\label{sec:probform}

\subsection{Governing equations}
\label{sec:goveq}
\begin{figure}[h!]
	\centering
	\includegraphics[width=0.5\textwidth]{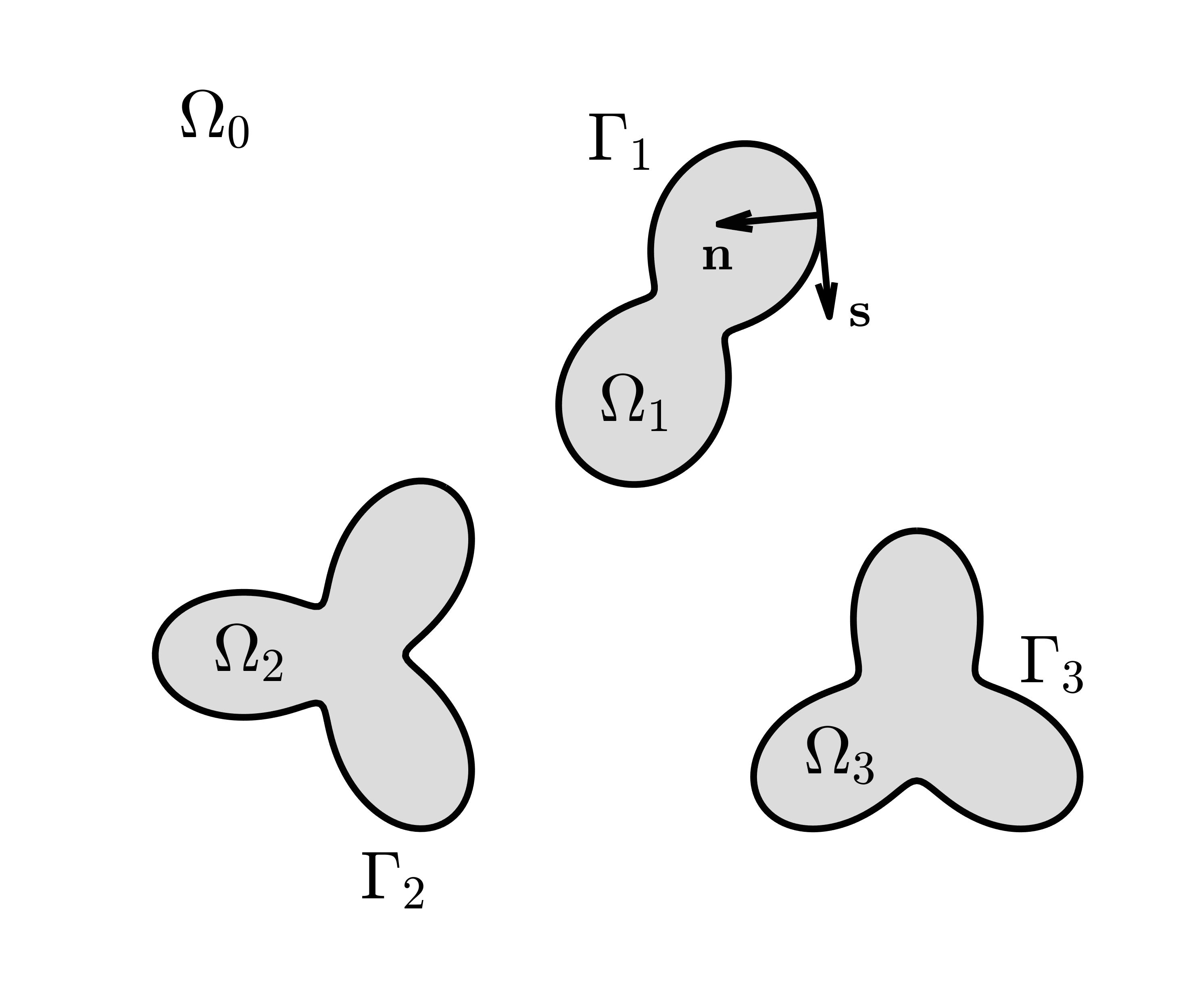}
	\caption{Example of domain consisting of three droplets.}
	\label{fig:form_domain}
\end{figure}

\noindent On the micro scale, the fluid velocity $\mathbf{u}$ and pressure $p$ of a flow can be computed by solving Stokes equations. Here, an infinite expanse of fluid, $\Omega_0$, with viscosity $\mu_0$ is considered. The fluid contains $n$ droplets, denoted by $\Omega_k$, for $k = 1, \hdots , n$, each with viscosity $\mu_k$. The interface between droplet $\Omega_k$ and bulk fluid $\Omega_0$ is denoted $\Gamma_k$, see Figure~\ref{fig:form_domain} for an example domain. The Stokes equations take the form
\begin{align}
	\begin{split}
		\mu_0 \Delta \mathbf{u}_0 = \nabla p_0,& \text{ for } \mathbf{x}\in\Omega_0,  \\
		\mu_k \Delta \mathbf{u}_k = \nabla p_k,& \text{ for } \mathbf{x}\in\Omega_k, \; k=1,\hdots,n.
	\end{split}
	\label{eq:form_stokeseq}
\end{align}
For droplets, $\mu_k \neq 0$ and the velocity is continuous over the interfaces. Also, the solution fulfils the normal stress balance over the interfaces,
\begin{align}
	-\left(p_0 - p_k\right)\mathbf{n} + 2\left(\mu_0\mathbf{e}_0 - \mu_k\mathbf{e}_k\right)\cdot\mathbf{n} = - \sigma\kappa\mathbf{n} + \nabla_s\sigma \text{ on } \Gamma,
	\label{eq:form_normalstress}
\end{align}
where $\Gamma = \bigcup_{k=1}^n \Gamma_k$. This condition states that the jump in normal stress over an interface is proportional to the curvature and Marangoni force. Furthermore, the kinematic interface condition states that the normal velocity of the interface, $\frac{d\mathbf{x}}{dt}\cdot\mathbf{n}$, is equal to the normal fluid velocity,
\begin{align}
	\dfrac{d\mathbf{x}}{dt}\cdot\mathbf{n} = \mathbf{u}\cdot\mathbf{n}.
	\label{eq:form_kincond}
\end{align}
Here, $\mathbf{n}$ is the inward facing normal, $\sigma$ surface tension coefficient, $\kappa = \nabla_s\cdot\mathbf{n}$ curvature and $\mathbf{e}$ the rate of strain tensor. Note that this curvature is negative for a circular drop. The term $\nabla_s\,\sigma$, where $\nabla_s$ is the surface gradient operator, gives the tangential stress (Marangoni force) which is the result of a non-uniform surface tension. The interfaces are discretised clockwise by $s\in[0,L_k(t)]$, where $L_k(t)$ is the length of the interface $\Gamma_k$ at time $t$.

The limit when $\mu_k = 0$ corresponds to the study of inviscid bubbles. The continuity condition of velocity across the interface can then be disregarded, as there is no velocity inside the bubbles. The kinematic condition \eqref{eq:form_kincond} is still necessary, as is the normal stress balance \eqref{eq:form_normalstress} modified to contain only pressure and strain tensor of the bulk fluid.

As the fluid domain $\Omega_0$ is unbounded, an additional condition at infinity is needed,
\begin{align*}
	\mathbf{u} \rightarrow \mathbf{u}_\infty, \text{ as } |\mathbf{x}|\rightarrow \infty.
\end{align*}
Here it is typical to impose a linear far-field flow, where
\begin{align}
	\mathbf{u}_\infty = \begin{pmatrix} Q_0 & B_0+\frac{G_0}{2} \\
										B_0-\frac{G_0}{2} & -Q_0
						\end{pmatrix} \cdot \mathbf{x}.
	\label{eq:form_farfield}
\end{align}
The two imposed linear flows considered in this paper are: \textit{extensional flow}, where $B_0=G_0=0$ and \textit{shear flow}, where $Q_0=0$ and $G_0=2B$.

The surfactants considered here are insoluble, i.e. they exist only on the interfaces of the droplets. Their concentration, $\rho(s,t)$, is described by the convection-diffusion equation \cite{stone}
\begin{align}
	\frac{D\rho}{Dt} + \rho\left(\nabla_s \cdot\mathbf{u}\right) = \frac{1}{D_\Gamma}\nabla_s^2\, \rho, \text{ for } \mathbf{x}\in\Gamma.
	\label{eq:form_surfacteq}
\end{align}
Here, $D_\Gamma$ is the diffusion coefficient along the interface and $\frac{D}{Dt}$ the material derivative. As the surfactants are insoluble, their mass on each interface is conserved,
\begin{align}
	\frac{d}{dt}\int_{\Gamma_k(t)} \rho(t)dS = 0, \quad k=1,\hdots,n.
	\label{eq:form_conservmass}
\end{align}

The surfactant concentration and the surface tension at an interface are coupled through an equation of state. Two common equations of states are the Langmuir equation of state,
\begin{align}
	\sigma(\rho) = \sigma_0 + RT\log\left(\rho_\infty-\rho\right)
	\label{eq:form_eqstate_langmuir}
\end{align}
and  the linearised equation of state
\begin{align}
	\sigma(\rho) = \sigma_0 - RT\rho,
	\label{eq:form_eqstate_linear}
\end{align}
as described in \cite{Pawar1996}. Here,  $R$ is the universal gas constant, $T$ the temperature and $\sigma_0$ the surface tension coefficient of a clean interface. Furthermore, $\rho_\infty$ is the maximum monolayer packing concentration of surfactants on $\Gamma(t)$. Generally, the Langmuir equation of state is considered more accurate. The linear equation of state is mostly used for problems with low saturation levels. It is also commonly used in many validation cases, which is the reason it is included here. To switch between the two is trivial.

\subsection{Nondimensionalisation}
\noindent The above equations are nondimensionalised using a characteristic length $r_0$, defined as the radius of an initial droplet. The characteristic velocity is defined as the ratio $\frac{\sigma_0}{\mu_0}$, where $\mu_0$ is the viscosity of the bulk and $\sigma_0$ the surface tension coefficient of a clean interface. Furthermore the characteristic pressure is defined as $\frac{\sigma_0}{r_0}$. For the remainder of this paper, all physical quantities are considered in their non-dimensional form.

Rewriting the equations, \eqref{eq:form_stokeseq} becomes
\begin{align}
	\begin{split}
		\Delta\mathbf{u}_0 = \nabla p_0,& \text{ for } \mathbf{x}\in\Omega_0, \\
		\lambda_k \Delta\mathbf{u}_k = \nabla p_k,& \text{ for } \mathbf{x}\in\Omega_k, \; k=1,\,\hdots,\,n,
	\end{split}
	\label{eq:form_stokesnondim}
\end{align}
and \eqref{eq:form_normalstress}
\begin{align}
	-\left(p_0-p_k\right)\mathbf{n} + 2\left(\mathbf{e}_0-\lambda_k\mathbf{e}_k\right) \cdot \mathbf{n} = -\sigma\kappa\mathbf{n} + \nabla_s\sigma, \text{ on } \Gamma.
	\label{eq:form_normalstressnondim}
\end{align}
Here, $\lambda_k = \frac{\mu_k}{\mu_0}$ is the viscosity ratio between droplet $k$ and bulk fluid. Also the imposed far-field flow is nondimensionalised, such that
\[
	\left(Q,\;G,\,B \right)~=~\frac{r_0\mu_0}{\sigma_0}\left(Q_0,\;G_0,\;B_0\right).
\]
The nondimensionalisation of the surfactant concentration depends on which equation of state is used:
\paragraph{Case 1} In the case of the Langmuir equation of state, it is natural to nondimensionalise $\rho$ with the maximum monolayer packing concentration $\rho_\infty$, and \eqref{eq:form_eqstate_langmuir} reads
\begin{align}
	\sigma = 1 + E\log(1-\rho),
	\label{eq:form_eqstate_langmuir_nd}
\end{align}
where $E=\frac{RT\rho_\infty}{\sigma_0}$ is the so-called elasticity number.

\paragraph{Case 2} When it comes to the linearised equation of state, it is instead common to nondimensionalise with the initial surfactant concentration, $\rho_0$. Then, \eqref{eq:form_eqstate_linear} becomes
\begin{align}
	\sigma = 1 - E\rho,
	\label{eq:form_eqstate_linear_nd}
\end{align}
with $E=\frac{RT\rho_0}{\sigma_0}$.
\\ \\
In both cases, \eqref{eq:form_surfacteq} becomes
\begin{align}
	\frac{D\rho}{Dt} + \rho\left(\nabla_s \cdot\mathbf{u}\right) = \frac{1}{Pe_\Gamma}\nabla_s^2\, \rho, \text{ for } \mathbf{x}\in\Gamma.
	\label{eq:form_surfacteq_nondim}
\end{align}
where $Pe_\Gamma=\frac{\sigma_0 r_0}{\mu_0 D_\Gamma}$ is the Peclet number.
\\ \\
In summary, to evolve the deforming interfaces of surfactant-covered droplets in time, the following coupled system needs to be solved for $\mathbf{x}\in\Gamma(t)$:
\begin{align}
	\begin{split}
		\frac{d\mathbf{x}}{dt} &= \mathbf{u}(\mathbf{x},\rho,t), \\
		\frac{D\rho}{Dt} &= -\rho\left(\nabla_s \cdot\mathbf{u}\right) + \frac{1}{Pe_\Gamma}\nabla_s^2\, \rho.
	\end{split}
	\label{eq:form_system}
\end{align}
The solution to this system will be explained in steps. First, the boundary integral equation (BIE) formulation to compute $\mathbf{u}$ is described below in \S\ref{sec:bieform}.  In \S\ref{sec:methodBIE} a numerical method to solve the BIE accurately is explained. Then, \S\ref{sec:nummethod} explains the numerical solution of the surfactant concentration, the coupling of the system and its evolution in time.

%% file: srcfiles/bie.tex

\subsection{Boundary integral formulation}
\label{sec:bieform}
\noindent In two dimensions, it is convenient to regard this problem in the complex plane. The spatial variable $\mathbf{x}$ then corresponds to the complex variable $z=x+iy$. On the interface, $z$ is regarded as a function of the parameter $\alpha$ and time, $z(\alpha,t)$, where $\alpha\in[0,2\pi]$. The two parametrisations $s$ and $\alpha$ are linked through $\alpha = 2\pi s/L_k(t)$ for each interface $\Gamma_k(t)$, where $s$ was introduced below \eqref{eq:form_kincond}.

The integral equation used to compute the velocity on each interface is obtained through the Sherman-Lauricella formulation as described in \cite{kropinski2001}. This approach stems from the fact that Stokes equations in two dimensions can be reduced to the biharmonic equation and that all physical quantities of the problem, such as velocity, pressure and vorticity, can be expressed as combinations of a pair of Goursat functions \cite{Langlois2014}. If the Goursat functions are constructed in a particular way, a Fredholm equation of the second kind is obtained to solve for a complex-valued density $\mu(z)$ defined on the interfaces $\Gamma(t)$,
\begin{align}
\begin{split}
	\mu(z) &+ \beta(z)\mathcal{T}(z,\mu) + \beta(z)\int_{\Gamma(t)}\mu(\tau)|d\tau| \\
	=& -\sigma(\alpha,t)\frac{\gamma(z)}{2}\frac{\partial z}{\partial \alpha} - \beta(z)\left(B-iQ\right)\overline{z}, \quad z\in\Gamma(t),
\end{split}
\label{eq:bie_mu}
\end{align}
where $\mathcal{T}(z,\mu)$ is the complex variable formulation of the stresslet
\begin{align*}
 	\mathcal{T}(z,\mu) = \frac{1}{\pi}\int_{\Gamma(t)} \mu(\tau)\Im\left(\frac{d\tau}{\tau-z}\right) +\frac{1}{\pi}\int_{\Gamma(t)} \overline{\mu(\tau)} \frac{\Im\left(d\tau(\overline{\tau}-\overline{z})\right)}{\left(\overline{\tau}-\overline{z}\right)^2},
\end{align*}
see \cite{kropinski2001} for details. The third term in \eqref{eq:bie_mu} evaluates to zero as a result of the area conservation of the droplets, and is used to remove rank deficiency in the case of inviscid bubbles. Here, notation as in \cite{Ojala2015} is used, where $\beta(z) \coloneqq \frac{1-\lambda_k}{1+\lambda_k}$ and $\gamma(z) \coloneqq \frac{1}{1+\lambda_k}$ for $z\in\Gamma_k(t)$. In the last term on the right hand side, $B$ and $Q$ are part of the imposed far-field flow in \eqref{eq:form_farfield}. Furthermore, $\sigma(\alpha,t)~\eqqcolon~\sigma(\rho(\alpha,t))$ is the surface tension coefficient obtained from the surfactant concentration through either equation of state, \eqref{eq:form_eqstate_langmuir_nd} or \eqref{eq:form_eqstate_linear_nd}. Note that the expression in \eqref{eq:form_normalstressnondim} has been integrated once, thus no differentiation of $\sigma(\alpha,t)$ is needed.

Once $\mu(z)$ is obtained for all $z\in\Gamma(t)$, the velocity can be evaluated for $z\in\Omega_0\cup\Gamma\cup\Omega_k$, as
\begin{align}
\begin{split}
	u(z) &= u_1(z) + iu_2(z) = -\frac{1}{\pi}\int_{\Gamma(t)}\mu(\tau)\Re\left(\frac{d\tau}{\tau-z}\right) \\ 
	-& \frac{1}{\pi i}\int_{\Gamma(t)}\overline{\mu(\tau)}\frac{\Im\left(d\tau(\overline{\tau}-\overline{z})\right)}{\left(\overline{\tau}-\overline{z}\right)^2} + (Q+iB)\overline{z} - \frac{iG}{2}z,
\end{split}
\label{eq:bie_u}
\end{align}
where the last two terms on the right hand side represent the far-field velocity $\mathbf{u}_\infty$ in complex form. Note that the first integral in \eqref{eq:bie_u} is singular and interpreted in a principal value sense.

%% file: srcfiles/meth_bie.tex

\section{Solving the BIE numerically}
\label{sec:methodBIE}

\noindent Here, a summary of the numerical method to solve the BIE is given. More details are available in \cite{Ojala2015}. The solution procedure described is valid for any instance of time $\bar{t}$. 

To compute the integrals \eqref{eq:bie_mu} and \eqref{eq:bie_u} accurately, a high-order discretisation of $\Gamma(\bar{t})$ is needed. Here, an explicit representation of the interfaces is used; $z_k(\alpha)$, for droplets $k=1,\hdots,n$, where $\alpha\in[0,2\pi)$. On interface $k$, each interface is split into panels and $N_k$ discretisation points are placed on a composite 16-point Gauss-Legendre grid. The total number of discretisation points becomes $N=\sum_{k=1}^n N_k$. 

\subsection{Solve for $\mu(z)$}
\noindent To compute $\mu_i \approx \mu(z_i)$, where $z_i$ are the Gauss-Legendre points on the interfaces at time $\bar{t}$, a Nystr\"{o}m method is used. The discretised version of \eqref{eq:bie_mu} is
\begin{align}
	\begin{split}
	\mu_i +& \dfrac{\beta_i}{\pi} \sum_{j=1}^N \mu_j M_{ij}^{(1)} + \dfrac{\beta_i}{\pi}\sum_{j=1}^N \overline{\mu_j} M_{ij}^{(2)} + \beta_i \sum_{j=1}^N \mu_j w_j|z^\prime_j| \\
	&= -\sigma_i \dfrac{\gamma_i}{2} \dfrac{z^\prime_i}{|z^\prime_i|} - \beta_i(B-iQ)\overline{z_i}, \; i=1,\hdots,N,
	\end{split}
	\label{eq:numm_mudisc}
\end{align}
where $z^\prime_i = \left.\frac{\partial z}{\partial\alpha}\right\vert_i$. Here, $w_j$ are the Gauss-Legendre weights associated with $z_j$. Also, $\sigma_i~=~\sigma(\alpha_i)$, $\beta_i = \beta(z_i)$ and $\gamma_i = \gamma(z_i)$. Moreover,
\begin{align*}
		M_{ij}^{(1)} = w_j \Im\left\{ \dfrac{z^\prime_j}{z_j-z_i}\right\}, \; \;
		M_{ij}^{(2)} = w_j \dfrac{\Im\left\{ z^\prime_j(\overline{z_j}-\overline{z_i}) \right\}}{(\overline{z_j}-\overline{z_i})^2},
\end{align*}
for all $j\neq i$. The limits when $j=i$ are finite and given by 
\begin{align*}
	M_{ii}^{(1)} = w_i \Im\left\{ \dfrac{z_i^{\prime\prime}}{2z^\prime_i} \right\}, \; \; M_{ii}^{(2)} = w_i \dfrac{\Im\left\{z_i^{\prime\prime}\overline{z_i^\prime} \right\}}{2(\overline{z_i^\prime})^2}.
\end{align*}
Solving \eqref{eq:numm_mudisc} to obtain $\mu_i$, $i=1,\hdots,N$, boils down to solving a system $A\mathbf{\mu} = \mathbf{b}$, where $\mu = (\mu_1 \hdots \mu_N)^T$, $\mathbf{b}$ represents the right hand side in \eqref{eq:numm_mudisc} and $A\mu$ the expression in the left hand side. The matrix $A$ is dense, but as \eqref{eq:bie_mu} is a Fredholm equation of the second kind, the condition number of $A$ is typically small and does not increase with refinement. Thus, there is no need for a preconditioner when solving \eqref{eq:numm_mudisc} with an iterative algorithm such as \texttt{GMRES}. Furthermore, the matrix-vector multiplications are sped up using the fast multipole method, \texttt{FMM} \cite{Greengard1997}.

\subsection{Compute $u(z)$}
\noindent Once $\mu_i$ is obtained for all discretisation points $z_i$, the velocity $u_i \approx u(z_i)$ can be computed through \eqref{eq:bie_u}. To handle the singular integral, singularity subtraction is used; when evaluating $u_i$ the integral 
\[
\dfrac{1}{\pi}\int_\Gamma \mu_i\,\Re\left\{ \dfrac{d\tau}{\tau-z_i} \right\}
\]
is added and subtracted from the expression in \eqref{eq:bie_u}. Using calculus of residues this removes the singularity and leads to the discretised expression
\begin{align}
\begin{split}
	u_i = -\dfrac{w_i}{\pi}\mu_i^\prime &- \dfrac{1}{\pi}\sum_{j=1, i\neq j}^N (\mu_j-\mu_i) w_j \Re\left\{ \dfrac{z_j^\prime}{z_j-z_i} \right\} \\
	&-\dfrac{1}{i\pi} \sum_{j=1}^N \overline{\mu_j}M_{ij}^{(2)} + (Q+iB)\overline{z_i} -\dfrac{iG}{2}z_i,
\end{split}	
\label{eq:numm_udisc}
\end{align} 
where $M_{ij}^2$ is defined previously.

%% file: srcfiles/meth_est.tex

\subsection{Error estimates}
\label{sec:est}
\noindent Regarding the equations \eqref{eq:bie_mu} and \eqref{eq:bie_u}, they both contain terms $(z_j-z_i)^{-1}$. For $M^{(1)}$ and $M^{(2)}$ the limits when $i=j$ are well-defined, and for the term $\Re\left\{ \frac{z_j^\prime}{z_j-z_i}\right\}$ singularity subtraction is used, see \eqref{eq:numm_udisc}. However, when droplets get close to each other, i.e. when $\|z_i-z_j\|\ll 1$ for points $z_i$ and $z_j$ on different interfaces, the integrals are said to become near-singular. Analytically this is not an issue, but numerically large errors are introduced as the integrand gets harder to resolve.

\begin{figure}[h!]
	\subfigure[Domain $\Omega$, source points and solution $u(z)$.]{%
	\centering
	\includegraphics[width=0.49\textwidth]{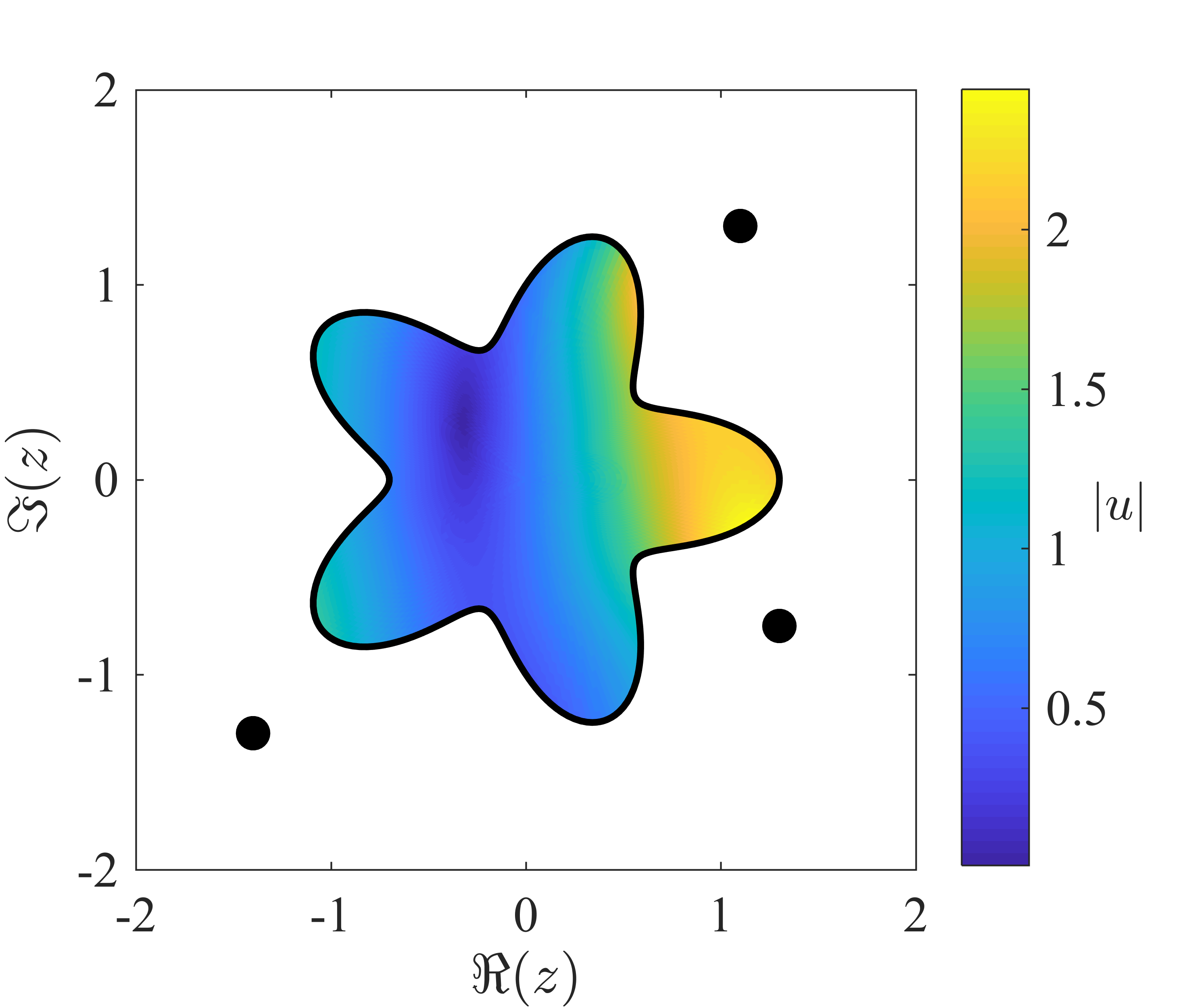}
		\label{fig:est_domu}}%
	\hfill
	\subfigure[Error in first quadrant of domain.]{%
	\includegraphics[width=0.49\textwidth]{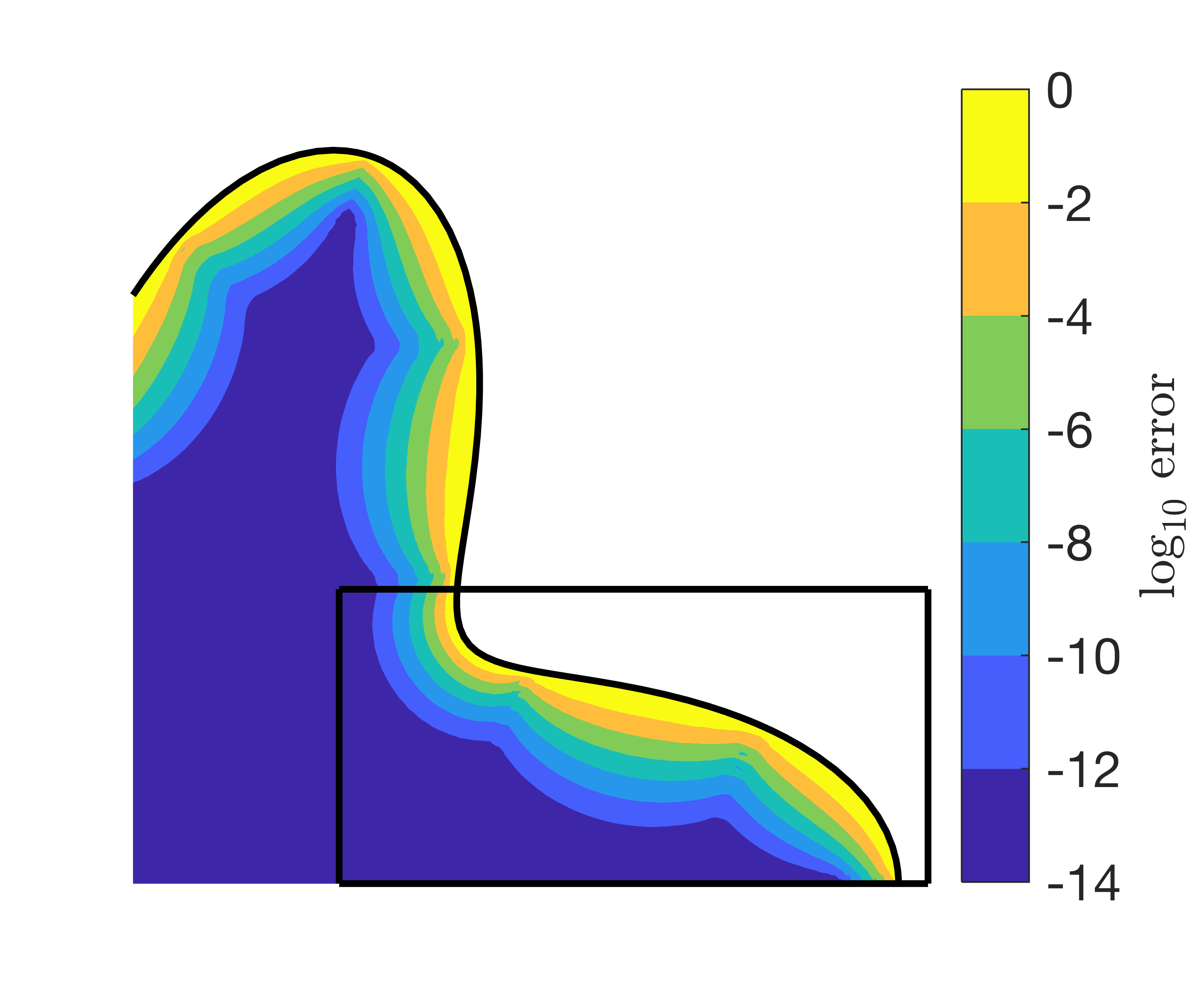}\label{fig:est_err}}%
	\caption{Example of analytical solution $u(z)$ and error to \eqref{eq:est_u} given sources ($x_1$,$f_1$), ($x_2,f_2$) and ($x_3$,$f_3$) (black markers) in a domain $\Omega$. Solution is computed with a composite 16-point Gauss-Legendre quadrature with 25 panels.}
\end{figure}

To demonstrate how these errors grow, Stokes equations are solved in a fixed domain with Dirichlet boundary conditions. This problem contains the same near-singular behaviour as \eqref{eq:bie_u}. Given any analytical solution to Stokes equations (for example, that generated by a point source forcing as shown in Figure~\ref{fig:est_domu}), the Dirichlet velocity data can be evaluated at the boundary and used as boundary conditions for the numerical method. The boundary integral formulation and exact problem setting is described in \ref{sec:app_est}, the domain and solution $u(z)$ are shown in Figure~\ref{fig:est_domu}.

\begin{figure}[h!]
	\centering
	\subfigure[$25$ panels]{%
	\includegraphics[width=0.49\textwidth,trim={1cm 2cm 0 2cm},clip]{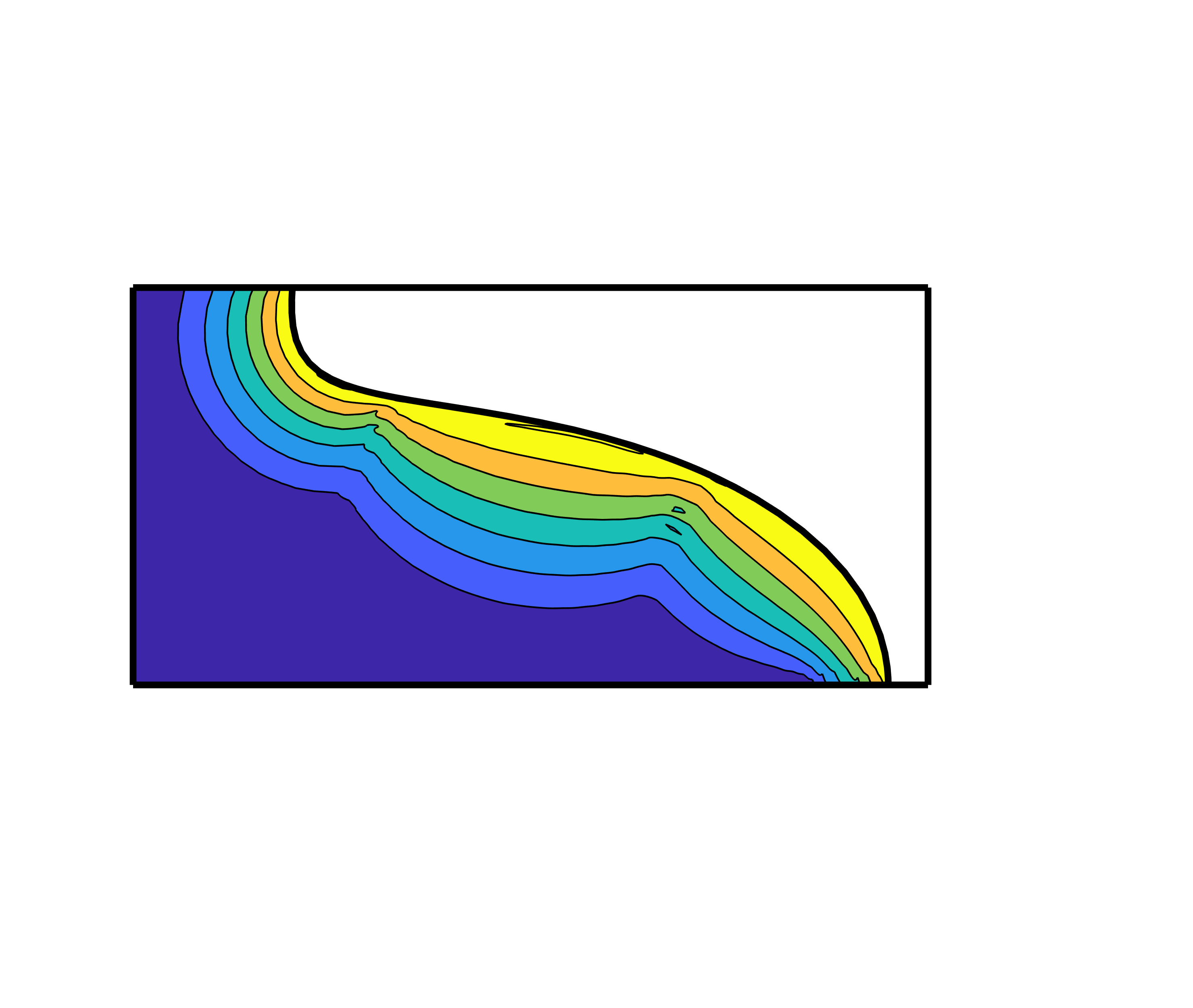}}%
	\;
	\hspace{-1cm}
	\subfigure[$50$ panels]{%
	\includegraphics[width=0.49\textwidth,trim={1cm 2cm 0 2cm},clip]{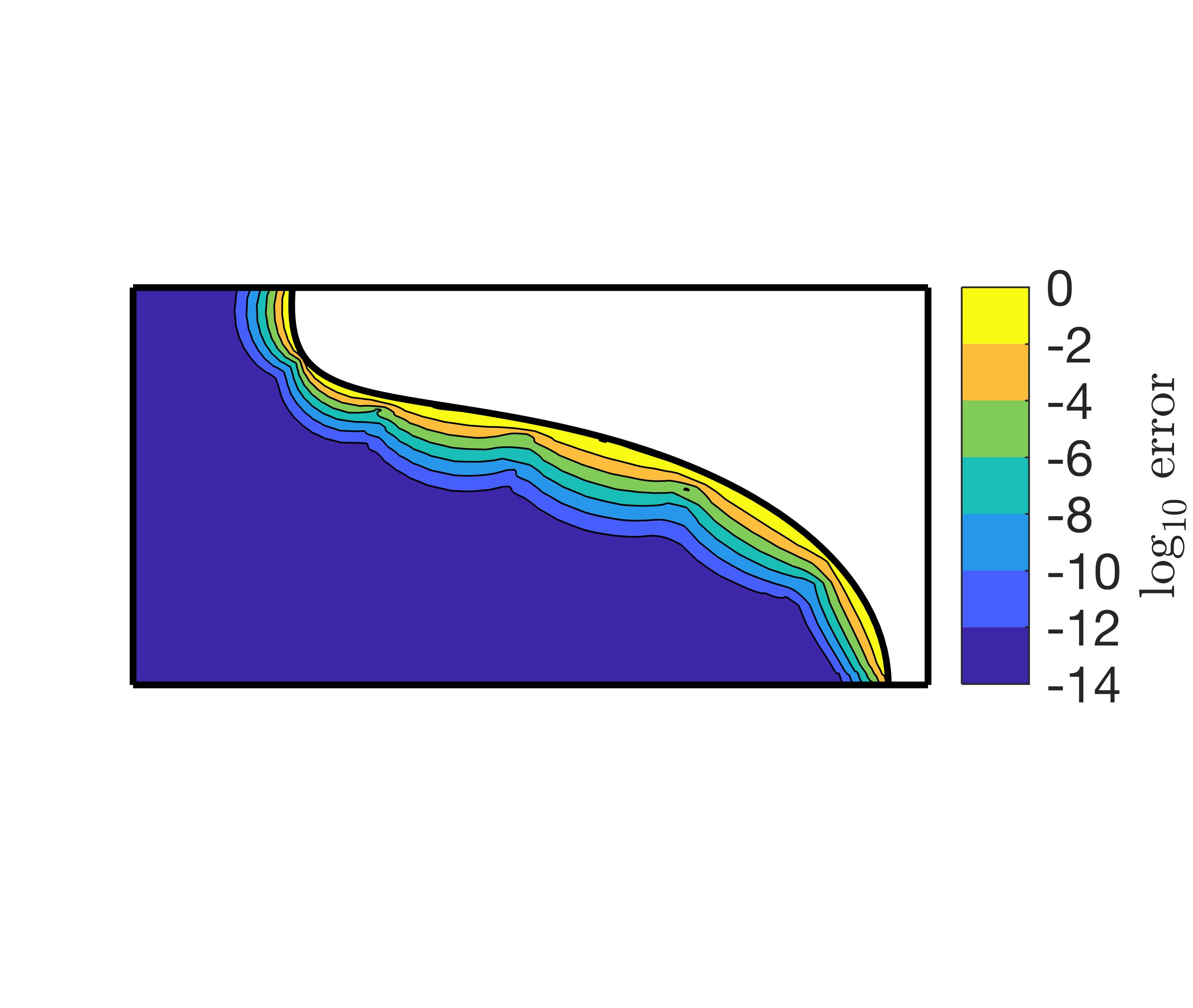}}%
	\caption{Computed error estimates in black for error levels $10^{-p}$, where $p=14,12,\hdots,2,0$. Measured errors in colour.}
	\label{fig:est_est}
\end{figure}

The relative error of $u(z)$ compared to the exact solution, using $25$ Gauss-Legendre panels, is shown in Figure~\ref{fig:est_err}. Here only the first quadrant is shown, as the error behaves identical in the other three. It is clear that at an evaluation point $z_0$ close to the boundary, errors are large and other treatments of the integrals are needed for these cases. In this paper a special quadrature is used, which is explained in \S\ref{sec:specq}. How close to the panel the errors become large depends on the refinement of the discretisation.  In order to know when special treatment is needed, the quadrature errors of such near-singular integrals can be estimated. This was originally done by \citeauthor{AfKlinteberg2017}, where estimates for the quadrature errors for Laplace's and Helmholtz equations were derived in \cite{AfKlinteberg2017,klintebergadapt}. Using the same approach based on contour integration and calculus of residues for integrals of the type appearing in \eqref{eq:bie_u}, the errors when computing $u(z)$ can be estimated also for Stokes equations. Details of this can be found in \ref{sec:app_est}. The estimates follow the error levels remarkably well, as is shown in Figure~\ref{fig:est_est}. Both errors and estimates for $25$ and $50$ panels are shown, and it is clear that refining the interface discretisation makes the region of large errors narrower, but will not eliminate it.

%% file: srcfiles/meth_spec.tex

\subsection{Special quadrature}
\label{sec:specq}
\noindent To improve the accuracy of the computation of $\mu(z,\bar{t})$ and $u(z,\bar{t})$ near the interfaces, the special quadrature method of \cite{Ojala2015} will be employed. It was originally introduced by \citeauthor{ojalahelsing} \cite{ojalahelsing} for Laplace's equation and extended to Stokes equations in \cite{Ojala2015}. This is a local method which regards point-panel pairs. For each evaluation point the error when using a standard Gauss-Legendre quadrature is estimated on each panel. If it is too large over a specific panel, $\Gamma_i$, the integral over that panel will instead be treated semi-analytically. An overview of the special quadrature can be found in \ref{sec:app_spec}.

\begin{figure}[h!]
\centering
\subfigure[Error in whole domain.]{%
	\includegraphics[width=0.49\textwidth,trim={1cm 1cm 2cm 1cm},clip]{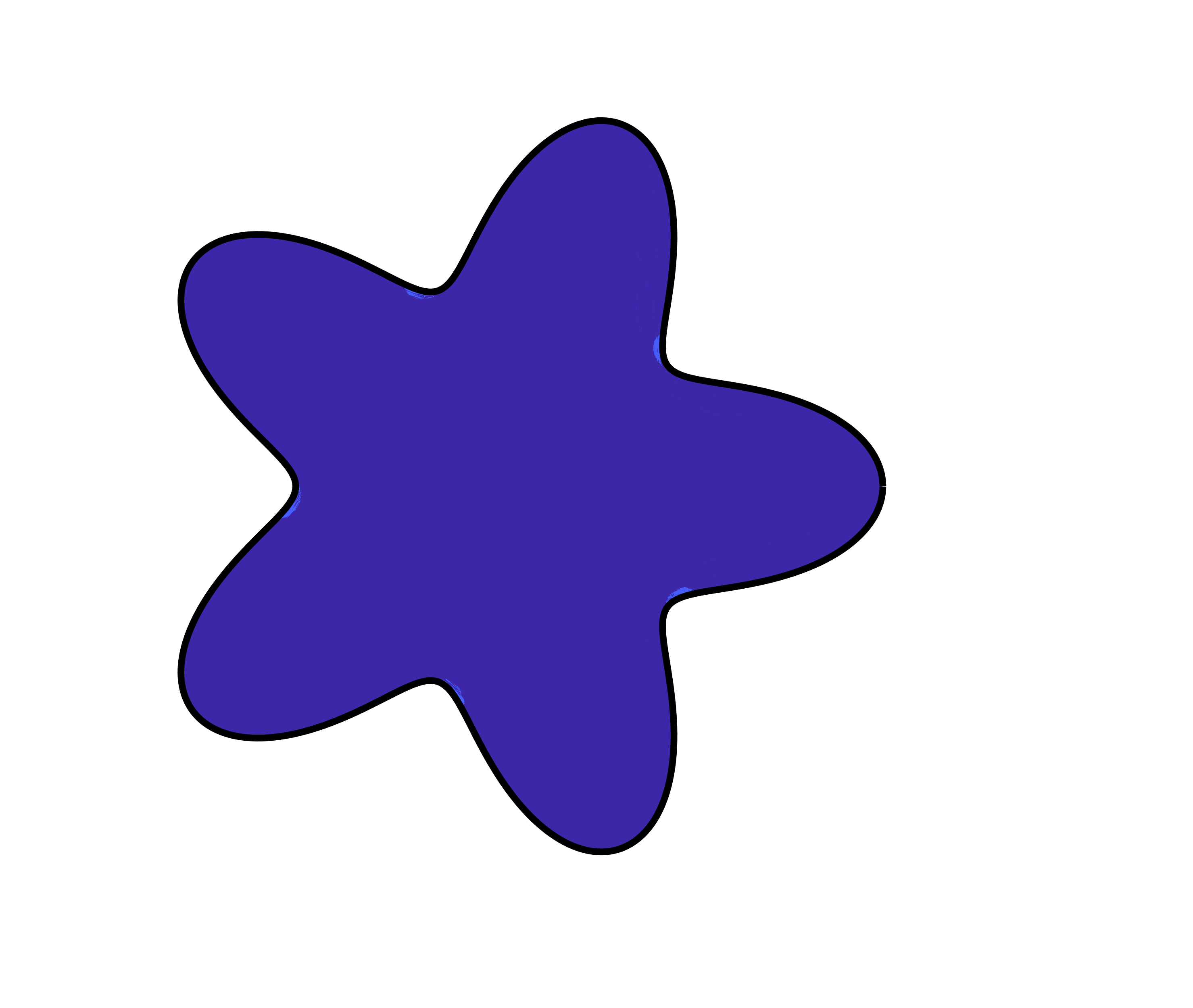}}%
	\;
\subfigure[Error in first quadrant.]{%
	\includegraphics[width=0.49\textwidth]{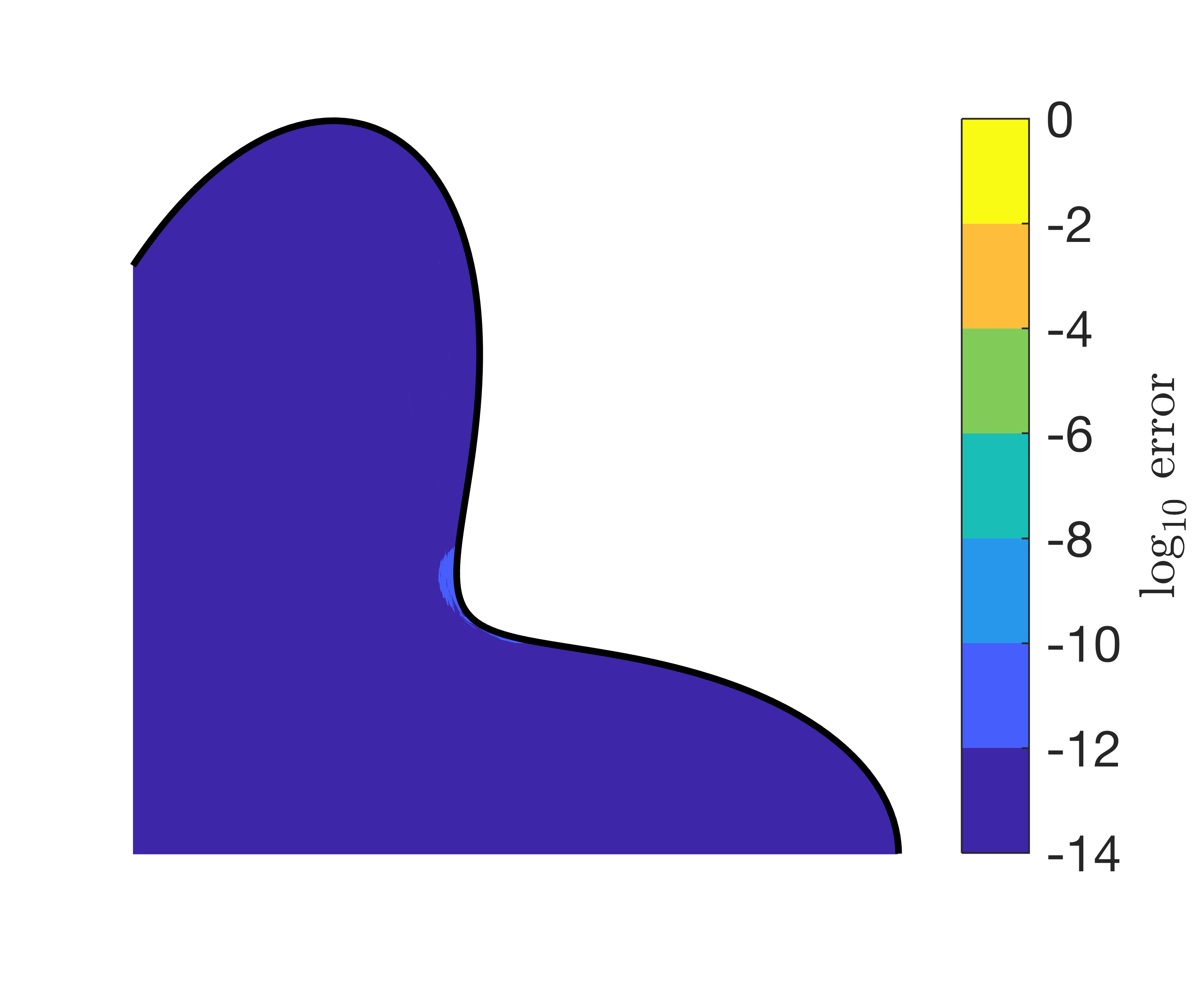}}%
	\caption{Logarithm of error when solving Stokes equations using special quadrature for evaluating at points close to the boundary, cf. Figure~\ref{fig:est_err}. Solution computed with $50$ Gauss-Legendre panels.}
\label{fig:numm_specquad}
\end{figure}

Solving the same problem of Stokes equations in a domain with non-deforming boundaries as in \S\ref{sec:est}, $u(z)$ obtained by standard quadrature is corrected in the regions where the error is large. The reduction of error when using the special quadrature can be seen in Figure~\ref{fig:numm_specquad}. This error should be compared to that in Figure~\ref{fig:est_err}, where standard 16-point composite Gauss-Legendre quadrature is used. It is clear that errors can be kept at a very low level of order $10^{-10}$ or less also for evaluation points close to the interface.

%% file: srcfiles/nummethods.tex

\section{Numerical method in time}
\label{sec:nummethod}
\noindent Several steps are needed to simulate the deformation of surfactant-covered droplets through the system in \eqref{eq:form_system}. Here, a method of lines approach is used, where the discretisation in space generates a system of ODEs to solve in time. This system can then be solved using a numerical method for initial value problems using a time-integration scheme as described in \S\ref{sec:methods_time}.

To evolve the interfaces in time, as well as to compute the surfactant concentration, an equidistant spatial discretisation is used. However, to compute the integrals in the BIE formulation described in \S\ref{sec:methodBIE}, it is beneficial to use a composite Gauss-Legendre discretisation. Thus, a hybrid method using both discretisations is employed. The interface of each droplet is discretised with $N_k$ points. For the uniform discretisation, these points are placed equispaced in arc length. The quadrature weights associated with this discretisation are the standard trapezoidal weights. For the composite Gauss-Legendre discretisation, the interface is divided into $N_k/16$ panels of equal length, each then discretised with a 16-point Gauss-Legendre quadrature rule. To go between the two discretisations, the points on each panel are interpolated.

Here follows a short overview of the steps needed at each instance in time, all of which will be explained in more detail further on in \S\ref{sec:methods_GL}-\ref{sec:methods_surf}. In order to evolve the droplets in time, the interface velocitites are obtained through the following steps:
\begin{enumerate}
	\item Upsample the uniform discretisation and interpolate to the composite 16-point Gauss-Legendre grid.
	\item Compute the fluid velocity at the boundaries accurately by solving the BIE and applying special quadrature where needed, see \S\ref{sec:methodBIE}.
	\item Interpolate the velocity back to the uniform grid and downsample.
	\item Compute the appropriate tangential velocity to keep the discretisation points equidistant in arc length.
\end{enumerate}
The velocity is then fed into a time-integration scheme for propagation of the droplet interfaces in time. Each time the fluid velocity is computed, the surface tension coefficient is needed. Since this is dependent on the surfactant concentration through \eqref{eq:form_eqstate_langmuir_nd} or \eqref{eq:form_eqstate_linear_nd}, also the surfactant concentration needs to be propagated in time through \eqref{eq:form_surfacteq_nondim}, see \S\ref{sec:methods_surf}.

\subsection{Interpolation to and from the composite Gauss-Legendre grid}
\label{sec:methods_GL}
\noindent First, upsampling is needed for stability and is achieved by zero-padding the Fourier coefficients of the equidistant discretisation with a factor two. Then, the interpolation to the composite Gauss-Legendre grid is done via a non-uniform \texttt{FFT} (\texttt{nuFFT}). Once the velocity and special quadrature corrections have been computed, a 16-point Gauss-Legendre interpolation on each panel is used to go back to the uniform discretisation. Finally, the velocity vector is downsampled and a Krasny filter is applied to the velocity vector, where all Fourier modes below the absolute value of $10^{-12}$ are set to zero. See \cite{Ojala2015} for more details.

\subsection{Modifying the tangential velocity}
\label{sec:methods_modT}
\noindent With the velocity at the interfaces represented on the equidistant discretisation, denoted $u$, a new modified velocity for the discretisation points is computed, $\frac{d\mathbf{x}}{dt}$. The reason for this is to avoid clustering of discretisation points as the interface evolves in time and it was originally introduced in \cite{kropinski2001}. The idea is to decompose the velocity $u\left(z(s)\right)$ into its normal and tangential components,
\begin{align*}
	u\left(z(s)\right) = \left[u_n(s) + iu_t(s) \right]n(s),
\end{align*}
where $n\left(s\right)$ is the complex representation of the normal, $s\in[0,L_k(t)]$ and $u_n$ and $u_s$ are real valued. The modified velocity is then defined as
\begin{align}
	\frac{d\mathbf{x}}{dt}(s) = \left[u_n(s)+i\tilde{u}_t(s)\right] n(s).
	\label{eq:numm_tanv}
\end{align}
The new tangential velocity $\tilde{u}_t$ can be chosen to dynamically preserve $s$ equidistant in arc length. It can be derived as (see reference above)
\begin{align*}
	\tilde{u}_t(s) = \dfrac{s}{2\pi}\int_0^{2\pi} \Im\left\{ \dfrac{z^{\prime\prime}(q)}{z^{\prime}(q)}\right\}u_n(q)dq - \int_0^s \Im\left\{ \dfrac{z^{\prime\prime}(q)}{z^{\prime}(q)}\right\}u_n(q)dq,
\end{align*}
which actually boils down to an antiderivative taken with an FFT. The notation of this tangential velocity is adapted from \cite{Ojala2015}. Note that the modified velocity $\frac{d\mathbf{x}}{dt}$ still fulfils the kinematic condition \eqref{eq:form_kincond}.

\subsection{Spatial adaptivity}
\label{sec:methods_spatialadap}
\noindent Using the modified tangential velocity as explained above keeps the interface discretisation points equidistant in arc length over time, i.e. with equal spacing $\Delta s$. As the droplets deform, the spacing $\Delta s$ might increase or decrease. To keep the spatial resolution constant, a spatial adaptivity is implemented where more discretisation points are added to the interfaces if $\Delta s$ grows more than a set threshold. Likewise, in case when $\Delta s$ becomes too small, discretisation points are removed. These operations are both performed on the equidistant discretisation.

\subsection{Solving the surfactant equation}
\label{sec:methods_surf}
\noindent The material derivative in \eqref{eq:form_surfacteq_nondim} can be expanded utilising the interface velocity $u\left(z(s)\right)$ and the velocity of the discretisation points $\frac{d\mathbf{x}}{dt}(s)$, $s~\in~[0,L_k(t)])$. The convection-diffusion equation for the surfactant concentration becomes
\begin{align}{}
	\dfrac{\partial \rho}{\partial t} + \left(\mathbf{u} - \dfrac{d\mathbf{x}}{dt}\right) \cdot \nabla_s \rho + \rho\nabla_s\cdot \mathbf{u} = \dfrac{1}{Pe_\Gamma} \nabla_s^2\rho, \; \mathbf{x} \in \Gamma.
	\label{eq:numm_surf}
\end{align}
Using that the normal components of $\frac{d\mathbf{x}}{dt}$ and $u$ are the same and that the surface gradient operator can be written as $\nabla_s(\cdot)= \mathbf{s}\frac{\partial}{\partial s}(\cdot)$, \eqref{eq:numm_surf} can for each droplet $k$ be rewritten as
\begin{align}
	\dfrac{\partial\rho}{\partial t} = \tilde{u}_t\dfrac{\partial\rho}{\partial s} - \dfrac{\partial \left(\rho u_t\right)}{\partial s} - \rho u_n \kappa + \dfrac{1}{Pe_\Gamma}\dfrac{\partial^2\rho}{\partial s^2}, \; s\in[0,L_k(t)].
	\label{eq:numm_surfs}
\end{align}
Here, the first two terms in the right hand side correspond to advection, the third stretching and the last term represents diffusion along the interface. To reformulate \eqref{eq:numm_surfs} in terms of $\alpha\in[0,2\pi)$, the equal-arc length parameter $s^k_\alpha(t)$ is defined as
\begin{align}
	s_\alpha^k(\alpha,t) = \frac{1}{2\pi}L_k(t),
\end{align}
see \cite{kropinskilushi2011}. The convection-diffusion equation becomes
\begin{align}
	\dfrac{\partial\rho}{\partial t} = \dfrac{\tilde{u}_t}{s^k_\alpha(t)}\dfrac{\partial\rho}{\partial \alpha} - \dfrac{1}{s^k_\alpha(t)}\dfrac{\partial \left(\rho u_t\right)}{\partial \alpha} - \rho u_n \kappa + \dfrac{1}{Pe_\Gamma \,s_\alpha^k(t)^2}\dfrac{\partial^2\rho}{\partial \alpha^2}, \; \alpha\in[0,2\pi),
	\label{eq:numm_surfalpha}
\end{align}
for all droplets $k=1,\hdots,n$.

Using the method of lines approach, the right hand side in \eqref{eq:numm_surfalpha} is computed whilst keeping the time $t$ continuous. Define $\rho^k_i(t) \approx \rho(\alpha_i,t)$, where $\alpha_i\in[0,2\pi)$ are the parameter discretisation points for the equidistant discretisation on the interface, $z^k_i = z(\alpha_i,t)$, for the droplet $k$. For each droplet, the PDE can easily be reduced to a system of ODEs for the Fourier coefficients of $\rho^k_i(t)$, $i=1,\hdots,N_k$, obtained by FFTs. In the following only the case of one droplet will be regarded, $k=n=1$, but the method extends trivially to multiple droplets.

It is advantageous to treat the convection and stretching terms explicitly, whilst the diffusion term should be treated implicitly to avoid stiffness. Thus, \eqref{eq:numm_surfalpha} becomes
\begin{align*}
	\dfrac{\partial\rho_i}{\partial t} = f_{E}\left(z_i,\rho_i,t \right) + f_{I} \left(z_i,\rho_i,t\right), \; i=1,\hdots,N,
\end{align*}
where the terms have been collected for explicit and implicit treatment respectively,
\begin{align*}
	f_E(z_i,\rho_i,t) &= \dfrac{\tilde{u}_t}{s^k_\alpha(t)}\dfrac{\partial\rho}{\partial \alpha} - \dfrac{1}{s^k_\alpha(t)}\dfrac{\partial \left(\rho u_t\right)}{\partial \alpha} - \rho u_n \kappa, \\
	f_I(z_i,\rho_i,t) &= \dfrac{1}{Pe_\Gamma \,s_\alpha^k(t)^2}\dfrac{\partial^2\rho}{\partial \alpha^2}.
\end{align*}
The corresponding system of ODEs in Fourier space is
\begin{align*}
	\dfrac{\partial \widehat{\rho}_j}{\partial t} = (\widehat{f_E})_j + (\widehat{f_I})_j, \; j=-\frac{M}{2}, \hdots , \frac{M}{2}-1,
\end{align*}
where $\widehat{(\cdot)}_j$ denotes the Fourier coefficient of mode $j$. For a droplet in two dimensions, the interface discretisation is one dimensional. In one dimension, the implicit term, $\widehat{(f_I)}_j$, is trivial to compute as the second derivative corresponds to a multiplication of $-j^2$ for each Fourier mode. The parts in $\widehat{(f_E)}_j$ are computed pseudo-spectrally. To avoid aliasing the discretisations on all interfaces are upsampled by factor $\frac{3}{2}$, i.e. $M=\frac{3}{2}N$ and a Krasny filter is applied, setting all Fourier modes below the absolute value of $10^{-12}$ are set to zero.

\subsection{Adaptive time-stepping}
\label{sec:methods_time}
\noindent Including the modified tangential velocity and the pseudo-spectral method for the surfactant concentration, the system \eqref{eq:form_system} that needs to be evolved in time can be rewritten as
\begin{align}
	\begin{split}
		\left.\dfrac{dz}{dt}\right\vert_i &= \left(u_n + i\tilde{u}_t\right)_in_i= g(z_i,\sigma_i), \;  i=1,\hdots,N\\
		\dfrac{\partial \widehat{\rho}_j}{\partial t} &= (\widehat{f_E})_j + (\widehat{f_I})_j, \; j=-\frac{M}{2}, \hdots , \frac{M}{2}-1,
	\end{split}
\end{align}
for each droplet $k$. The system is coupled, as the surface tension coefficient $\sigma_i$ depends on $\rho_i$ through some equation of state, e.g. \eqref{eq:form_eqstate_langmuir_nd} or \eqref{eq:form_eqstate_linear_nd}.

As mentioned in \S\ref{sec:methods_surf}, the diffusion term in \eqref{eq:numm_surfalpha} needs to be treated implicitly for stability whilst the other terms in the right hand side preferably are treated explicitly. The evolution equation for the interface should be treated explicitly. Furthermore, adaptivity in time is required in order to handle multiple droplets getting very close to each other.

A range of potential time-stepping schemes for this problem were studied in \cite{Palsson2017}. Unlike the case for vesicles \cite{Quaife2016}, for deforming droplets the time step has a restriction of order one, which enables larger freedom when selecting a suitable time-integration scheme. The design criteria was an accurate but efficient scheme, which would keep the number of solutions of Stokes equations at a minimum as it is the most costly part of the method. The conclusion was that a combination of an explicit midpoint method for $z$ together with an implicit-explicit Runge-Kutta of order two for $\rho$ together with adaptivity was advantageous to use. The authors observed that the adaptive time step should be adjusted considering errors both in interface position $z$ and surfactant concentration $\rho$. An explanation of each part of the scheme follows:

\paragraph{Interface position} The explicit method used to solve
\begin{align*}
	\dfrac{dz}{dt} = g(z,\sigma),
\end{align*}
is the standard Runge-Kutta explicit midpoint rule. At each time step, the local error in $z^{n+1}$, $r_z^{n+1}$, is computed by comparing $z^{n+1}$ with that of a lower order method, $\tilde{z}^{n+1}$. The scheme in full is shown in the Butcher tableau in Table~\ref{tab:num_midp}. The local error for the time step $t_n$ to $t_{n+1}$ is computed as
\begin{align*}
	r^{n+1}_z = \dfrac{\|z^{n+1}-\tilde{z}^{n+1}\|_\infty}{\|z^{n+1}\|_\infty}.
\end{align*}
\begin{table}[h!]
\centering
\begin{tabular}{c | c c}
	$0$ &  &  \\
	$1/2$ & $1/2$ & \\
	\hline
	 & $0$ & $1$ \\
	 & $1$ & $0$ \\
\end{tabular}
\caption{Butcher tableau, adaptive method using explicit midpoint and explicit Euler for interface position $z$.}
\label{tab:num_midp}
\end{table}

\paragraph{Surfactant concentration} The implicit-explicit Runge-Kutta of order two used here (IMEX2) is described by its Butcher tableau in Table~\ref{tab:num_imex2} and explained in detail in \cite{Ascher1997}. Adaptivity in $\rho$ can be achieved in two possible ways, either by comparing to a lower order IMEX scheme or by using the concept of mass conservation from \eqref{eq:form_conservmass}. Here, adaptivity based on mass conservation is considered. The time-discretised version of the mass conservation reads
\begin{align*}
	\int_{\Gamma^{n+1}}\rho^{n+1}dS  - \int_{\Gamma^{n}}\rho^ndS = 0,
\end{align*}
where the integrals can be computed with spectral accuracy using the trapezoidal rule. To measure the local error in $\rho$ of the time step from $t_n$ to $t_{n+1}$, regard
\begin{align*}
	r_{\rho}^{n+1} = \dfrac{\|\int_{\Gamma^{n+1}}\rho^{n+1}dS-\int_{\Gamma^n}\rho^ndS\|_\infty}{\|\int_{\Gamma^n}\rho^ndS\|_\infty}.
\end{align*}
\begin{table}[h!]
\centering
\begin{tabular}{c | c c}
	$0$ &  &  \\
	$1/2$ &  & $1/2$ \\
	\hline
	 & $0$ & $1$ \\
\end{tabular}
\hspace{5mm}
\begin{tabular}{c | c c}
	$0$ &  &  \\
	$1/2$ & $1/2$ & \\
	\hline
	 & $0$ & $1$ \\
\end{tabular}
\caption{Butcher tableau, IMEX2 scheme for $\rho$}
\label{tab:num_imex2}
\end{table}

\paragraph{Coupled method} The Runge-Kutta stages of the methods for $z$ and $\rho$ are the same. This means that at every stage, i.e. at $t_n$ and $t_{n+1/2}$, information between the equations can be exchanged. This coupling, together with the two-second order methods, give as scheme which is second order in total. However, the most important feature is the combined adaptivity. The total local error of a time step is computed as the maximum of the local errors in $z$ and $\rho$,
\begin{align*}
	r^{n+1} = \max\left(r_z^{n+1},r_\rho^{n+1}\right).
\end{align*}
As is typical in adaptive time-stepping methods, if the local error $r^{n+1}>T$ for some tolerance $T$ for a time step of size $dt^n_{old}$, the time step is retaken with updated size. To meet the tolerance, the new time step to go from $t_n$ to $t_{n+1}$, denoted $dt^n_{new}$, needs to fulfil
\begin{align*}
	\left(\dfrac{dt_{new}^n}{dt_{old}^n}\right)^2 \approx \dfrac{T}{r^{n+1}}.
\end{align*}
This gives a new time step of
\begin{align*}
	dt_{new}^n = dt_{old}^n \left( 0.9\dfrac{T}{r^{n+1}}\right)^{\frac{1}{2}},
\end{align*}
where $0.9$ is a safety factor, \cite{Hairer1987}. The expected second order convergence rate is confirmed in \cite{Palsson2017}.

%% file: srcfiles/validation4.tex

\let\oldnl\nl
\newcommand{\nonl}{\renewcommand{\nl}{\let\nl\oldnl}}
\newcommand{\nextnr}{\stepcounter{AlgoLine}\ShowLn}

\section{Validations}
\label{sec:valid}
\noindent To validate the coupled method described in \S\ref{sec:nummethod}, the conformal-mapping techniques of \citeauthor{crowdy2005pair} \cite{crowdy2005pair} and \citeauthor{siegelvalid} \cite{siegelvalid} are used to compute the evolution of surfactant-covered bubbles in an extensional flow for specific cases. The aim of the validation is to test the coupling of the solution for interface position and surfactant concentration as well as the dynamics over time to high accuracy. With this aim, three cases for validation have been selected:
\begin{enumerate}
	\item the exact solutions of surfactant-covered bubbles and droplets in extensional flow at steady state
	\item semi-analytic solutions for a pair of clean, deforming bubbles in extensional flow
	\item semi-analytic solutions for a pair of surfactant-covered bubbles deforming in extensional flow.
\end{enumerate}
The first case uses an analytical solution by \citeauthor{siegelvalid} and tests the dynamics of surfactant-covered bubbles and droplets over a long time. The second method was introduced by \citeauthor{crowdy2005pair} and regards a pair of clean bubbles in extensional flow.  This case tests especially the special quadrature of \S\ref{sec:specq} when the two bubbles are pushed close to each other. The third case builds on the second, with the addition of insoluble surfactants. This tests the time-dependent dynamics of the whole method put together.

In this section, an overview of how to compute the validation results for all three cases is described. Note that case two and three are both described in \S\ref{sec:cleanpair}, where the former simplifies to no surfactants $\rho$ and a constant surface tension $\sigma$.

\input{srcfiles/valid_steady}

\input{srcfiles/valid_pair}

%% file: srcfiles/valid_steady.tex

\subsection{Surfactant-covered bubble or drop in steady state}
\label{sec:validsteady}
\noindent For single surfactant-covered bubbles, analytical solutions to interface position $z$ and surfactant concentration $\rho$ exist in certain cases.
The case regarded here is that of a bubble in an extensional flow, $\mathbf{u}_\infty = Q(x,-y)$, at \textit{steady state}, for which exact solutions were originally obtained by \citeauthor{siegelvalid} \cite{siegelvalid}. A steady state is defined as the time $T$ when $\mathbf{u}\cdot\mathbf{n}=0$. These solutions exist for bubbles deforming from an initially circular shape with uniform surfactant concentration and where there is no diffusion of surfactants along the interface, i.e. $Pe_\Gamma = \infty$. To couple surfactant concentration to surface tension, the linear equation of state \eqref{eq:form_eqstate_linear_nd} is used. For details and derivations on how to obtain these solutions, the reader is referred to the original reference. Here, it will be explained how to compute these solutions to use as a validation case. Given a deformation at steady state, i.e. $z$ at time $T$, these solutions are used to analytically compute $\rho$ at time $T$ as well as the Capillary number $Q$ that is needed to obtain the set deformation at steady state. Using a numerical method to simulate the deformation and surfactant evolution under this Capillary number, the solution at steady state can then be compared to the analytical solution to determine the accuracy.

To describe the bubble interface, a conformal mapping is used,
\begin{align}
	z(\zeta,t) = \dfrac{a(t)}{\zeta} + b(t)\zeta,
	\label{eq:siegel_confmap}
\end{align}
where $a(t)<0$ and $b(t)$ is real. Furthermore, the condition of constant area gives $a^2-b^2 = 1$. At the bubble interface, $|\zeta|=1$, i.e. $\zeta=e^{i\nu}$ for $\nu\in[0,2\pi]$. At steady state, it is required that
\begin{align}
	Q  = -\dfrac{Ab}{(1+b^2)^{1/2}},
	\label{eq:siegel_Q}
\end{align}
where
\begin{align}
	A = \dfrac{\int_0^{2\pi}B(b,\nu)^{1/2}d\nu-2\pi E}{\int_0^{2\pi}B(b,\nu)d\nu},
	\label{eq:siegel_A}
\end{align}
with $E$ the elasticity number in \eqref{eq:form_eqstate_linear_nd}, and
\begin{align}
	B(b,\nu) = 1 + 2b^2-2(1+b^2)^{1/2}b\cos(2\nu).
	\label{eq:siegel_B}
\end{align}
The deformation of a bubble is measured as
\begin{align}
	D=\dfrac{R_{max}-R_{min}}{R_{max}+R_{min}},
	\label{eq:siegel_D}
\end{align}
where $R_{max}$, $R_{min}$ are the maximum and minimum radial distance of a point on the interface to the centre respectively.

To compute the surfactant concentration at steady state, \eqref{eq:form_surfacteq_nondim} is used. In steady state, this can be simplified to
\begin{align}
	\rho \Re\left[ \dfrac{u(z)+\overline{z}_\nu}{|z_\nu|}\right] = \dfrac{1}{Pe_\Gamma}\left(\dfrac{\rho_\nu}{|z_\nu|}\right) = 0,
\end{align}
where the second equality comes from having no diffusion, $Pe_\Gamma=\infty$. This means, that on any part of the interface at steady state, either $\rho = 0$ or $u(z) = 0$ (no slip condition). Assuming the extensional flow is sufficiently low such that the concentration of surfactants is non-zero everywhere, it fulfils
\begin{align}
	\rho(\nu) = \dfrac{1-A|z_\nu|}{E}.
	\label{eq:siegel_rho}
\end{align}
An algorithm to compute the steady state surfactant concentration $\rho(\nu)$ and Capillary number $Q$ given the steady state deformation can be found in Algorithm~\ref{alg:steady}. The deformation at steady state is determined through the conformal mapping parameters $a$ and $b$ in \eqref{eq:siegel_confmap} which are given as input to the algorithm. How the deformation $D$ depends on the Capillary number $Q$ for different elasticity numbers is shown in Figure~\ref{fig:plot4}.

To validate the results of a numerical method against these exact solutions, the results at steady state are compared when simulating an initially circular bubble with uniform surfactant concentration $\rho_0 = 1$ in an extensional flow with Capillary number $Q$ and $Pe_\Gamma=\infty$, i.e. no diffusion of surfactants.

\IncMargin{1em}
\RestyleAlgo{boxruled}
\SetCommentSty{emph}
\SetKwComment{Comment}{}{}
\DontPrintSemicolon
\LinesNumbered
\begin{algorithm}[H]
\SetKwInOut{Input}{Input}\SetKwInOut{Output}{Output}
\Input{Conformal mapping parameters $a$ and $b$ at steady state.}
\Output{Deformation $D$. Capillary number $Q$ for which this deformation is obtained at steady state. Interface position $z(\nu_j)$ and surfactant concentration $\rho(\nu_j)$ at steady state, for $\nu_j=\frac{2\pi j}{M}$, $j\in[0,M-1]$.}
\vspace{2mm}
Discretise $\nu_j = \dfrac{2\pi j}{M}$, $j\in[0,M-1].$ \;
Compute $z(\nu_j)$ through \eqref{eq:siegel_confmap}.\;
Find $R_{max}$, $R_{min}$ and compute deformation $D$ in \eqref{eq:siegel_D}. \;
Compute $B(b,\nu_j)$ in \eqref{eq:siegel_B}.\;
Compute $A$ in \eqref{eq:siegel_A} using trapezoidal rule. \;
Compute $\rho(\nu_j)$ with \eqref{eq:siegel_rho}. \;
Compute $Q$ from \eqref{eq:siegel_Q}. \;
\caption{Steady state computation}
\label{alg:steady}
\end{algorithm}\DecMargin{1em}

As was discussed by \citeauthor{Milliken1993a} in \cite{Milliken1993a}, in the case of no diffusion along the interface the viscosity ratio does not affect steady state deformation. When $Pe_\Gamma = \infty$, the interfacial velocity is dominated by the high Marangoni stresses and the effect of viscosity ratio is therefore negligible.  Thus, the results above can be used to validate the steady state deformation and surfactant concentration also for drops, where $\lambda \neq 0$.

%% file: srcfiles/valid_pair.tex

\subsection{A pair of bubbles in extensional flow}
\label{sec:cleanpair}
\noindent This method of using complex-variable formulations to semi-analytically compute the deformation of a pair of bubbles in an extensional flow was introduced by \citeauthor{crowdy2005pair} \cite{crowdy2005pair}. There, the specific problem setting was a pair of bubbles, $\lambda_k=0$ for $k=1,2$, without surfactants. Here, the method is extended to also include surfactant-covered bubbles, using a linear equation of state \eqref{eq:form_eqstate_linear_nd}. The imposed far-field flow is the extensional flow $\mathbf{u}_\infty = Q(x,-y)$, where $Q$ is the Capillary number as previously defined. This method considers the special case when two bubbles are reflectionally symmetric about both $x$- and $y$-axes, see Figure~\ref{fig:crowdy_z} (left).

The method is based on the reformulation of Stokes equations using complex-variable methods. The bubble boundaries can then be parametrised in terms of a conformal mapping. The numerical algorithm boils down to solving a system of ODEs for the time evolution of the conformal mapping parameters. At each instance of time, the flow field is computed on the unit disc in the conformal-mapping space, using Laurent series.

For a detailed explanation of each of these steps, the reader is referred to \cite{crowdy2005pair}. Here, only the numerical approach will be explained in detail. The equations follow those of \cite{crowdy2005pair}, with two modifications: firstly, the nondimensionalisation has been adjusted to allow for non-uniform surface tension, $\sigma$, along the interface. This is in order to make the extension to surfactant-covered bubbles. Secondly, the convection-diffusion equation for surfactants \eqref{eq:form_surfacteq_nondim} is coupled to the problem. When considering the original case of clean bubbles, the surface tension is set to be constant and \eqref{eq:form_surfacteq_nondim} does not need to be solved.

The fluid domain and bubble interfaces are described by the conformal mapping,
\begin{align}
	z\left( \zeta,t \right) = \dfrac{b(t)}{\zeta-\sqrt{\phi(t)}} + \hat{z}\left(\zeta,t\right),
	\label{eq:crowdy_confmap}
\end{align}
where
\begin{align}
	\hat{z}\left(\zeta,t\right) = \sum_{n=-\infty}^\infty a_n(t)\zeta^n.
	\label{eq:crowdy_zhat}
\end{align}
Furthermore, the relations
\begin{align}
	a_0(t) = \dfrac{b(t)}{2\sqrt{\phi(t)}}, \quad a_{-n}(t)=-\phi(t)^na_n(t), \; n\geq 1,
	\label{eq:crowdy_an}
\end{align}
are known, $0<\phi(t)<1$ and $b(t)$ is real. On the bubble interfaces, $\zeta=e^{i\nu}$, for $\nu\in[0,2\pi]$. In order to also consider the surface tension in the $\zeta$-plane, the composite functions $\sigma(\zeta,t) = \sigma(z(\zeta,t),t)$ for surface tension and $\rho(\zeta,t)=\rho(z(\zeta,t),t)$ for surfactants are defined.

In the numerical implementation, all Laurent series expansions are truncated to only include terms $\left\{\zeta^n,\, n=-N_V,\hdots,N_V\right\}$. When zero-padding is needed to avoid aliasing, a factor two is used: giving a new truncation limit $M_V=2N_V$. The upper bubble interface is discretised by $\zeta_j = e^{i\nu_j}$, for $j=1,\hdots,2M_V+1$, where $\nu_j$ is equidistant in $[0,2\pi)$ with $2M_V+1$ discretisation points.

\subsubsection{Computing the flow field}
\noindent At a given time $\bar{t}$, the interface position of bubbles, $z(\zeta,\bar{t})$, its conformal mapping parameters $b(\bar{t})=:b$, $\phi(\bar{t}) =: \phi$, $\left\{a_n(\bar{t}),\, n=-N_V,\hdots,N_V\right\}$ and $\sigma(\zeta,\bar{t})$ are known.
\\ \\
In \cite{crowdy2005pair}, it is shown that the flow field at $\zeta$ at time $\bar{t}$ can be described by functions $F(\zeta,\bar{t})$, $G(\zeta,\bar{t})$ and $C(\bar{t})$. On the upper bubble interface, $|\zeta|=1$, these functions are determined by solving the equation
\begin{align}
	A(\zeta)F_0 + iF(\zeta,\bar{t}) + P(\zeta) -iG(\zeta^{-1},\bar{t}) - i C(\bar{t}) = E(\zeta),
	\label{eq:crowdy_FGCreform}
\end{align}
where
\begin{align}
	\begin{cases}
		A(\zeta) = 2\sqrt{\phi}\left( \dfrac{i}{\zeta-\sqrt{\phi}}-\dfrac{iz(\zeta,\bar{t})}{z_\zeta(\zeta^{-1},\bar{t})\left(\zeta^{-1}-\sqrt{\phi}\right)^2}\right), \\
		E(\zeta) = \sigma(\zeta,\bar{t})\dfrac{\zeta z_\zeta(\zeta,\bar{t})}{2|z_\zeta(\zeta,\bar{t})|}+\dfrac{iQb}{\zeta^{-1}-\sqrt{\phi}}, \\
		P(\zeta) = B(\zeta)F_\zeta(\zeta^{-1},\bar{t}), \; \text{ where } \; B(\zeta) = \dfrac{iz(\zeta,\bar{t})}{z_\zeta(\zeta^{-1},\bar{t})}.
	\end{cases}
	\label{eq:crowdy_ABE}
\end{align}
Both $F(\zeta,\bar{t})$ and $G(\zeta,\bar{t})$ can be expanded into Laurent series,
\begin{align*}
	F(\zeta,\bar{t}) = \sum_{n=-\infty}^{\infty} F_n\zeta^n, \quad G(\zeta,\bar{t}) = \sum_{n=-\infty}^{\infty} G_n\zeta^n,
\end{align*}
for real coefficients $\{F_n\}$ and $\{G_n\}$.

Expanding \eqref{eq:crowdy_FGCreform} into a Laurent series on both sides of the equality sign, and using the orthogonality of $\zeta=e^{i\nu}$, gives $2N_V+1$ equations:
\begin{align}
	A_nF_0 + iF_n + P_n - iG_{-n} -i\delta_nC(\bar{t})  = E_n, \quad n\in[-N_V,N_V],
	\label{eq:crowdy_CFGn}
\end{align}
where $\delta_n := 1$ for $n=0$ and $0$ otherwise. To compute the flow field at a given time $\bar{t}$, the Laurent coefficients $\{F_n, \; \forall n\in[-N_V,N_V]\}$ and $\{G_n,\; \forall n\in[-N_V,N_V]\}$ need to be determined together with $C(\bar{t})$.
The coefficients $\{F_n\}$, $\{G_n\}$ obey
\begin{align}
	F_{-n} &=-\phi^nF_n, \; n\geq 1, \text{ and } \label{eq:crowdy_Fn} \\
	G_{-n} &=-\phi^nG_n, \; n\geq 1, \; G_0 = \dfrac{Qb}{2\sqrt{\phi}}, \;\label{eq:crowdy_Gn}
\end{align}
thus the unknowns are reduced to $\{F_n, \; \forall n\in[0,N_V]\}$, $\{G_n,\; \forall n\in[1,N_V]\}$ and $C(\bar{t})$. Note that since $Q$ is known, $G_0$ is known for time $\bar{t}$. Therefore, in \eqref{eq:crowdy_CFGn}, $G_0$ should be moved to the right hand side for $n=0$.

The number of unknowns in the problem is thus $2N_V+2$, and the vector of unknowns is
\begin{align}
	x = \begin{pmatrix}
		F_0 & F_1 & \hdots & F_n & G_1 & \hdots & G_n & C(\bar{t})
	\end{pmatrix}^T.
	\label{eq:crowdy_x}
\end{align}
From \eqref{eq:crowdy_CFGn}, $2N_V+1$ equations are obtained. The system is completed with the equation
\begin{align}
	G_\zeta\left(\sqrt{\phi},\bar{t}\right) &= \sum_{n=1}^{N_V} 2n\phi^{\frac{n-1}{2}}G_n = Q\hat{z}_\zeta\left(\sqrt{\phi},\bar{t}\right).
	\label{eq:crowdy_Gz}
\end{align}
which originates from the presence of an irrotational extensional flow at infinity. As typically $\phi<1$, the expressions above need to be computed using the relations between $a_n$ and $a_{-n}$ \eqref{eq:crowdy_an} as well as between $G_n$ and $G_{-n}$ \eqref{eq:crowdy_Gn} for stability.

Put together, to compute $x$ in \eqref{eq:crowdy_x} it is necessary to solve a system $Mx=b$ of size $2N_V+2$, where $Mx$ corresponds to the left hand side in \eqref{eq:crowdy_CFGn} together with \eqref{eq:crowdy_Gz}, and
\begin{align}
	b = \begin{pmatrix} E_{-N_V} & \hdots & E_{-1} & E_0+iG_0 &E_1 & E_{N_V} & Q\hat{z}_\zeta(\sqrt{\phi},\bar{t}) \end{pmatrix}^T.
	\label{eq:crowdy_RHS}
\end{align}
The algorithm to compute $x$ is described in Algorithm~\ref{alg:crowdy_1} with $M_V=2N_V$. Despite the fact that the system matrix $M$ has a reasonable condition number at around $10^5$ the system $Mx=b$ is sensitive, especially to the value of $\phi$. The authors note that not all conventional approaches to solve the linear system will give a satisfying solution.  It has been found numerically, however, that treating the matrix as sparse helps with stability.
\\ \\
\IncMargin{1em}
\RestyleAlgo{boxruled}
\SetCommentSty{emph}
\DontPrintSemicolon
\LinesNumbered
\begin{algorithm}[H]
\SetKwInOut{Input}{Input}\SetKwInOut{Output}{Output}
\Input{Conformal mapping parameters: $b$, $\phi$ and $\{a_n(\bar{t}),n\in[-N_V,N_V]\}$. Surface tension $\sigma(\zeta_j,\bar{t})$, $j=1,\hdots 2M_V+1$.}
\Output{Coefficients $\{F_n,n\in[0,N_V]\}$ and $\{G_n,n\in[1,N_V]\}$, $C\left(\bar{t}\right)$.}
\vspace{2mm}
Zero-pad $\{a_n\}_{n=-N_V}^{N_V}$ to $\{a_m\}_{m=-M_V}^{M_V}$.\;
Compute $z(\zeta_j,\bar{t})$ through \eqref{eq:crowdy_confmap}. Similarly, compute $z(\zeta_j^{-1},\bar{t})$, $z_\zeta(\zeta_j,\bar{t})$ and $z_\zeta(\zeta_j^{-1},\bar{t})$. \;
Compute $A(\zeta_j)$, $B(\zeta_j)$ and $E(\zeta_j)$, $\forall j\in[1,2M_V+1]$, in \eqref{eq:crowdy_ABE}. \;
Compute Laurent series coefficients for $A(\zeta)$, $B(\zeta)$ and $E(\zeta)$ using \texttt{FFT}s. \;
Truncate Laurent series $\{A_m\}_{-M_V}^{M_V}$, $\{E_m\}_{-M_V}^{M_V}$ to include only $2N_V+1$ terms. \;
Create right hand side $b$ to solve $Mx=b$, constructed according to \eqref{eq:crowdy_RHS}. \;
Create system matrix $M$: \;
\nonl The first $2N_V+1$ rows correspond to \eqref{eq:crowdy_CFGn}. To find the terms in $P_n$, the following expression is used
	\begin{align}
		P_n = \sum_{k=-N_V}^{N_V}B_{n+k-1}kF_k, \; \forall n\in[-N_V,N_V].
	\end{align}
	For all terms, the relations \eqref{eq:crowdy_Fn} and \eqref{eq:crowdy_Gn} are used. \;
\nonl The last row of the matrix corresponds to \eqref{eq:crowdy_Gz}. \;
Solve $Mx=b$. \;
Define function coefficients: $\{F_n\}_0^{N_V} =(x_1,\; \hdots,\; x_{N_V+1})$, $\{G_n\}_1^{N_V} = (x_{N_V+2},\; \hdots,\; x_{2N_V+1})$ and $C(\bar{t})=x_{2N_V+2}$.
\caption{Compute flow field functions: $F$, $G$ and $C(\bar{t})$}
\label{alg:crowdy_1}
\end{algorithm}\DecMargin{1em}

\subsubsection{Evolving the conformal mapping parameters in time}
\noindent In order to evolve the deforming bubble interfaces, $z(\zeta,t)$, $|\zeta|=1$, the conformal mapping parameters $b(t)$, $\phi(t)$ and $\{a_n(t)\}_{n=-N_V}^{N_V}$ are integrated in time.

Discretising in time, at time $t_k$ define $b^k \approx b(t_k)$, $\phi^k \approx \phi(t_k)$, $\{a_n^k\} \approx \{a_n(t_k)\}$ etc. The differential equations to evolve the bubble interfaces are
\begin{align}
	\begin{split}
	\dfrac{d\hat{z}}{dt}(\zeta,t) =& \zeta\hat{z}_\zeta(\zeta,t) - \dfrac{b(t)\zeta I(\zeta,t) - b(t)\sqrt{\phi(t)}I(\sqrt{\phi(t)},t)}{(\zeta-\sqrt{\phi(t)})^2} - 2F(\zeta,t) \\
	\hdots&+\dfrac{b(t)I(\sqrt{\phi(t)},t)+b(t)\sqrt{\phi(t)}I_\zeta(\sqrt{\phi(t)},t)}{\zeta-\sqrt{\phi(t)}}  := f(\zeta,t),
	\end{split} \label{eq:crowdy_dzhat} \\
	\dfrac{d\phi}{dt} &= -2\phi(t)I(\sqrt{\phi(t)},t) =: g(t), \label{eq:crowdy_drho}
\end{align}
where $I(\zeta,t_k)$ is known and $\hat{z}$ is defined in \eqref{eq:crowdy_zhat}. Regarding \eqref{eq:crowdy_dzhat} and taking the Laurent series expansion of both sides, this can be written as a system of ODEs, as
\begin{align}
	\dfrac{d}{dt}a_n = f_n(t), \; n\in[-N_V,N_V],
	\label{eq:crowdy_dzhatn}
\end{align}
where $\{f_n(t)\}$ are the Laurent series coefficients of $f(\zeta,t)$ above. The function $I(\zeta,t)$ is obtained through the kinematic condition \eqref{eq:form_kincond} by first computing
\begin{align}
	D(\zeta) = \dfrac{1}{2|z_\zeta|} + \Re\left[\dfrac{iC}{\zeta z_\zeta}\right],
	\label{eq:crowdy_D}
\end{align}
and taking a Laurent series expansion $D(\zeta) = \sum_{n=-\infty}^\infty D_n\zeta^n$. The coefficients of $I$ are obtained through
\begin{align}
	\begin{cases}
	I_n = 2\left(\dfrac{1+\phi^n}{1-\phi^{2n}}\right)D_n, \; n \geq 1\\
	I_0 = D_0 \\
	I_{-n} = -\left(\dfrac{2\phi^n}{1-\phi^n}\right)D_n, \; n\geq 1.
	\end{cases}
	\label{eq:crowdy_dn}
\end{align}
The last conformal mapping parameter, $b(t)$, can be decided through the constant area constraint, i.e.
\begin{align}
	\mathcal{A}(t) = -\dfrac{1}{2i}\oint_{|\zeta|=1} \overline{z}(\zeta^{-1},t)z_\zeta(\zeta,t)d\zeta = \pi.
	\label{eq:crowdy_A}
\end{align}
Each bubble will have an area of $\pi$ due to the non-dimensionalisation. Expanding \eqref{eq:crowdy_A} the following expression for $a$ is obtained:
\begin{align}
\begin{split}
	b^2&\oint_{|\zeta|=1}\dfrac{d\zeta}{(\zeta^{-1}-\sqrt{\phi})(\zeta-\sqrt{\phi})^2} + b\oint_{|\zeta|=1}\dfrac{\hat{z}(\zeta^{-1},t)}{(\zeta-\sqrt{\phi})^2}-\dfrac{\hat{z}_\zeta(\zeta,t)}{(\zeta^{-1}-\sqrt{\phi})}d\zeta \\
	&- \oint_{|\zeta|=1}\hat{z}(\zeta^{-1},t)\hat{z}_\zeta(\zeta,t)d\zeta = 2i\pi.
\end{split}
\label{eq:crowdy_A2}
\end{align}
How to compute the $f(\zeta,t)$ and $g(t)$ in \eqref{eq:crowdy_dzhat} and \eqref{eq:crowdy_drho} is described in Algorithm~\ref{alg:crowdy_2}.
\\ \\
\IncMargin{1em}
\RestyleAlgo{boxruled}
\SetCommentSty{emph}
\SetKwComment{Comment}{}{}
\DontPrintSemicolon
\LinesNumbered
\begin{algorithm}[H]
\SetKwInOut{Input}{Input}\SetKwInOut{Output}{Output}
\Input{Conformal mapping parameters at time $t^k$: $b^k$, $\phi^k$, $\{a_n^k,n\in[-N_V,N_V]\}$. Flow functions and surface tension at time $t^k$: $C(t_k)$, $\left\{F_n^k\right\}_{n=0}^{N_v}$, $\left\{G_n^k\right\}_{n=1}^{N_V}$ and $\sigma_j^k$, $j=1,\hdots,2N_V+1$.}
\Output{$\{f_n(t_k)\}_{n=-N_V}^{N_V}$, $g(t_k)$.}
\vspace{2mm}
Zero-pad $\{a_n^k\}_{n=-N_V}^{N_V}$ to $\{a_m^k\}_{m=-M_V}^{M_V}$.\;
Compute $z_\zeta(\zeta_j,t_k)$, $\hat{z}_\zeta(\zeta_j,t_k)$ and $\hat{z}(\zeta^{-1},t_k)$ through \eqref{eq:crowdy_confmap}, $\forall j\in[1,2M_V+1]$.\;
Compute $D(\zeta_j)$ in \eqref{eq:crowdy_D}, $\forall j\in[1,2M_V+1]$.\;
Compute $\{d_m\}_{m=-M_V}^{M_V}$ through \texttt{FFT} of $D(\zeta_j)$. \;
Compute $\{I_m\}_{m=-M_V}^{M_V}$ from $\{d_m\}_{-M_V}^{M_V}$ as in \eqref{eq:crowdy_dn}. \;
Using $\{I_m\}$, compute $I(\zeta_j,t_k)$, $I(\sqrt{\phi},t_k)$ and $I_\zeta(\sqrt{\phi},t)$. \;
Compute $f$ according to \eqref{eq:crowdy_dzhat}, compute Laurent series coefficients $\{f_m\}_{m=-M_V}^{M_V}$  by \texttt{FFT} and truncate to $2N_V+1$ terms, $\{f_n\}_{n=-N_V}^{N_V}$.\;
Compute $g$ through \eqref{eq:crowdy_drho}.\;
\caption{Computation of RHS in \eqref{eq:crowdy_dzhat}, \eqref{eq:crowdy_drho}.}
\label{alg:crowdy_2}
\end{algorithm}\DecMargin{1em}

\subsubsection{Computing the surfactant concentration}
\noindent On the bubble boundaries, $\zeta=e^{i\nu}$, the convection-diffusion equation for the surfactant concentration \eqref{eq:form_surfacteq_nondim} can be written in terms of $\nu\in[0,2\pi)$ as
\begin{align}
\begin{split}
	\left. \dfrac{\partial\rho}{\partial t}\right\rvert_\nu - \Re\left(\dfrac{\rho_\nu}{z_\nu}z_t\right) &+ \dfrac{1}{|z_\nu|}\dfrac{\partial}{\partial\nu}\Re\left(P(\nu,t)\right) - \dfrac{1}{|z_\nu|}\Im\left(\dfrac{z_{\nu\nu}}{z_\nu}\right)\Im\left(P(\nu,t)\right) \\
	& -\dfrac{1}{|z_\nu|}\dfrac{1}{Pe_\Gamma}\dfrac{\partial}{\partial\nu}\left(\dfrac{\rho_\nu}{|z_\nu|}\right) = 0,
\end{split}
\label{eq:valid_surf}
\end{align}
where
\begin{align}
	P(\nu,t) = \dfrac{u(\zeta,t)\overline{z}_\nu \rho(t)}{|z_\nu|}.
	\label{eq:valid_P}
\end{align}
The flow velocity $u(\zeta,t)=u(z(\zeta,t),t)$ is obtained from the normal stress balance \eqref{eq:form_normalstressnondim},
\begin{align}
	u(\zeta,t) = \dfrac{\sigma(\zeta,t)}{2}\dfrac{\zeta z_\zeta}{|z_\zeta|}+iC(t)-\dfrac{4i\sqrt{\phi}F_0}{\zeta-\sqrt{\phi}} -2iF(\zeta,t),
\label{eq:valid_usurf}
\end{align}
where $F(\zeta,t)$ and $F_0$ as defined previously. Furthermore, $z_t$ is obtained from the kinematic condition \eqref{eq:form_kincond}
\begin{align}
	z_t = \zeta z_\zeta I(\zeta,t) - -\dfrac{4i\sqrt{\phi}F_0}{\zeta-\sqrt{\phi}} -2iF(\zeta,t),
\label{eq:valid_zt}
\end{align}
for $I(\zeta,t)$ as defined in \eqref{eq:crowdy_dn}.

As in \S\ref{sec:methods_surf}, the convection and diffusion parts of \eqref{eq:valid_surf} need to be treated explicitly and implicitly respectively. The equation becomes
\begin{align}
	\left. \dfrac{\partial\rho}{\partial t}\right\rvert_\nu = f_{exp}(\nu,\rho,t) + f_{imp}(\nu,\rho,t),
\end{align}
where
\begin{align}
\begin{split}
f_{exp}(\nu,\rho,t) =&  \Re\left(\dfrac{\rho_\nu}{z_\nu}z_t\right) - \dfrac{1}{|z_\nu|}\dfrac{\partial}{\partial\nu}\Re\left(P(\nu,t)\right) \\
&+ \dfrac{1}{|z_\nu|}\Im\left(\dfrac{z_{\nu\nu}}{z_\nu}\right)\Im\left(P(\nu,t)\right),
\end{split}
\label{eq:valid_fexp}
\end{align}
and
\begin{align}
f_{imp}(\nu,\rho,t) = \dfrac{1}{|z_\nu|}\dfrac{1}{Pe_\Gamma}\dfrac{\partial}{\partial\nu}\left(\dfrac{\rho_\nu}{|z_\nu|}\right).
\label{eq:valid_fimp}
\end{align}
How to compute $f_{exp}$ is shown in Algorithm~\ref{alg:crowdy_fexp}. The function $f_{imp}$ is computed in a similar way. In the time-integration scheme, the solution is found through solving a system using \texttt{gmres}.
\\ \\
\IncMargin{1em}
\RestyleAlgo{boxruled}
\SetCommentSty{emph}
\SetKwComment{Comment}{}{}
\DontPrintSemicolon
\LinesNumbered
\begin{algorithm}[H]
\SetKwInOut{Input}{Input}\SetKwInOut{Output}{Output}
\Input{At time $t_k$: conformal mapping parameters; $b^k$, $\phi^k$, $\{a_n^k,n\in[-N_V,N_V]\}$, flow field; $\{F_n^k\}_0^{N_V}$, $\{G_n^k\}_1^{N_V}$ and $C(t_k)$, and surfactant concentration, $\rho^k_j$, $j\in[0,2N_V+1]$.}
\Output{$f_{exp,j}(t_k)$}
\vspace{2mm}
Zero-pad $\{a_n^k\}_{n=-N_V}^{N_V}$ to $\{a_m^k\}_{m=-M_V}^{M_V}$.\;
Using $\{a_m^k\}$, compute $z^k(\zeta_j)$, $z^k_\nu(\zeta_j)$ and $z_{\nu\nu}^k(\zeta_j)$ through \eqref{eq:crowdy_confmap}.\;
Zero-pad Fourier coefficients of $\rho_j^k$ to $\{\hat{\rho}_m^k\}_{m=-M_V}^{M_V}$ and compute $\rho_\nu^k(\zeta_j)$ from coefficients. \;
Compute $F(\zeta_j,t_k)$ from $\{F_n^k\}_0^{N_V}$. \;
Compute $u(\zeta_j,t_k)$ from \eqref{eq:valid_usurf}. \;
Compute $z_t(\zeta_j,t_k)$ from \eqref{eq:valid_zt}. \;
Compute $P(\nu_j,t_k)$ through \eqref{eq:valid_P}. Derivative $\frac{\partial}{\partial\nu}\Re(P(\nu,t))$ computed via \texttt{FFT}. \;
Compute $f_{exp}(\nu_j,\rho_j^k,t_k)$ from \eqref{eq:valid_fexp}.
\caption{Computation of $f_{exp}$ in \eqref{eq:valid_fexp}.}
\label{alg:crowdy_fexp}
\end{algorithm}\DecMargin{1em}

\subsubsection{Time-integration of the system}
\noindent To evolve the system with conformal mapping parameters and surfactant concentration a suitable time-integration method is needed. In this paper, the second order adaptive method in \S\ref{sec:methods_time} is used. For simplicity, in Algorithm~\ref{alg:crowdy_time} a first order method with fixed time step is employed to demonstrate the method. In the case of clean drops, $\sigma(\zeta,t)$ is set to be constant and the algorithm skips steps 9 to 12. At each time step, a Krasny filter is applied where all Fourier modes smaller than $10^{-12}$ are set to zero.
\\ \\
\IncMargin{1em}
\RestyleAlgo{boxruled}
\SetCommentSty{emph}
\SetKwComment{Comment}{}{}
\DontPrintSemicolon
\LinesNumbered
\begin{algorithm}[H]
\SetKwInOut{Input}{Input}\SetKwInOut{Output}{Output}
\Input{Conformal mapping parameters: $b^0$, $\phi^0$, $\{a_n^0,n\in[-N_V,N_V]\}$. Initial surfactant concentration $\rho_j^0$ and surface tension $\sigma_j^0$, $j=1,\hdots,2N_V+1$. Time interval $[t_0,t_s]$, time step $dt$.}
\Output{Conformal mapping parameters at time $t_s$: $b(t_s)$, $\phi(t_s)$, $\{a_n(t_s),n\in[-N_V,N_V]\}$. Surfactant concentration at time $t_s$: $\rho(\zeta_j,t_s)$, $j=1,\hdots,2N_V+1$.}
\vspace{2mm}
$N_t = \dfrac{t_s - t_0}{dt}$ \;
\vspace{0.5mm}
\For {$k \in[0,N_t-1]$}{
	Compute flow field functions at time $t_k$: \\ \nonl \hspace{3mm}
	$\left[\{F_n^k\}_{n=0}^{N_V},\{G_n^k\}_{n=0}^{N_V},C^k\right] =$ \texttt{algorithm2}$\left(b^k, \{a_n^k\}_{-N_V}^{N_V},\phi^k,\sigma_j^k\right).$ \;
	Compute RHS in \eqref{eq:crowdy_dzhat} and \eqref{eq:crowdy_drho}: \\ \nonl \hspace{3mm} $[\{f_n\}_{n=-N_V}^{N_V}, g]=$ \texttt{algorithm3}$\left(b^k,\phi^k,\{a_n^k\}_n, \sigma^k_j,\left\{F_n^k\right\}_n,\left\{G_n^k\right\}_n,C^k\right)$ \;
	Update $\{a_n\}$ in discretised version of \eqref{eq:crowdy_dzhatn}: \\ \nonl \hspace{3mm} $a_n^{k+1} = a_n^k + dt\cdot f_n$, $n\in[1,N]$.\;
	Update $\phi$ in discretised version of \eqref{eq:crowdy_drho}: \\ \nonl \hspace{3mm} $\phi^{k+1} = \phi^k + dt\cdot g$.\;
	Compute $b^{k+1}$ through \eqref{eq:crowdy_A2}.\;
	Compute $a_n^{k+1}$ for $n\leq 0$ through \eqref{eq:crowdy_an}.\;
	\Comment{Steps $9-12$ are skipped in case of clean bubbles:}
	Compute $f_{exp}(t_k)$ in \eqref{eq:valid_fexp}, $\forall j\in[1,2M_V+1]$: \\ \nonl \hspace{3mm}$f_{exp,j}^k=$ \texttt{algorithm4}$\left(b^k,\phi^k,\{a_n^k\},\{F_n^k\},\{G_n^k\},C(t_k),\rho_j^k \right)$.\;
	Compute $\rho^{k+1}_j$ from $\rho^{k+1}_j=\rho^k_j + dtf_{exp,j}^k + dtf_{imp,j}^{k+1}$ with \texttt{gmres}, $f_{imp}$ as in \eqref{eq:valid_fimp}.  \;
	Truncate zero-padding from $\rho_j^{k+1}$. \;
	Update surface tension coefficient $\sigma_j^{k+1}$ through equation of state. \;
	Update time $t_{k+1} = t_k + dt$.\;
}
\caption{First order method for time-integration of surfactant-covered bubbles using fixed time step.}
\label{alg:crowdy_time}
\end{algorithm}\DecMargin{1em}

%% file: srcfiles/results.tex

\section{Results}
\label{sec:results}
\noindent In this section the numerical method suggested in this paper is tested against the validation cases in \S\ref{sec:valid}. Each test case is described in detail below. The numerical method described in \S\ref{sec:nummethod} will be denoted the BIE method for the remainder of the paper.

\input{srcfiles/steady}

\input{srcfiles/cleanvalid}

\input{srcfiles/surfvalid}

\subsection{The Swiss roll}
\begin{figure}[h!]
\centering
\includegraphics[width=0.6\textwidth]{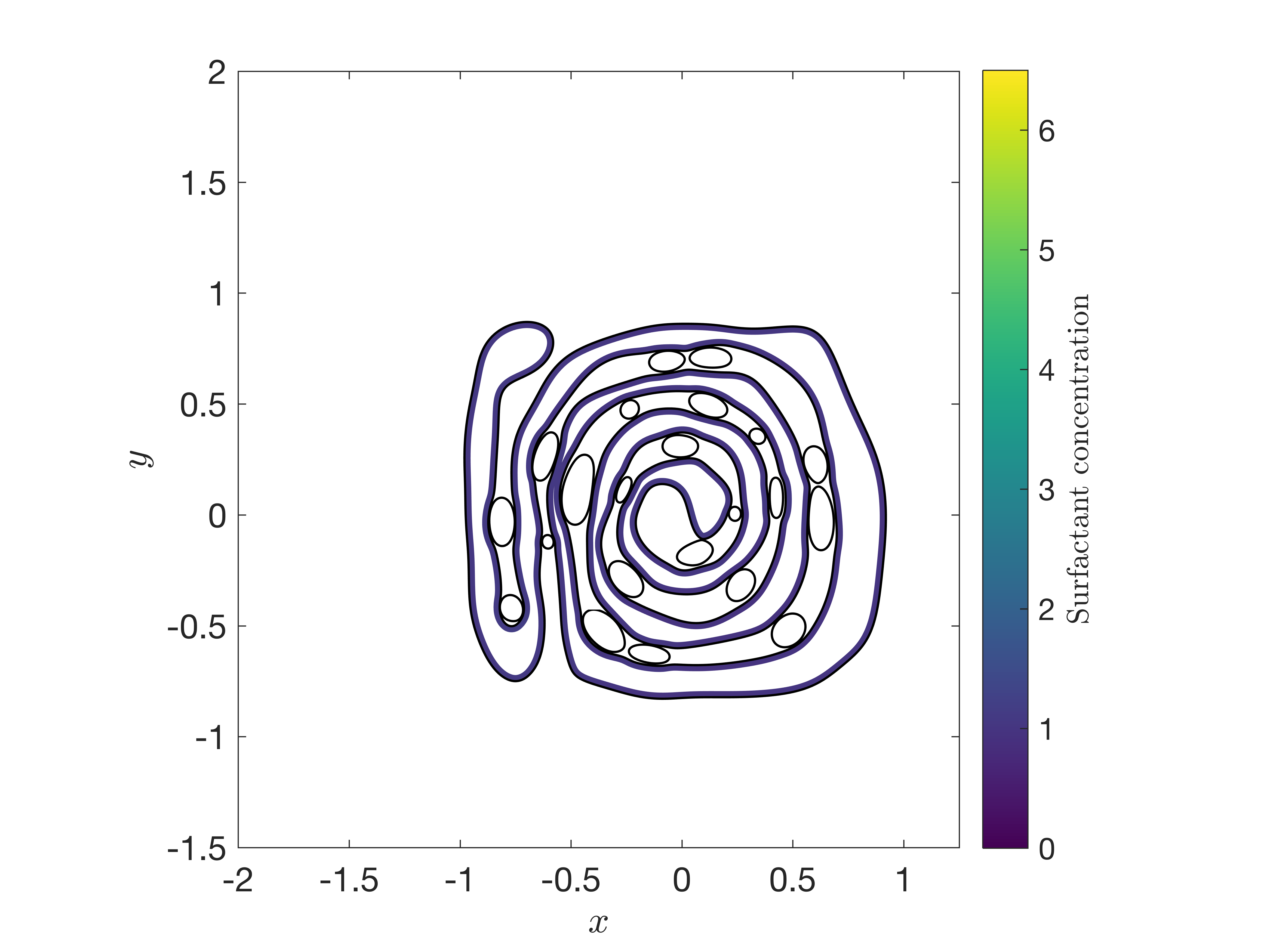}
\caption{Drop configuration for the swiss roll simulation.}
\label{fig:swissroll_domain}
\end{figure}

\noindent To further assess the BIE method, a more complicated simulation is set up to test the robustness of the method. The setup is inspired by that in \cite{Ojala2015}. The drop configuration is shown in Figure~\ref{fig:swissroll_domain}. No far-field flow is imposed, instead the drops will be allowed to deform until circular under surface tension. A drop is deemed circular when
\begin{align*}
  \left| 1 - \dfrac{\max{(|z-c|)}}{\text{mean}(|z-c|)}\right| < 10^{-4},
\end{align*}
where $c$ is the centre of each drop and $z$ its interface discretisation.

The largest drop (i.e. the roll) will be covered with surfactants with an initial uniform concentration $\rho_0=1$. Other parameters are $Pe=10$ and $E=0.1$. The authors note that the standard nondimensionalisation of length by initial radius is not appropriate in this case; instead the characteristic length of half of the length of a square box containing the initial drop configuration is used. Initially, the roll is covered with $375$ panels and the ellipses $37$ panels each. For this simulation $\lambda_k=1$ for all drops $k$ is considered, in order to make the simulation run in a reasonable time on a standard workstation. The simulation has also been run with 50 \% more points without any visible difference. Comparing the final circle radius and centre point between the two discretisations, they differ on the level of the circular tolerance imposed above. The simulation takes approximately $10^5$ time steps.

\begin{figure}%
\centering
\subfigure[][$t=2.5$]{%
\label{fig:ex3-a}%
\frame{\includegraphics[trim={4cm 2.75cm 5cm 2cm},clip,height=2in]{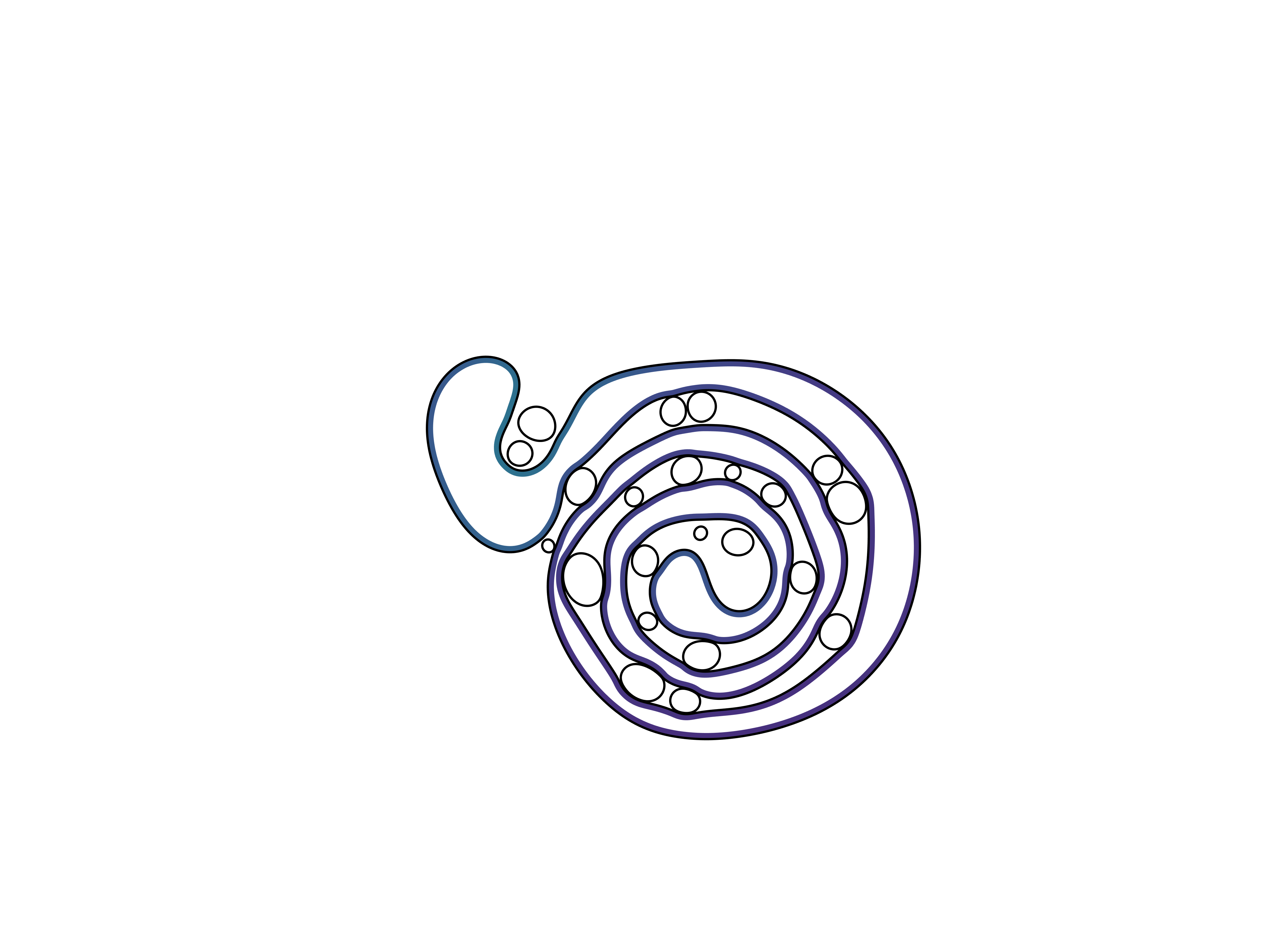}}}%
\hspace{0.5cm}
\subfigure[][$t=5$]{%
\label{fig:ex3-b}%
\frame{\includegraphics[trim={4cm 2.75cm 5cm 2cm},clip,height=2in]{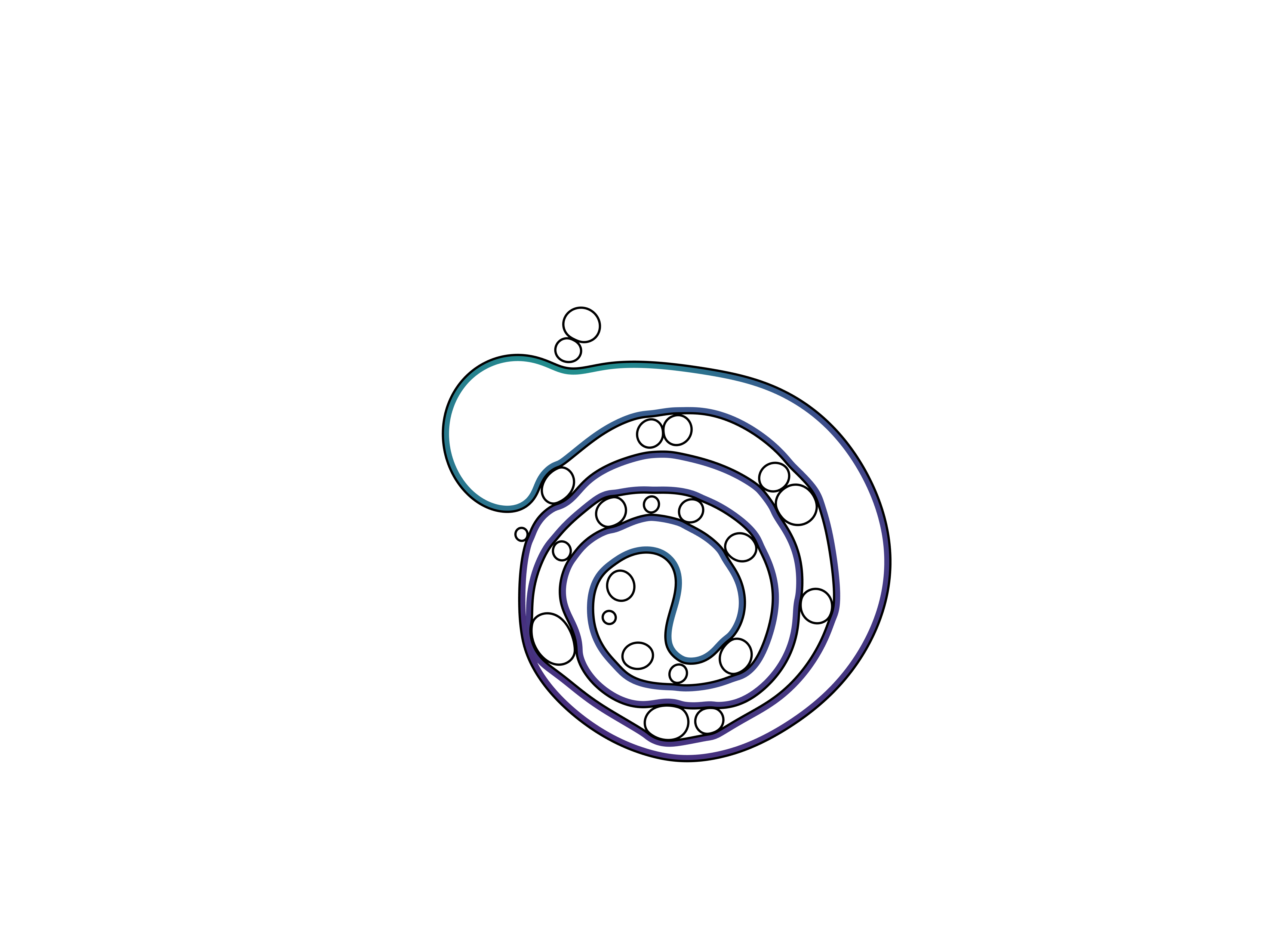}}} \\
\subfigure[][$t=10$]{%
\label{fig:ex3-c}%
\frame{\includegraphics[trim={4cm 2.75cm 5cm 2cm},clip,height=2in]{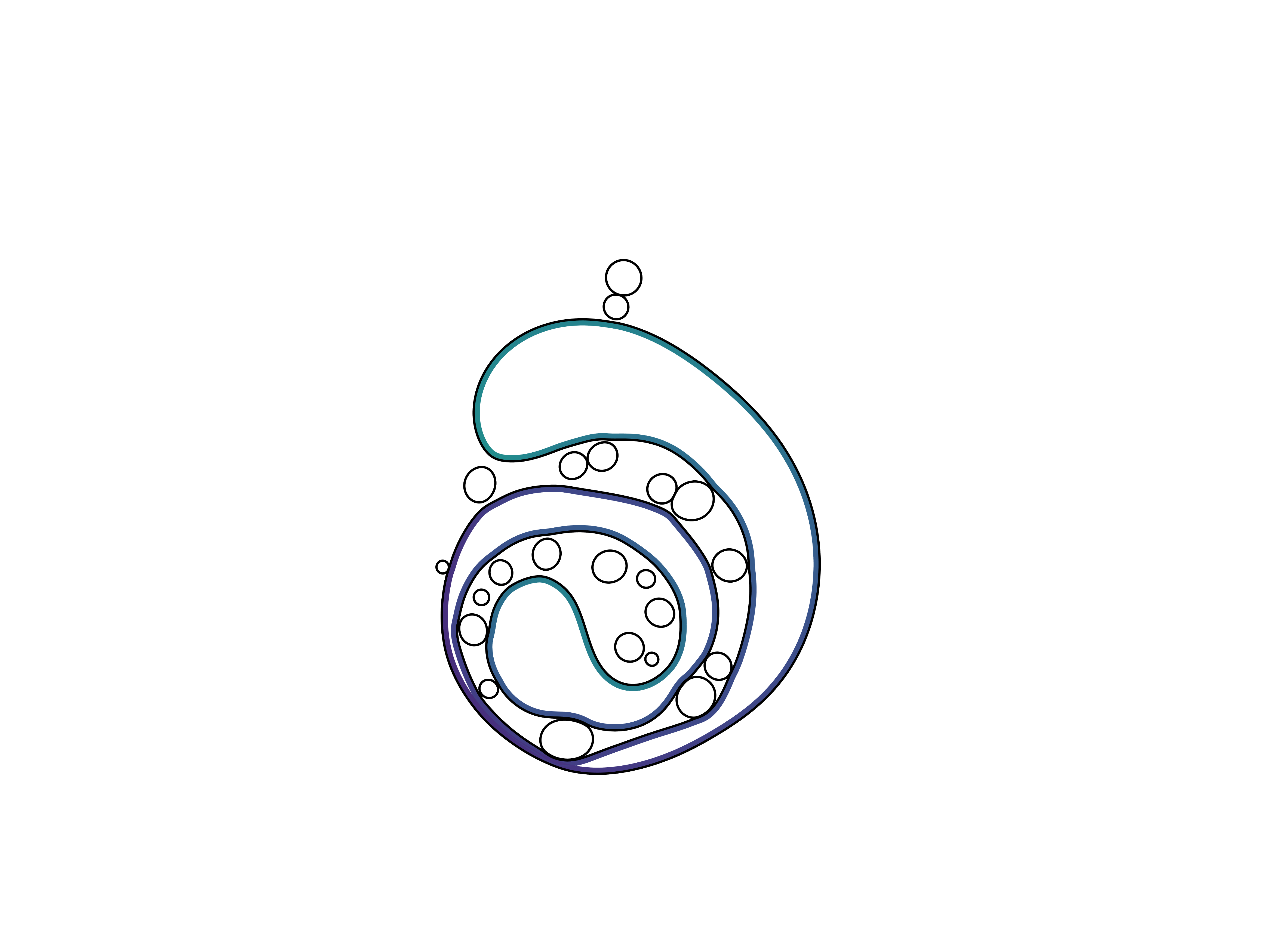}}}%
\hspace{0.5cm}
\subfigure[][$t=20$]{%
\label{fig:ex3-d}%
\frame{\includegraphics[trim={4cm 2.75cm 5cm 2cm},clip,height=2in]{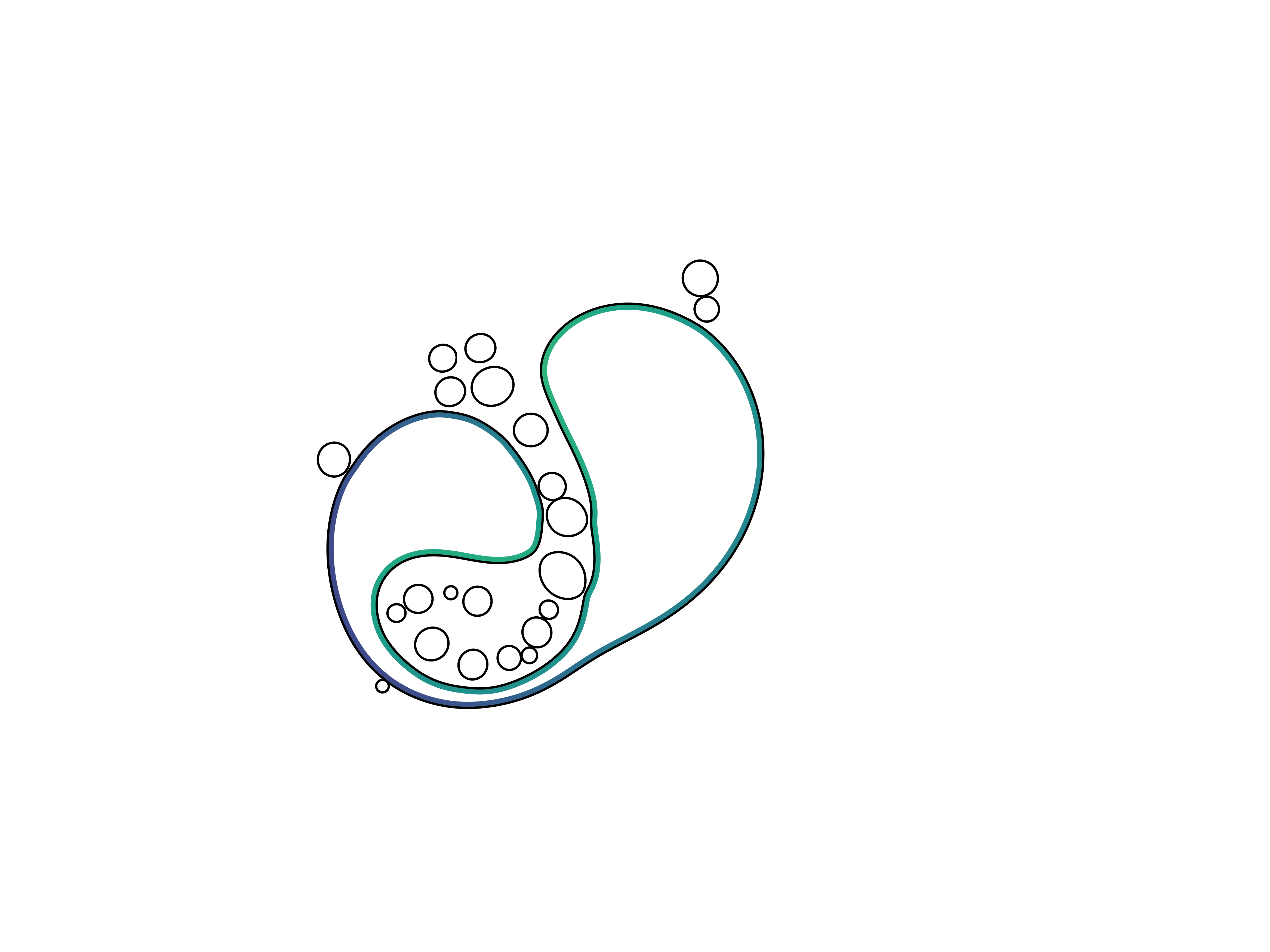}}} \\
\subfigure[][$t=30$]{%
\label{fig:ex3-e}%
\frame{\includegraphics[trim={4cm 2.75cm 5cm 2cm},clip,height=2in]{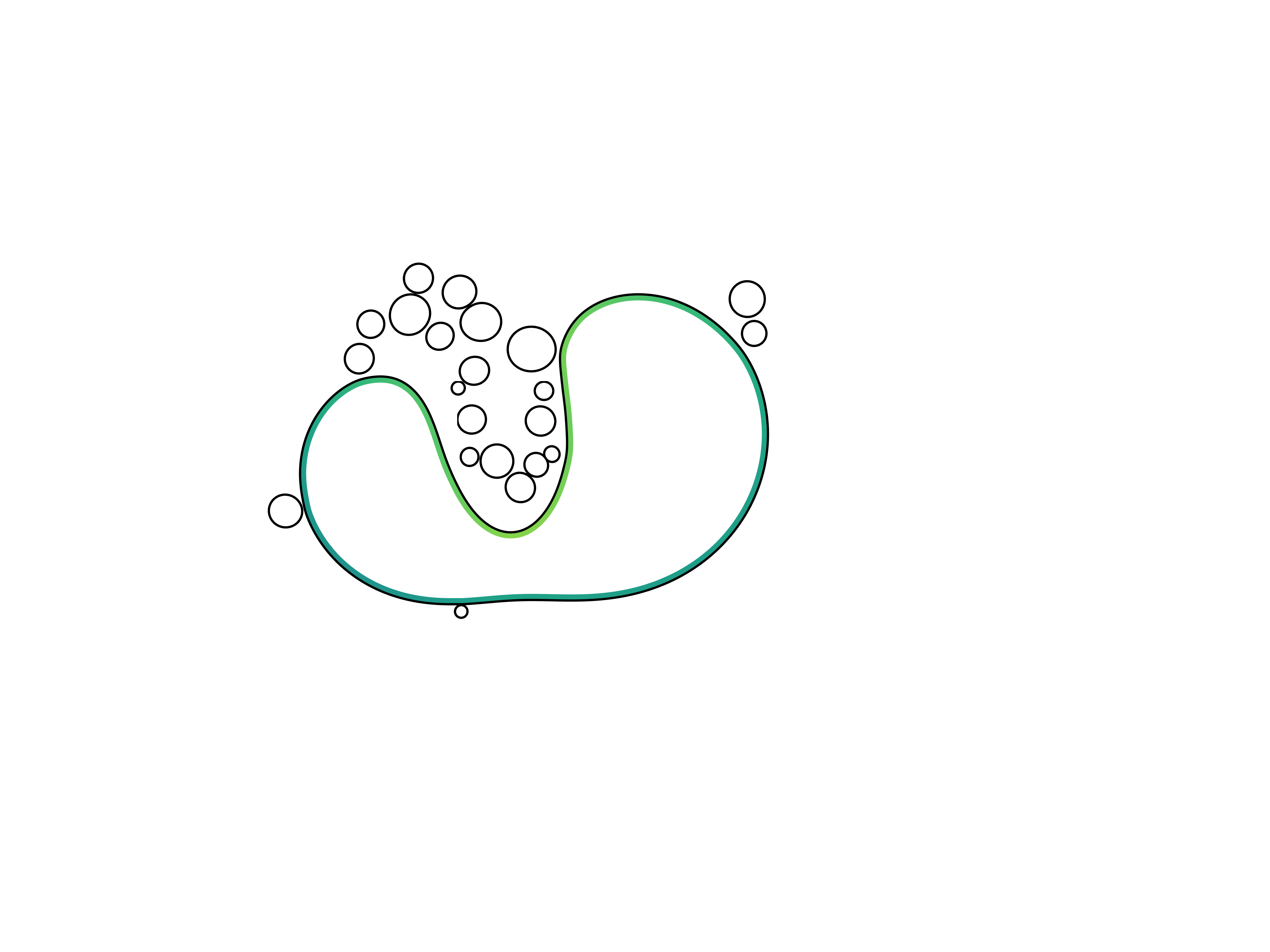}}}%
\hspace{0.5cm}
\subfigure[][$t=50$]{%
\label{fig:ex3-f}%
\frame{\includegraphics[trim={4cm 2.75cm 5cm 2cm},clip,height=2in]{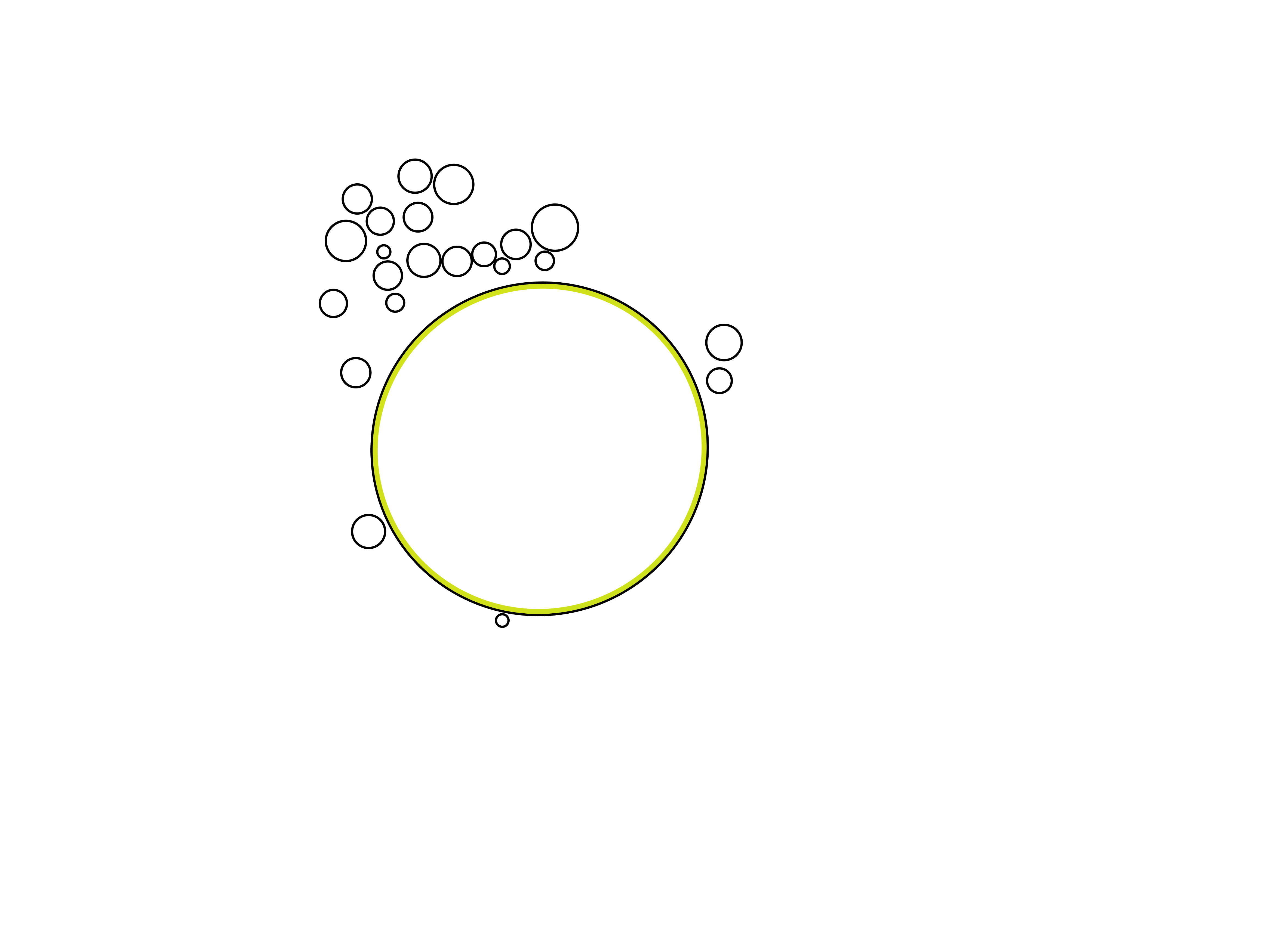}}}
\caption{Deformation over time for the swiss roll domain. Surfactant concentration in colour of the same scale as in Figure~\ref{fig:swissroll_domain}.
}%
\label{fig:swissroll_evolution}%
\end{figure}

How the drops deform in time is shown in Figure~\ref{fig:swissroll_evolution}. The drops reach their circular form at time $t=70$. The minimal distance between drops measured during the simulation is $9\cdot10^{-5}$. This can be compared to a corresponding simulation without surfactants, where the drops are circular at time $t=32$, with minimal distance $3\cdot10^{-4}$.

%% file: srcfiles/steady.tex

\subsection{Surfactant-covered bubble in steady state}
\begin{figure}[h!]
	\centering
	\includegraphics[width=0.445\textwidth]{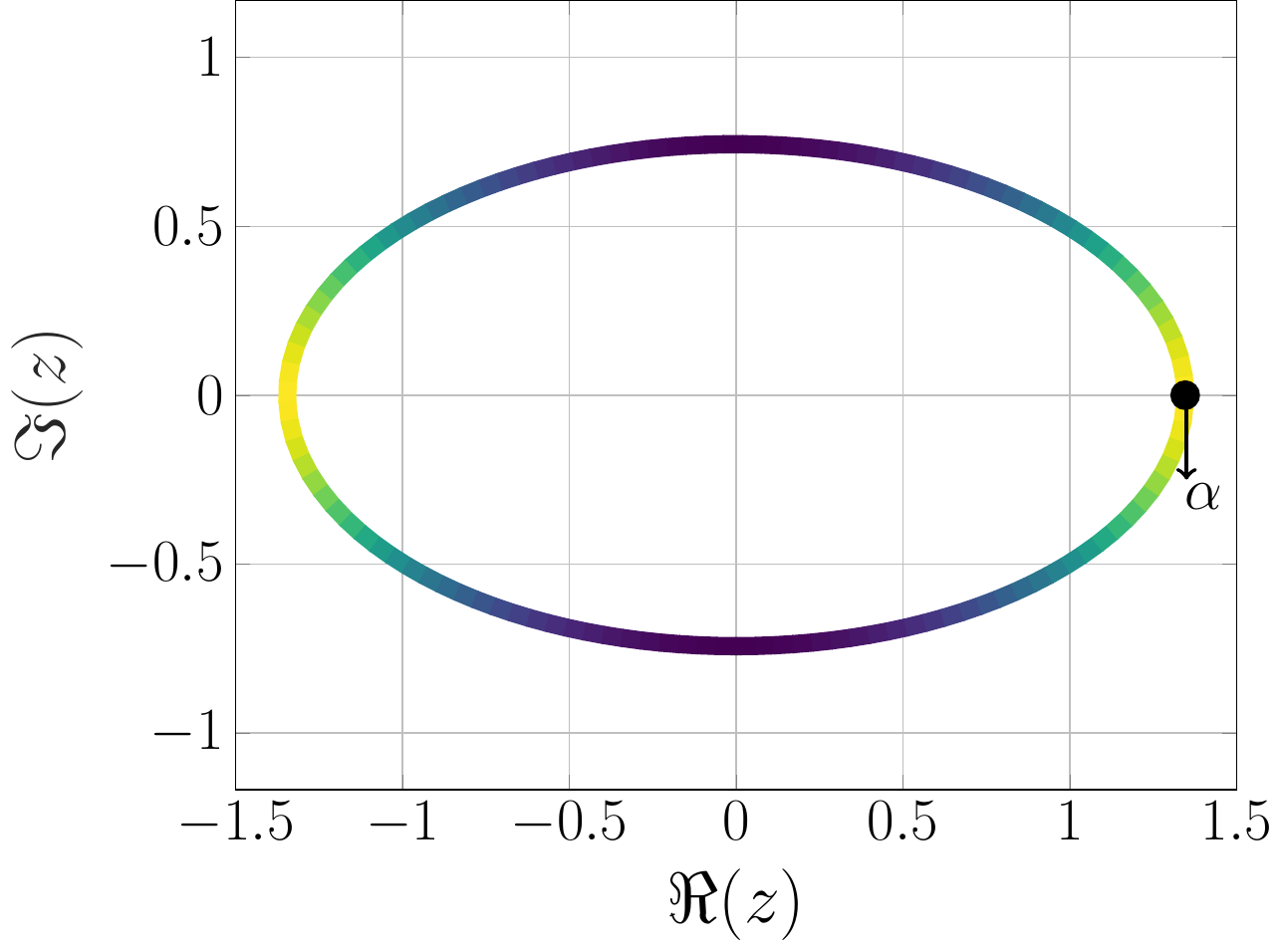}%
	\includegraphics[width=0.5\textwidth]{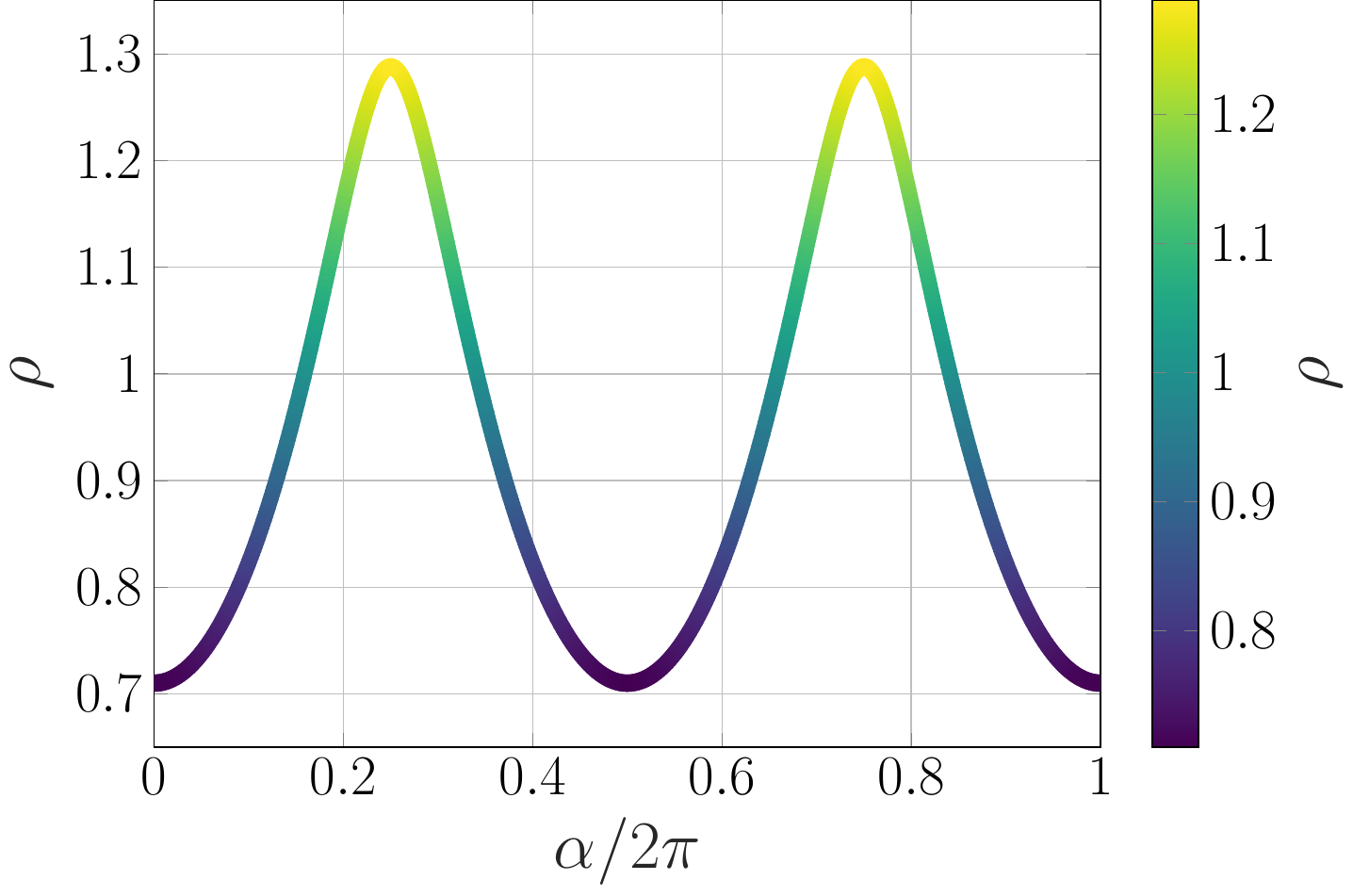}%
    \caption{Interface (left) and surfactant concentration (right) at steady state for Capillary number $Q=0.14$. Elasticity number $E=0.5$, initial surfactant concentration $\rho_0=1$, $Pe_\Gamma=\infty$ and viscosity ratio $\lambda=0$. Steady state was reached at time $T=46.35$.}
    \label{fig:plot12}
\end{figure}

\noindent Using the exact solutions of \S\ref{sec:validsteady}, the bubble deformation and surfactant concentration at steady state is tested. Given a Capillary number $Q$, the steady state interface position $z(\zeta,T)$ and surfactant concentration $\rho(\zeta,T)$ are known, where $T$ is the time of steady state and $\zeta=e^{i\nu}$ for $\nu\in[0,2\pi)$. An example of how such a steady state may look like can be seen in Figure~\ref{fig:plot12}, where the steady state corresponding to $Q=0.14$ is shown.
In a simulation, a bubble is said to have reached steady state when the normal velocity, $\mathbf{u}\cdot\mathbf{n}=0$. Practically, here the definition of steady state is that
\begin{align*}
	\left|\mathbf{u}\cdot\mathbf{n}\right|\leq 10^{-8}.
\end{align*}
It should be noted that for large Capillary numbers, $Q$, or high viscosity ratios, $\lambda$, this threshold is slightly too large and should be decreased since the deformation is very slow.

The analytical solutions from Algorithm~\ref{alg:steady} are denoted by $z^{V}(\zeta,T)$ and $\rho^{V}(\zeta,T)$ for interface position and surfactant concentration respectively. Defining the obtained solutions from the numerical method as $z(\alpha,T)$ and $\rho(\alpha,T)$, a comparison between the numerical results and the validation is obtained by studying the point-wise differences in $\alpha\in[0,2\pi)$. Two errors will be studied: in position, $e_z = |z-z^V|$ and surfactant concentration $\rho$, $e_\rho$. In order to compare point-wise in $\alpha\in[0,2\pi)$, a mapping from $\nu$ to the equal-arc length measure $\alpha^V$ is needed. This is defined as
\begin{align}
	\alpha^V(\nu) = \dfrac{S_\nu 2\pi}{L},
	\label{eq:alphaV}
\end{align}
where $S_\nu = \int_0^\nu |z^V_\nu|d\nu$. The solutions $z^V$ and $\rho^V$ are computed at non equidistant discrete points $\alpha_i^V=\alpha^V(\nu_i)$, for $\nu_i\in[0,2\pi)$, $i=1,\hdots,N_V$ for some $N_V$. Finally, to compare the results, $z(\alpha,T)$ and $\rho(\alpha,T)$ need to be interpolated to $\alpha^V_i$. This is done using their Fourier expansions.

\begin{figure}[h!]
	\centering
	\includegraphics[width=0.5\textwidth]{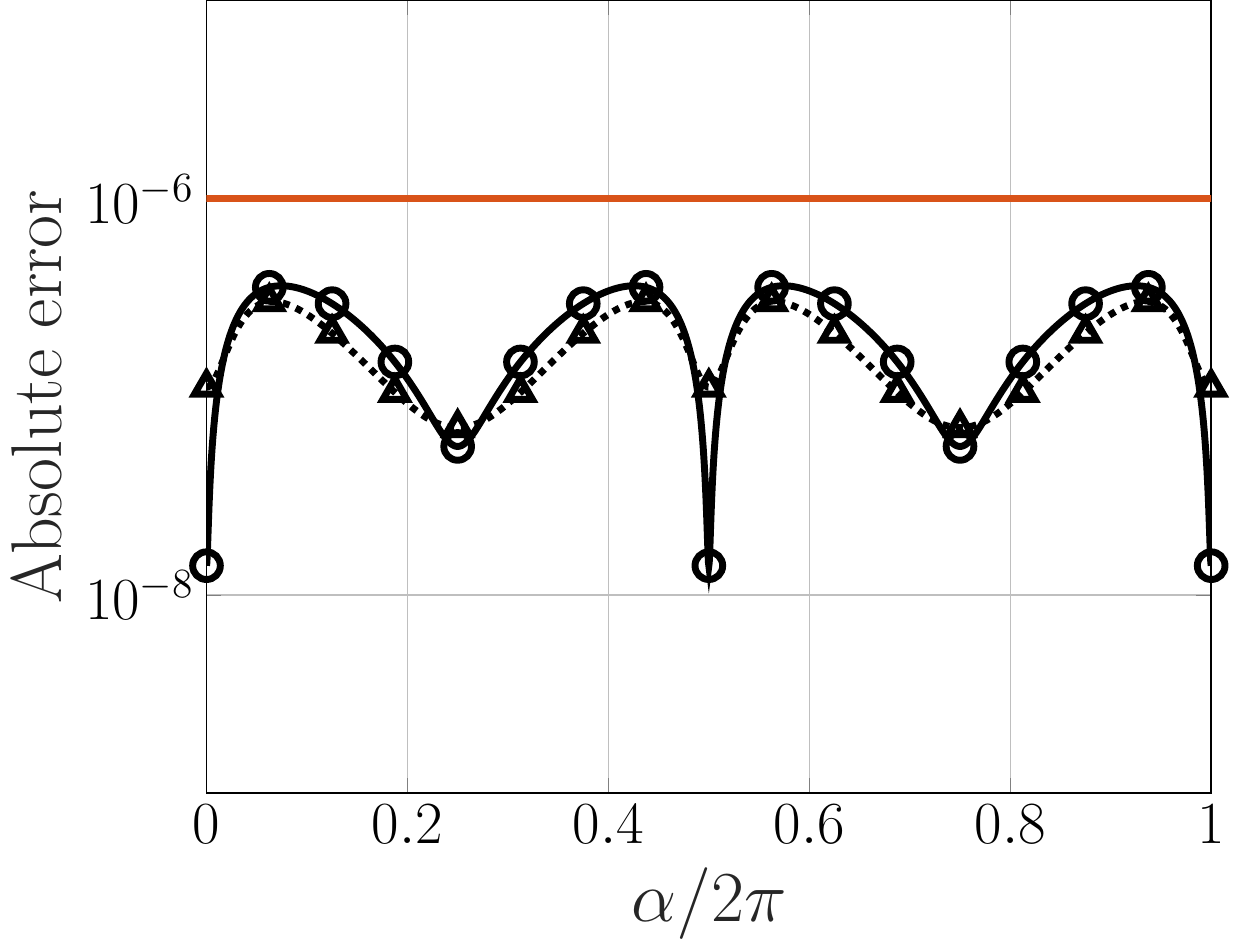}%
	\caption{Point-wise error in $z$ (marker $\circ$) and $\rho$ (marker $\triangle$) vs. $\alpha$ at steady state, for the same parameters as in Figure~\ref{fig:plot12}, grid spacing $\Delta s = 0.0016$. Red line corresponds to time-stepping tolerance set to $10^{-6}$.}
	\label{fig:plot3}
\end{figure}

In Figure~\ref{fig:plot3}, the point-wise absolute error between analytical solution and the BIE simulation can be seen for the case in Figure~\ref{fig:plot12}. The spatial discretisation is adaptive, keeping the spatial distance in arc length similar at all times; here approximately $\Delta s=0.008$. The simulation was run with time-step tolerance $10^{-6}$ and the errors are all below the set tolerance.

\begin{figure}[h!]
	\centering
	\includegraphics[width=0.6\textwidth]{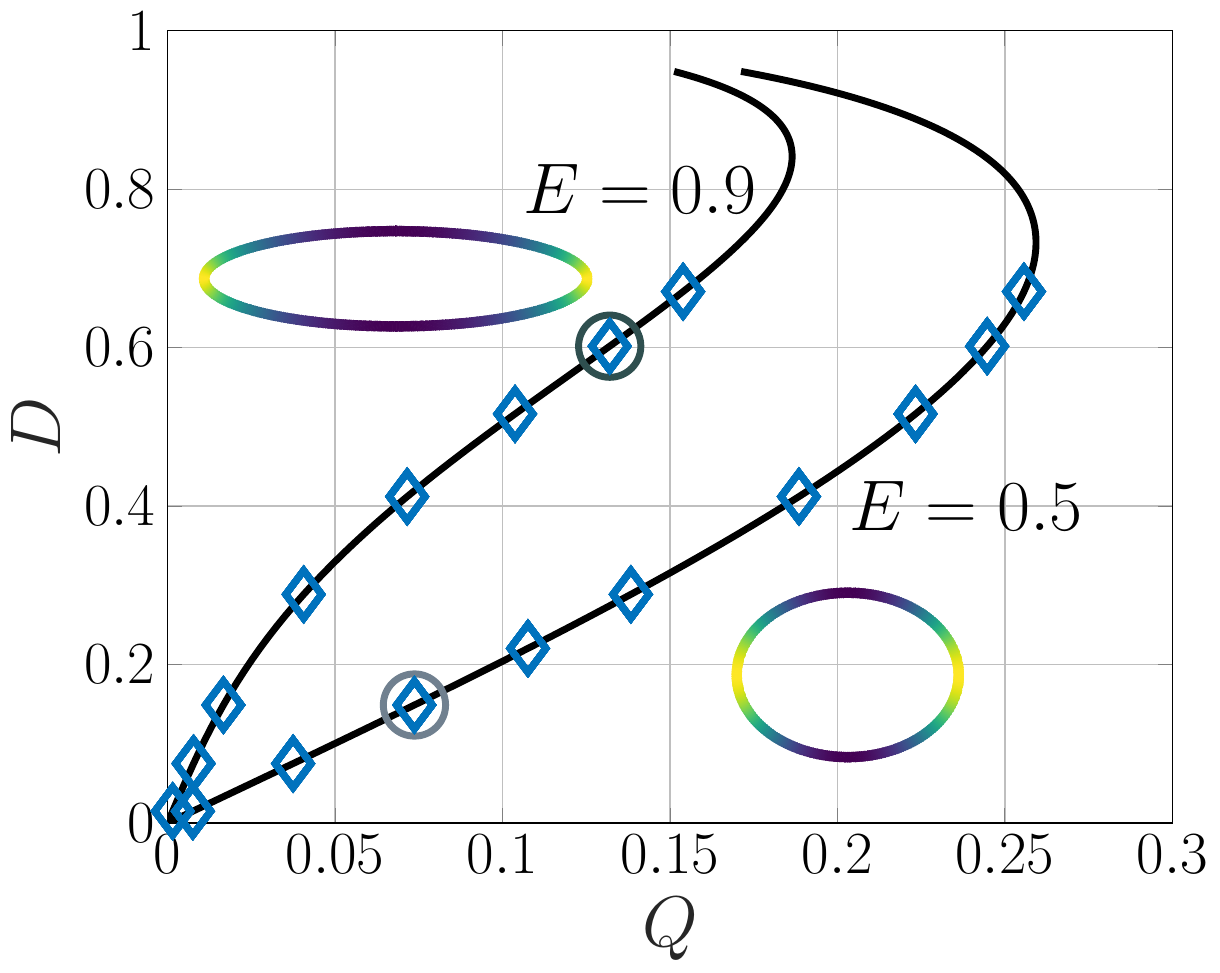}%
	\caption{Steady state deformation, $D$, vs. Capillary number $Q$ at different elasticity numbers $E$ for inviscid bubbles, $\lambda=0$. Black lines: validation solutions computed through the method in \S\ref{sec:validsteady}. Diamonds: numerical results from the BIE method. Circles show simulations chosen for the examples of deformed bubbles shown.}
	\label{fig:plot4}
\end{figure}

Using the method in \S\ref{sec:validsteady}, one can obtain a graph for how deformation $D$ depends on Capillary number $Q$. Such a graph for two different elasticity numbers, $E=0.5$ and $E=0.9$, is shown in Figure~\ref{fig:plot4}. The black lines represent the validation solutions computed through the method in \S\ref{sec:validsteady} and the blue diamonds show results from simulations using the BIE method of this paper. All simulations have been run from the same initial setting as the previous case, with time-step tolerance $10^{-6}$. The simulations differ only in the choice of Capillary number $Q$. Examples of bubble deformation and surfactant concentration are also shown, corresponding to $Q=0.07$ and $Q=0.13$, shown by circles in the figure.

\begin{figure}[h!]
	\centering
	\includegraphics[width=0.45\textwidth]{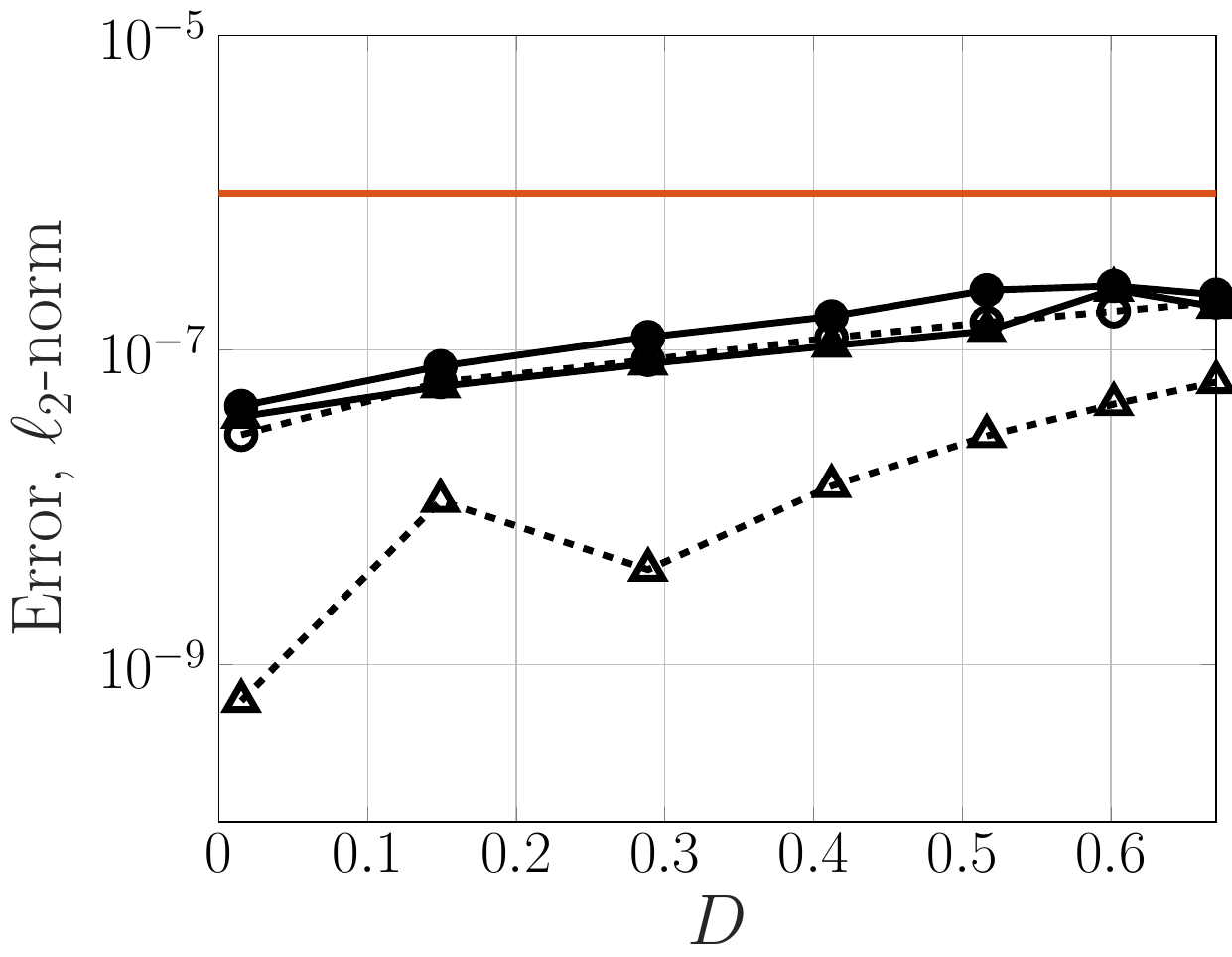}%
	\hspace{1mm}
	\includegraphics[width=0.45\textwidth]{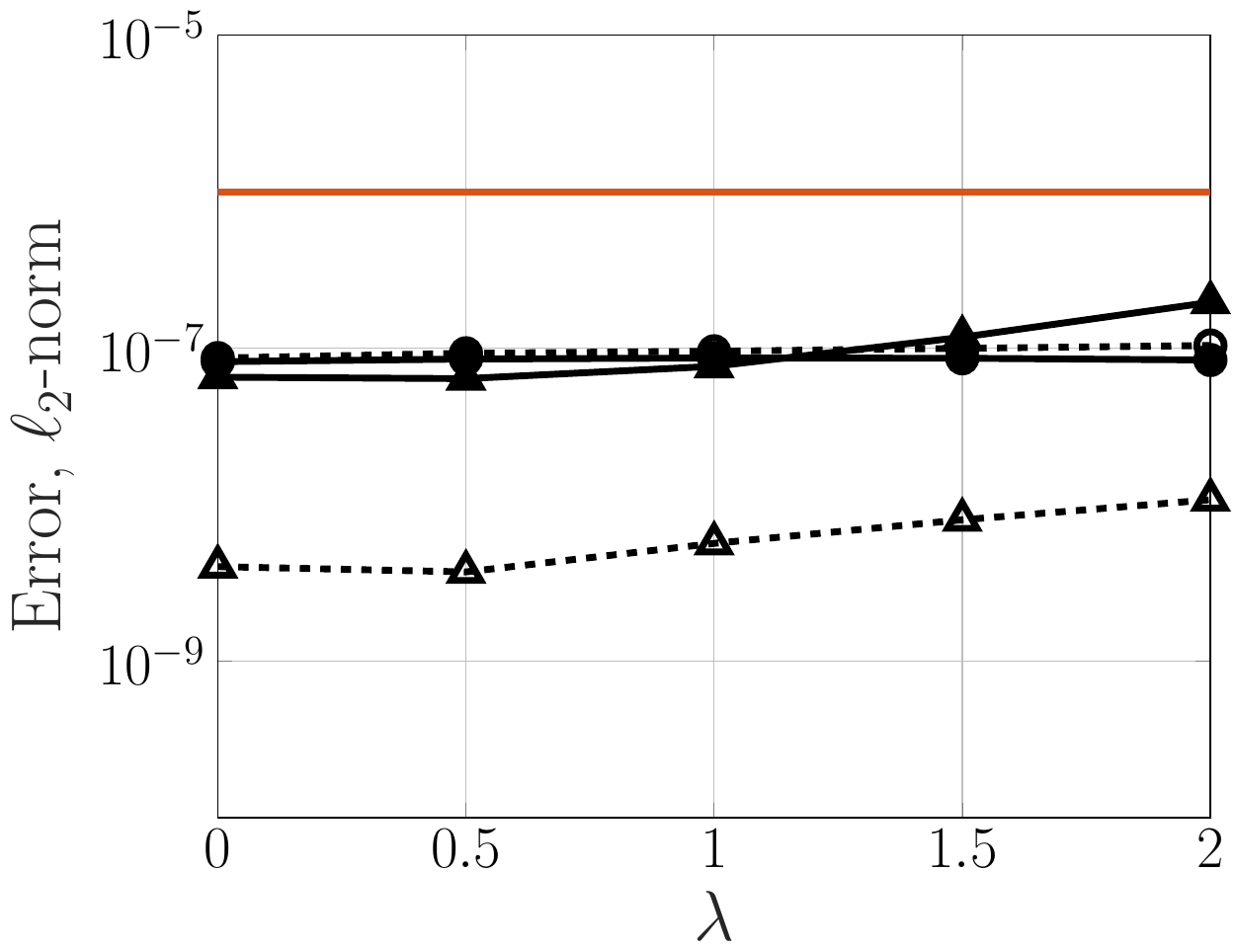}%
	\caption{$\ell_2$ error in $z$ and $\rho$ vs. (left) deformation $D$ for cases in Figure~\ref{fig:plot4} and (right) vs. viscosity ratio $\lambda$ for deformation $D=0.29$. Solid and dashed lines represent $E=0.5$ and $E=0.9$ respectively. Dots and triangles represent $e_z$ and $e_\rho$ respectively. Red solid line shows time-step tolerance $10^{-6}$.}
	\label{fig:plot56}
\end{figure}

To quantify the results in Figure~\ref{fig:plot4}, the $\ell_2$-error for each deformation $D$ is plotted in Figure~\ref{fig:plot56} (left). For each diamond in Figure~\ref{fig:plot4}, the simulation and validation data is compared and plotted vs. deformation $D$. Note that both cases $E=0.5$ and $E=0.9$ are shown (solid and dashed lines respectively), and that they for the same $D$ were obtained by different $Q$. It is shown that both errors ($e_z$ and $e_\rho$) stay below time-step tolerance.

As mentioned in \S\ref{sec:validsteady}, the viscosity ratio does not affect steady state deformation in the case of $Pe_\Gamma=\infty$. The results of this paper agree with this observation, as is shown in Figure~\ref{fig:plot56} (right). There, droplets of different viscosity ratios ranging from $\lambda = 0$ (corresponding to an inviscid bubble) to $\lambda = 2$ were deformed under Capillary number $Q=0.14$. As can be seen in the figure, both interface position and surfactant concentration coincides with that for bubbles up to time-step tolerance.

%% file: srcfiles/cleanvalid.tex

\subsection{A pair of clean bubbles in extensional flow}
\noindent Using the approach described in \S\ref{sec:cleanpair}, the semi-analytical solutions for a pair of bubbles deforming in an extensional flow are computed  and compared against the boundary integral method of this paper.

\begin{figure}[h!]
	\centering
	\includegraphics[width=0.45\textwidth]{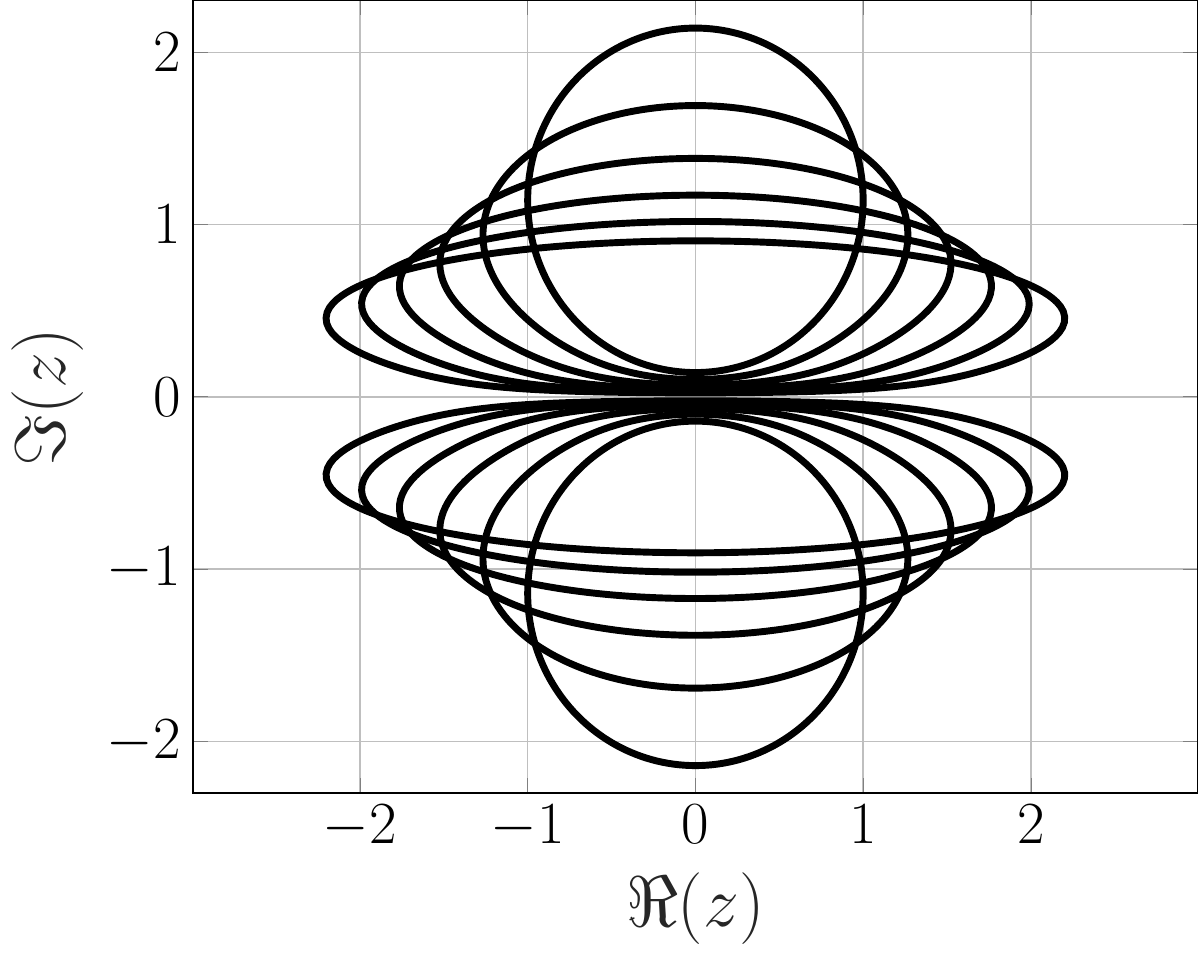}%
	\includegraphics[width=0.45\textwidth]{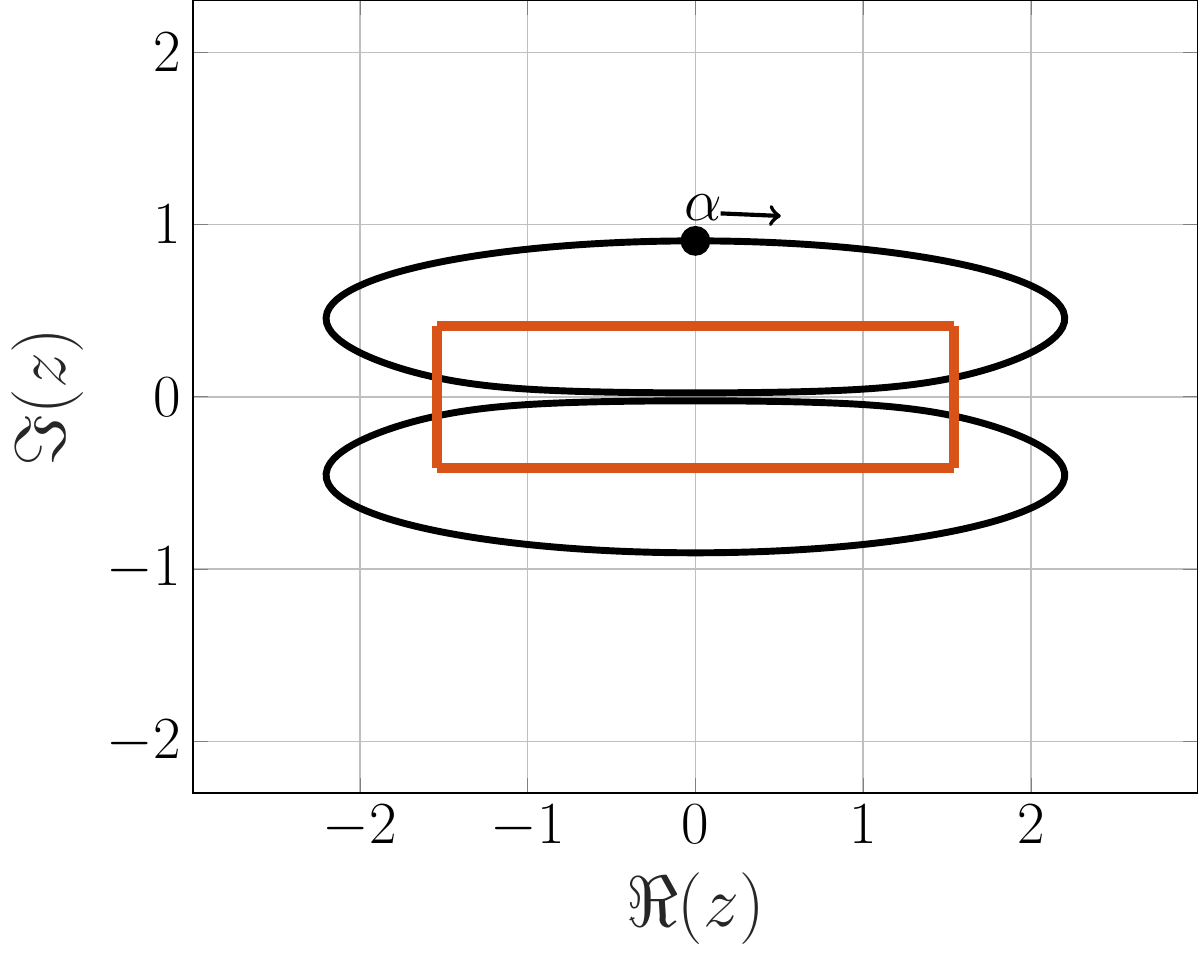}%
	\caption{Clean bubbles deforming in an extensional flow. Left: interface position from time $t=0$ to $t=1.5$, at $dt=0.3$ intervals, computed with BIE method. Right: final interface configuration at time $t=1.5$. Red box: cut-off region for method comparisons.}
	\label{fig:crowdy_z}
\end{figure}

The following case has been selected as it pushes the bubbles close to each other, thus providing an excellent test case for the special quadrature. The bubbles are initially circular with radius one, centred around $\pm 1.419i$ which corresponds to an initial $\phi(0)=0.35$. The bubbles are clean, i.e. there are no surfactants present in this problem, and they deform under an extensional flow with Capillary number $Q=0.5$ until time $t=1.5$. At the final time, the minimum distance between the bubbles is $0.04$. In Figure~\ref{fig:crowdy_z} the movement of the bubbles over time  and final deformation is shown (left and right respectively).

\begin{figure}[h!]
	\centering
	\includegraphics[width=0.45\textwidth]{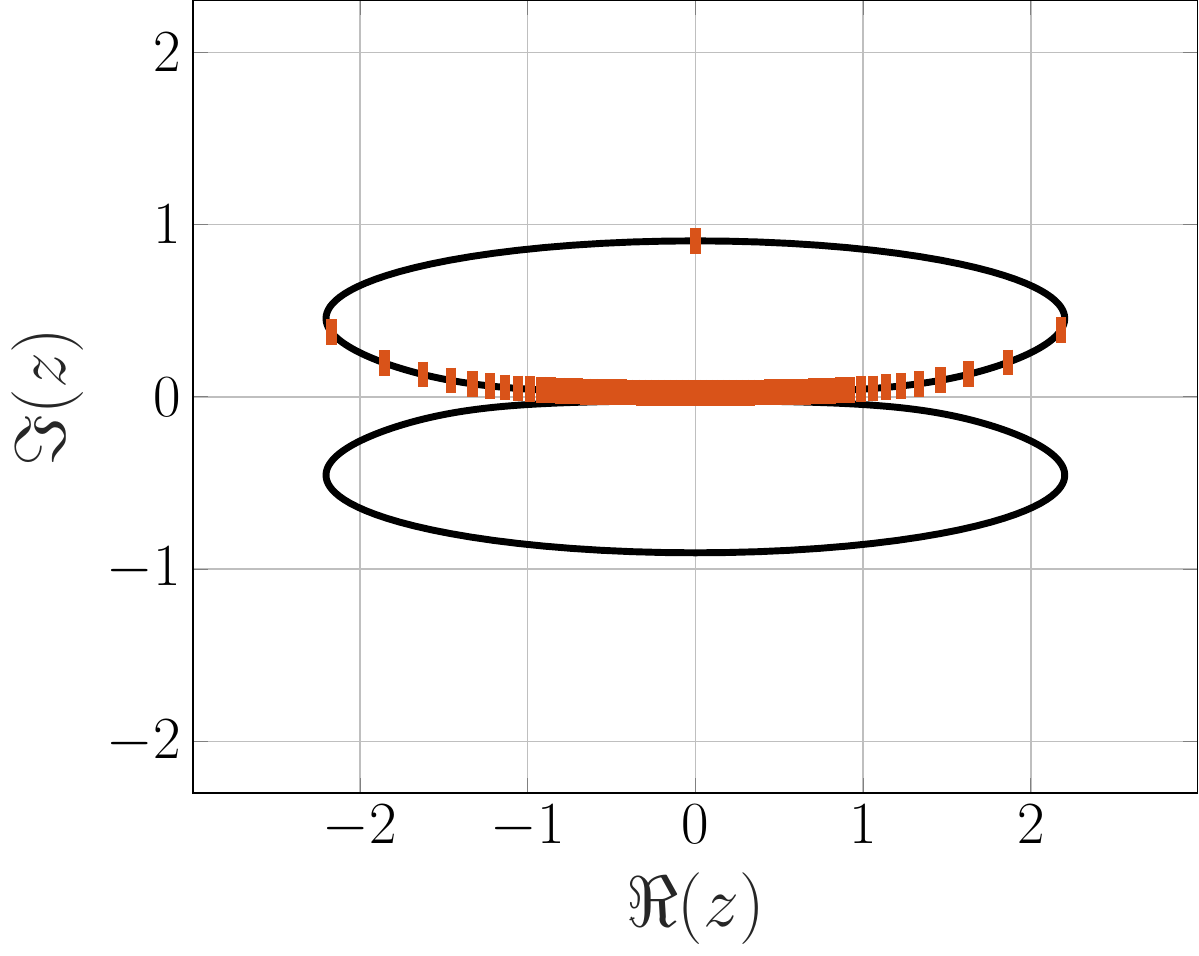}%
	\includegraphics[width=0.45\textwidth]{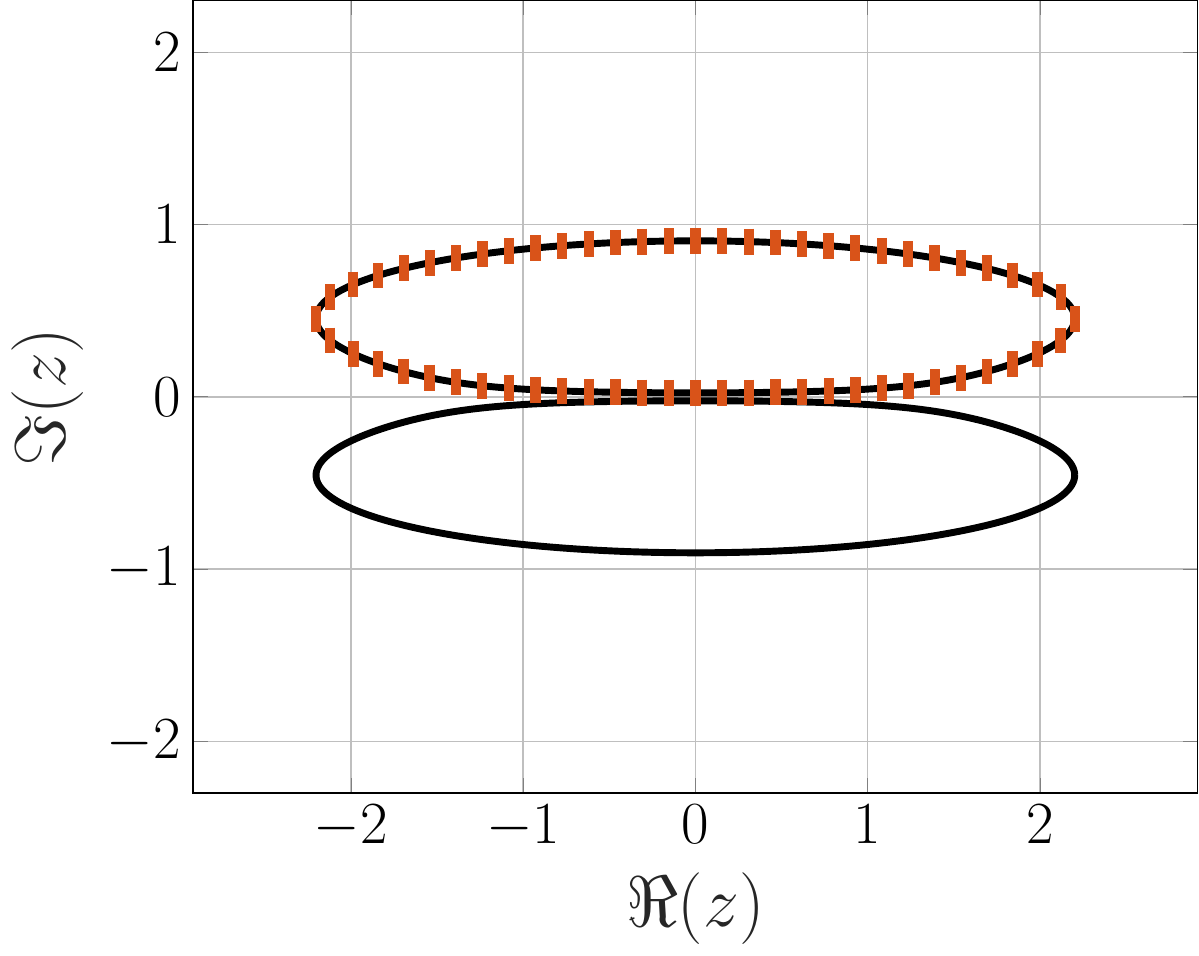}%
	\caption{Interface at time $t=1.5$ for clean bubbles. Red markers: every 16th point of validation discretisation for $N_V=1025$ (left) and every 16th point of uniform discretisation for BIE method, $N_k=960$, $k=1,2$. (right). }
	\label{fig:crowdy_bubble_PR1}
\end{figure}

The interface position obtained from the BIE method is denoted $z(\alpha,t)$, where $\alpha\in[0,2\pi)$. In Figure~\ref{fig:crowdy_bubble_PR1} (right) the discretisation along the interface is shown. Note that it is uniform around the interface. In contrast, the discretisation along the interface for the validation method is not uniform, see Figure~\ref{fig:crowdy_bubble_PR1} (left). There, the points are clustered where the distance between the bubbles is the smallest. When comparing the results of the two methods, only the part of the interfaces between $\alpha=\frac{2\pi}{3}$ and $\alpha=\frac{4\pi}{4}$ will be considered. This corresponds to the red box in Figure~\ref{fig:crowdy_z} (right). This cut-off is made because the discretisation of the validation method is very coarse outside of it, with large spatial errors. Also, since it is the area of near-interaction between the bubbles which is of greatest interest, as it is there the boundary integral method will potentially struggle the most, it is sufficient to consider this part.

To investigate how the two methods behave, they are compared against a reference solution, $z_{ref} = x_{ref}+iy_{ref}$. This solution is computed with the boundary integral method using $600$ Gauss-Legendre panels, i.e. $N_k=9600$, $k=1,2$, and a time-step tolerance of $tol=10^{-10}$. To compare the non-uniform discretisation of the validation method with the uniform discretisation of the BIE method, the map as in \eqref{eq:alphaV} is employed and the points are interpolated using a non-uniform \texttt{FFT}. Errors in both $x$- and $y$-coordinates will be considered and denoted $e_x=|x-x_{ref}|$ and $e_y=|y-y_{ref}|$ respectively, the combined error in $z$ is denoted $e_z = |e_x + ie_y|$.

\begin{figure}[h!]
	\centering
	\includegraphics[width=0.42\textwidth]{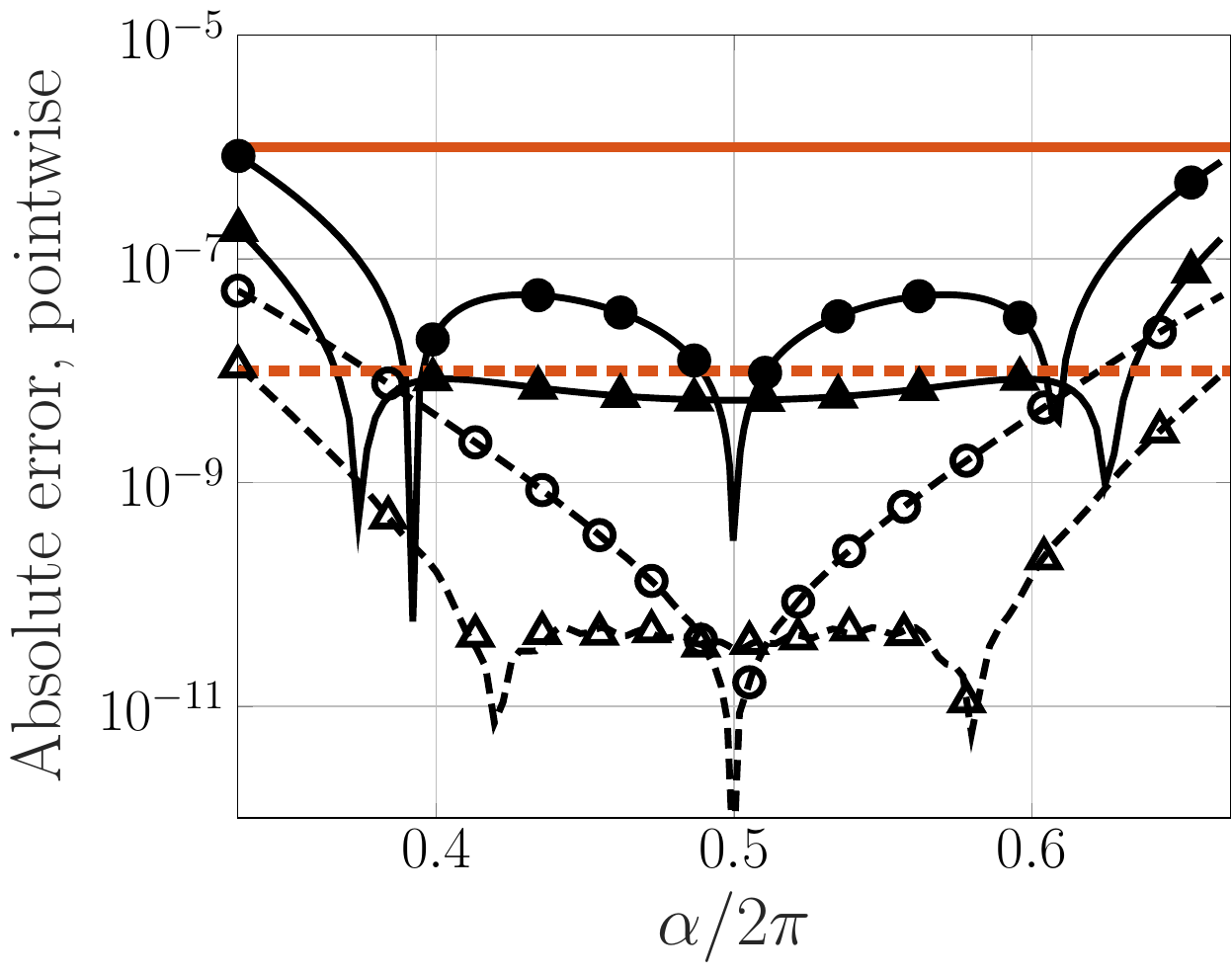}%
	\hspace{1.5mm}
	\includegraphics[width=0.42\textwidth]{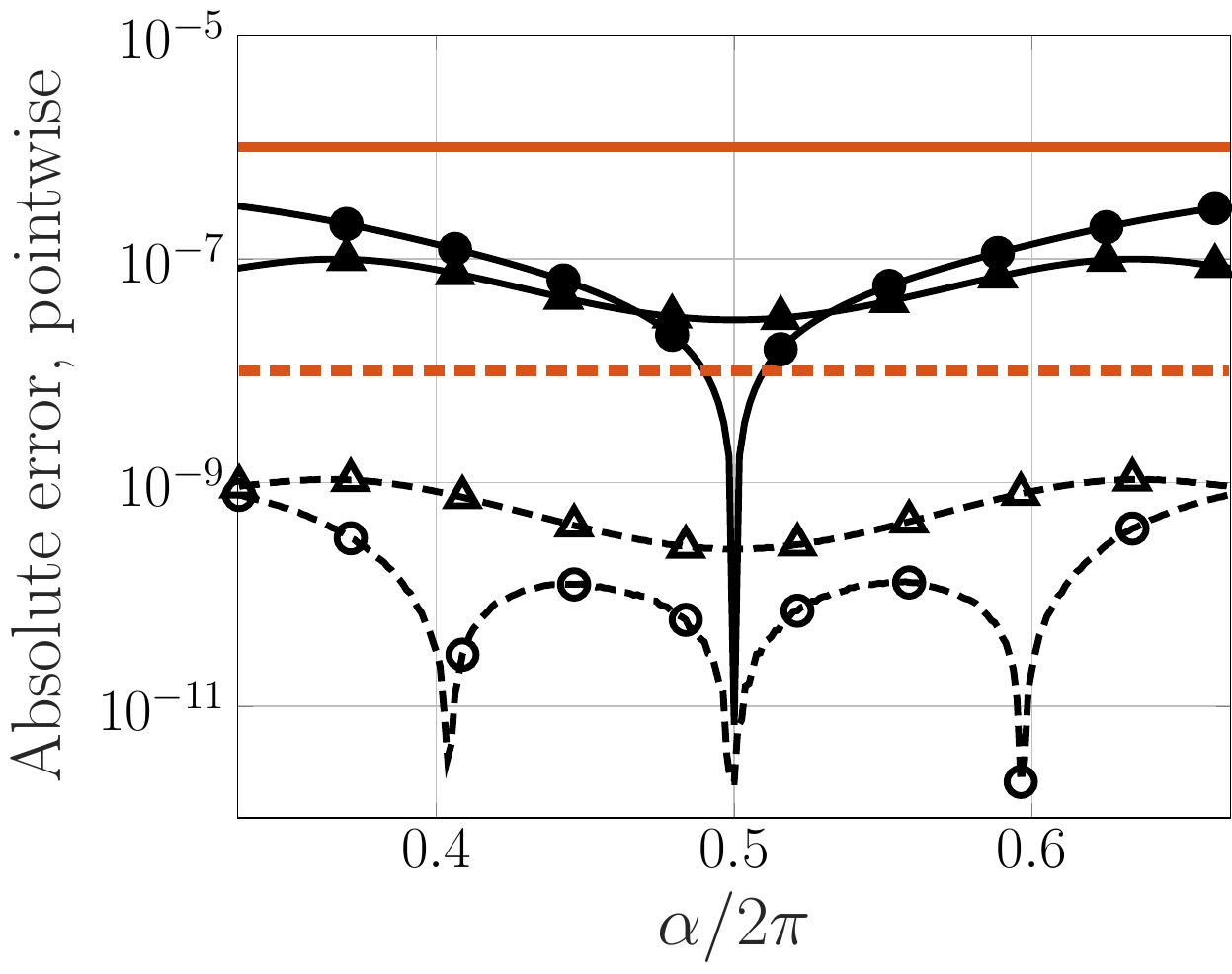}%
	\caption{Pointwise relative difference in $x$- and $y$-coordinates for validation method (left) and BIE method (right) against reference solution. Solid lines/markers shows time-stepping tolerance $10^{-6}$, dashed for tolerance $10^{-8}$. BIE method computed with $576$ and $800$ points per bubble for $tol_1$ and $tol_2$ respectively. Validation method computed with $2049$ and $9001$ points respectively. Error in $x$ denoted by $\circ$ and in $y$ by $\vartriangle$.}
	\label{fig:crowdy_differences_small}
\end{figure}

In Figure~\ref{fig:crowdy_differences_small} the point-wise absolute errors are shown in the region of interest for the validation method and the BIE method. The results have been computed using two time-step tolerances: $tol_1 = 10^{-6}$ (solid lines) and $tol_2 = 10^{-8}$ (dashed lines). For the validation method, $N_V~=~2049$ and $N_V=9001$ points were used for $tol_1$ and $tol_2$ respectively, and for the BIE method $N_k=576$, $k=1,2$, was used for $tol_1$ and $N_k=800$ for $tol_2$. It is shown that both methods have errors below the set tolerance $tol_1$, but that the validation method does not quite achieve the tolerance $10^{-8}$. This is due to the increase in $\Delta\alpha$ towards the end of the cut-off, which causes some of the points to have a larger error than the set tolerance.

\begin{figure}[h!]
	\centering
	\includegraphics[width=0.45\textwidth]{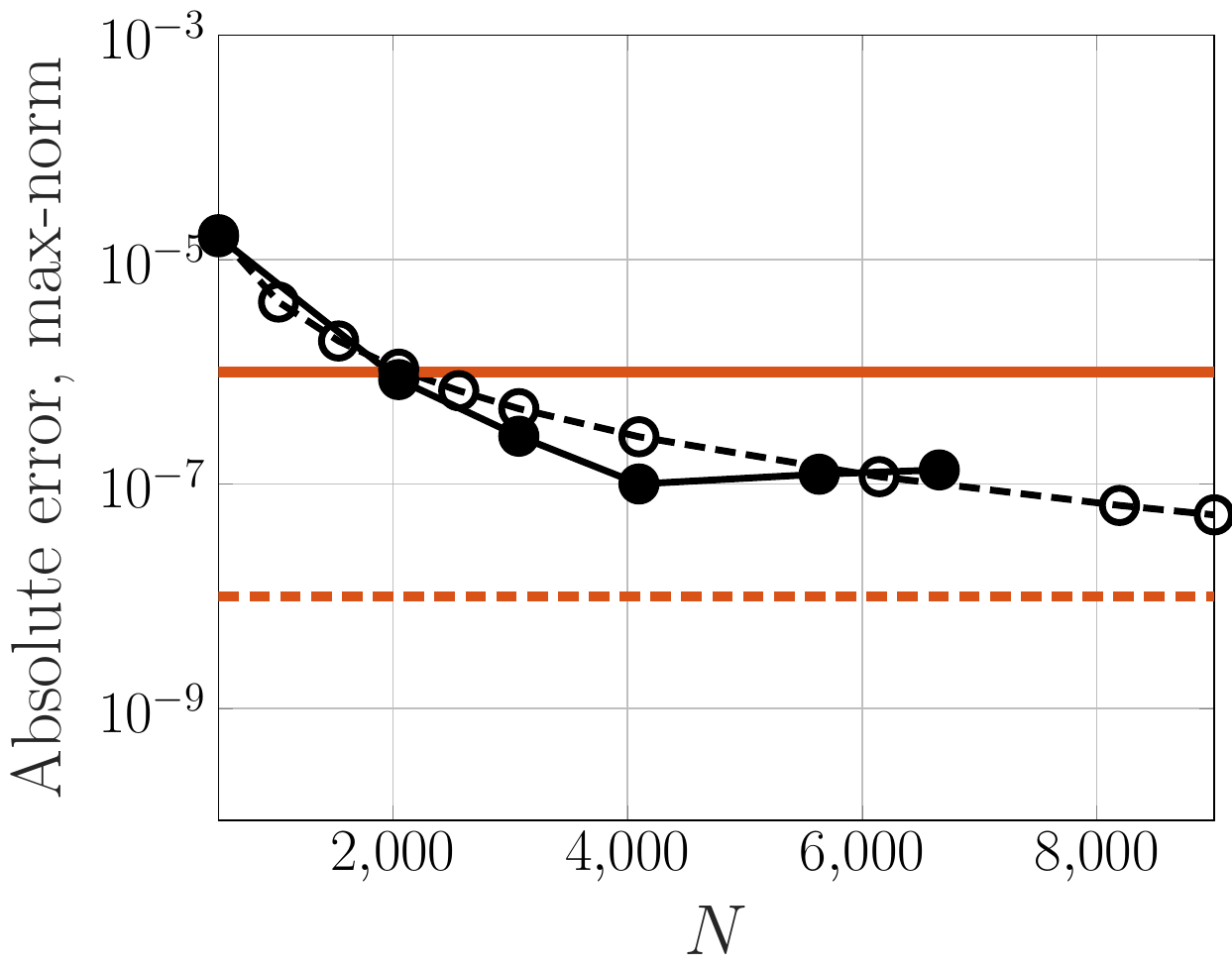}%
	\hspace{1.5mm}
	\includegraphics[width=0.465\textwidth]{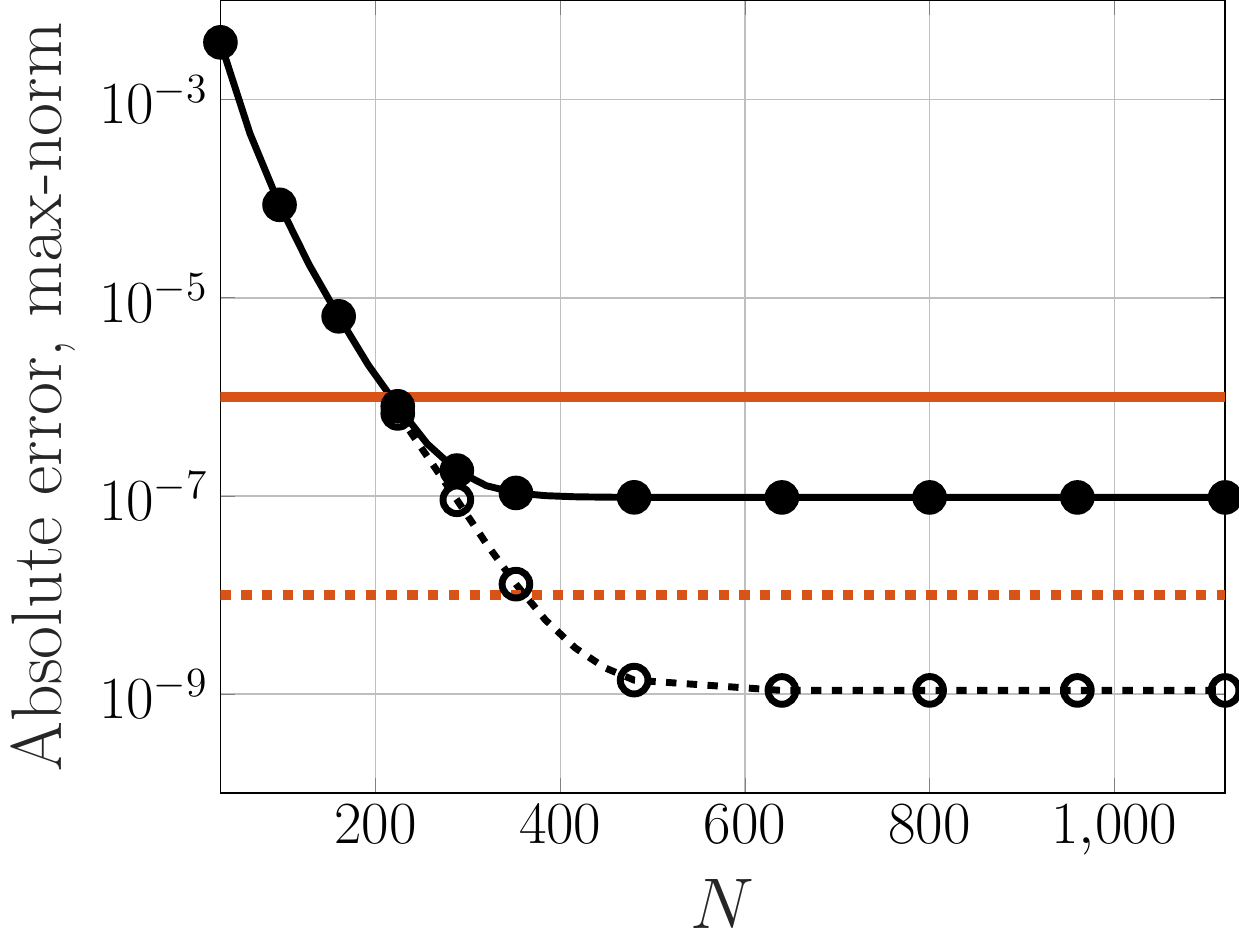}%
	\caption{Error as a function of number of discretisation points for validation method (left) and BIE method (right). Solid lines/markers show results for time-stepping tolerance $10^{-6}$ and dashed lines tolerance $10^{-8}$.  The spatial adaptivity of the BIE method was turned off to keep $N$ constant throughout the simulation.}
	\label{fig:crowdy_vsN}
\end{figure}

To quantify the behaviour of the two methods, the absolute errors in max-norm for different discretisations of the bubbles are plotted. Note that the number of discretisation points is constant for both methods, i.e. it does not increase over time as the bubbles get more stretched out. This especially means that the spatial adaptivity of the BIE method was turned off throughout the simulation. In Figure~\ref{fig:crowdy_vsN} the errors vs number of discretisation points $N$ are shown. It is clear that the BIE method (Figure~\ref{fig:crowdy_vsN} (right)) agrees with the reference solution to within time-step tolerance for both $tol_1$ and $tol_2$ when the spatial resolution is high enough. For $tol_1$, $36$ panels, i.e. $N_k=576$ is sufficient to reach the set tolerance. For $tol_2$, $46$ panels, i.e. $N_k=736$, are needed. For the validation method the spatial errors dominate the solution and it is first at around $N_V=2000$ points that the method converges to the tolerance $10^{-6}$. To achieve also $tol_2$, the spatial discretisation needs to be increased further.

%% file: srcfiles/surfvalid.tex

\subsection{A pair of surfactant-covered bubbles in extensional flow}
\noindent Similarly to a pair of clean bubbles, the semi-analytical solutions in \S\ref{sec:cleanpair} for a pair of surfactant-covered bubbles is computed and compared against the BIE method of this paper.

\begin{figure}[h!]
	\centering
	\includegraphics[width=0.45\textwidth]{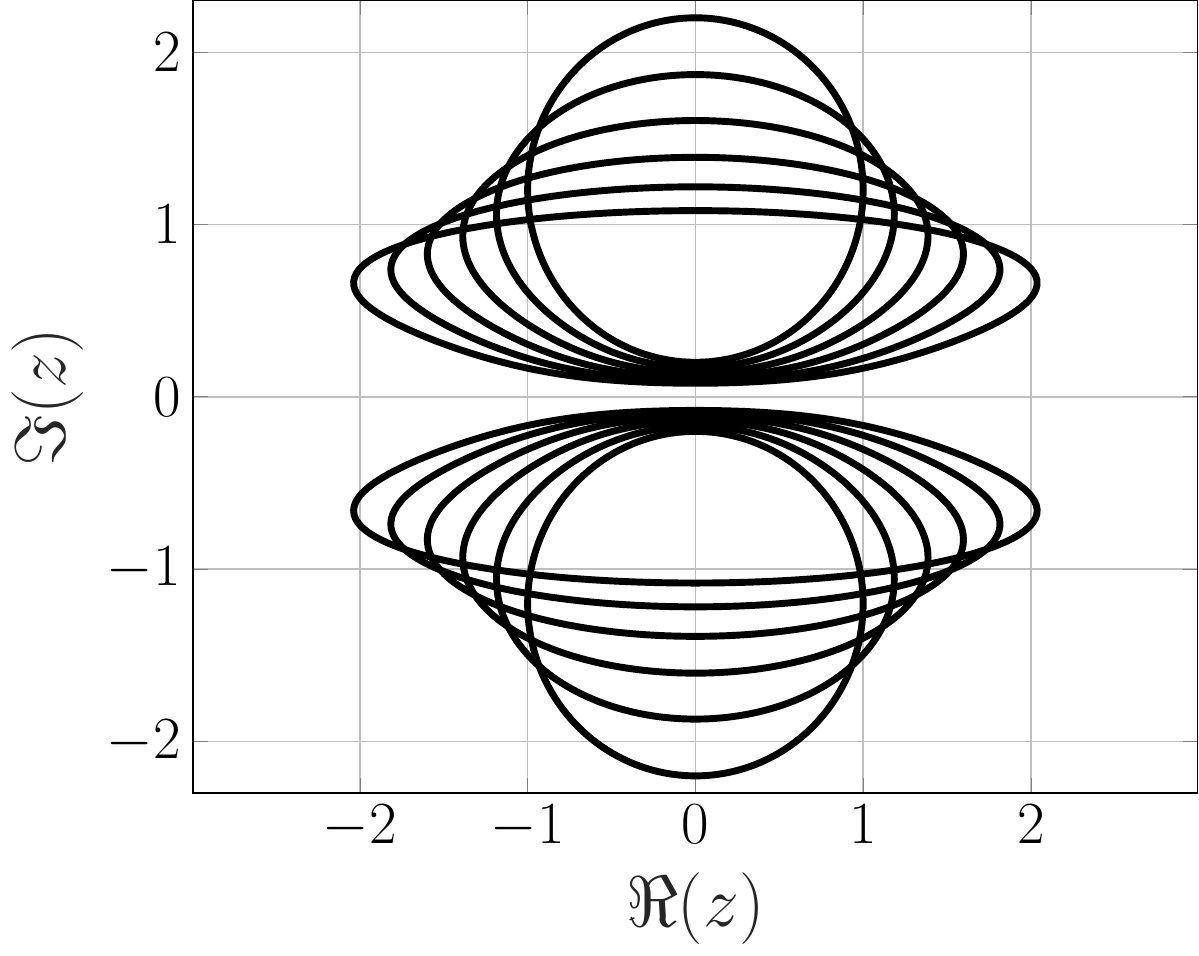}%
	\hspace{1mm}
	\includegraphics[width=0.45\textwidth]{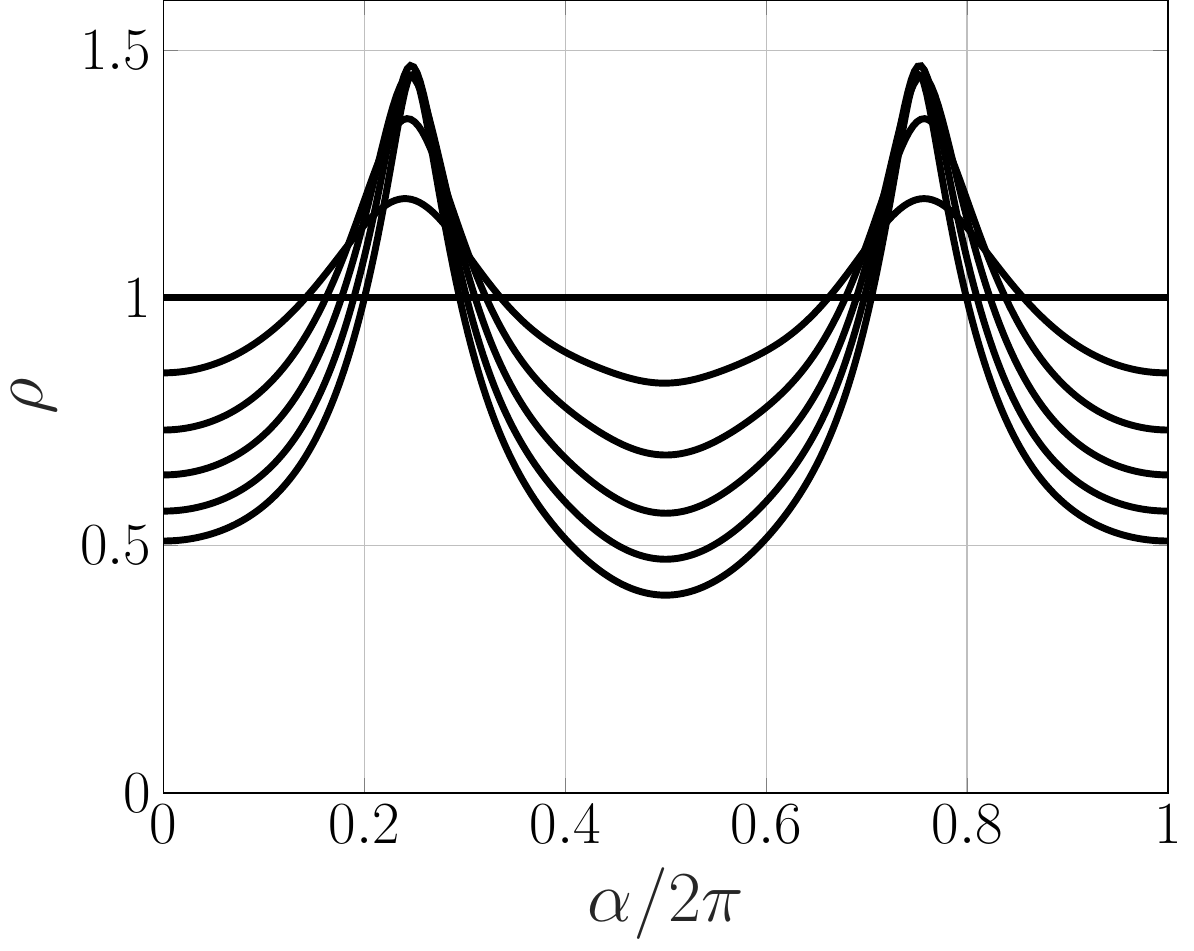}%
	\caption{Pair of surfactant-covered bubbles deforming over time. Left: interface position from time $t=0$ to $t=1$, at $dt=0.2$ intervals, computed with BIE method. Right: surfactant concentration at corresponding times.}
	\label{fig:surfvalid_alltime}
\end{figure}

\begin{figure}[h!]
	\centering
	\includegraphics[width=0.44\textwidth]{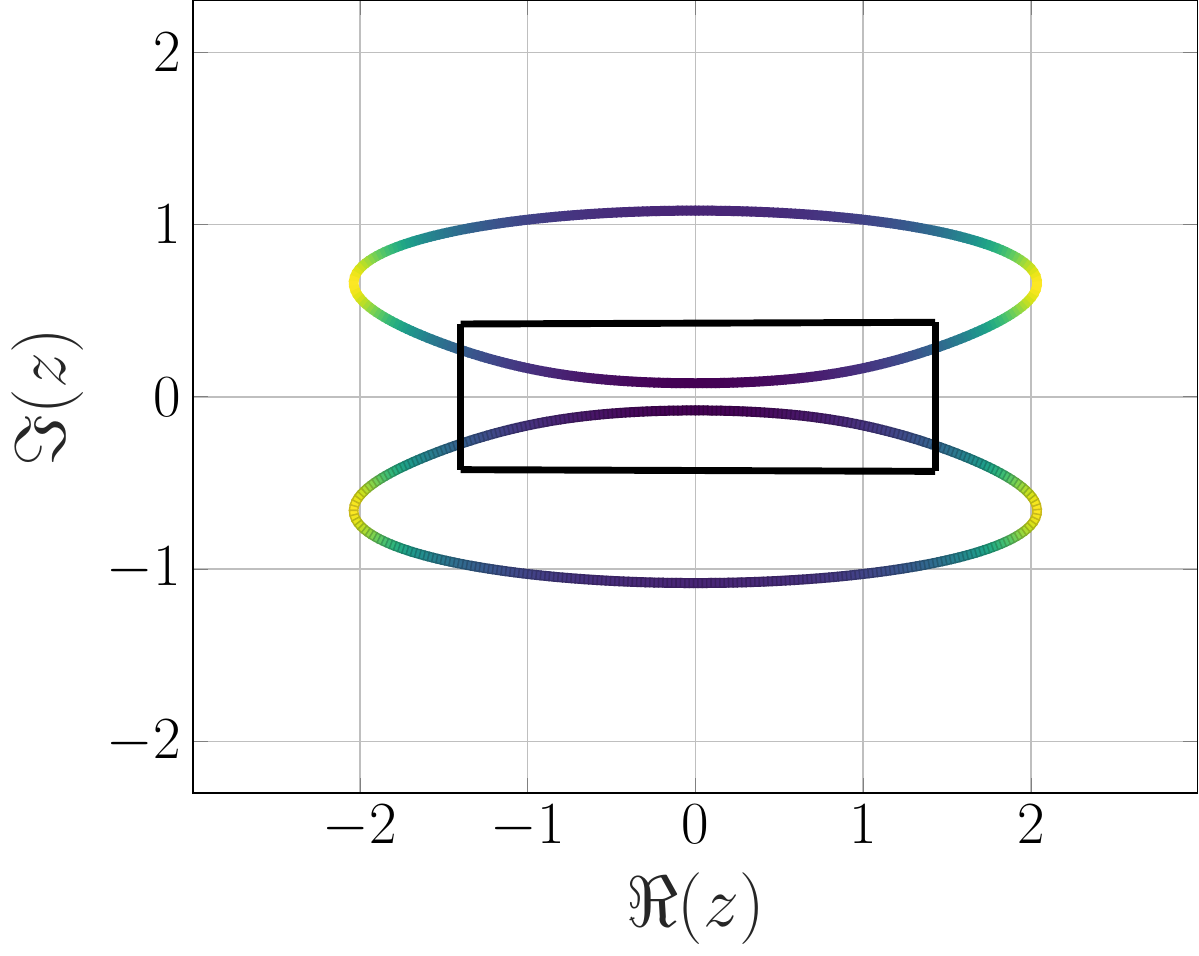}%
	\includegraphics[width=0.53\textwidth]{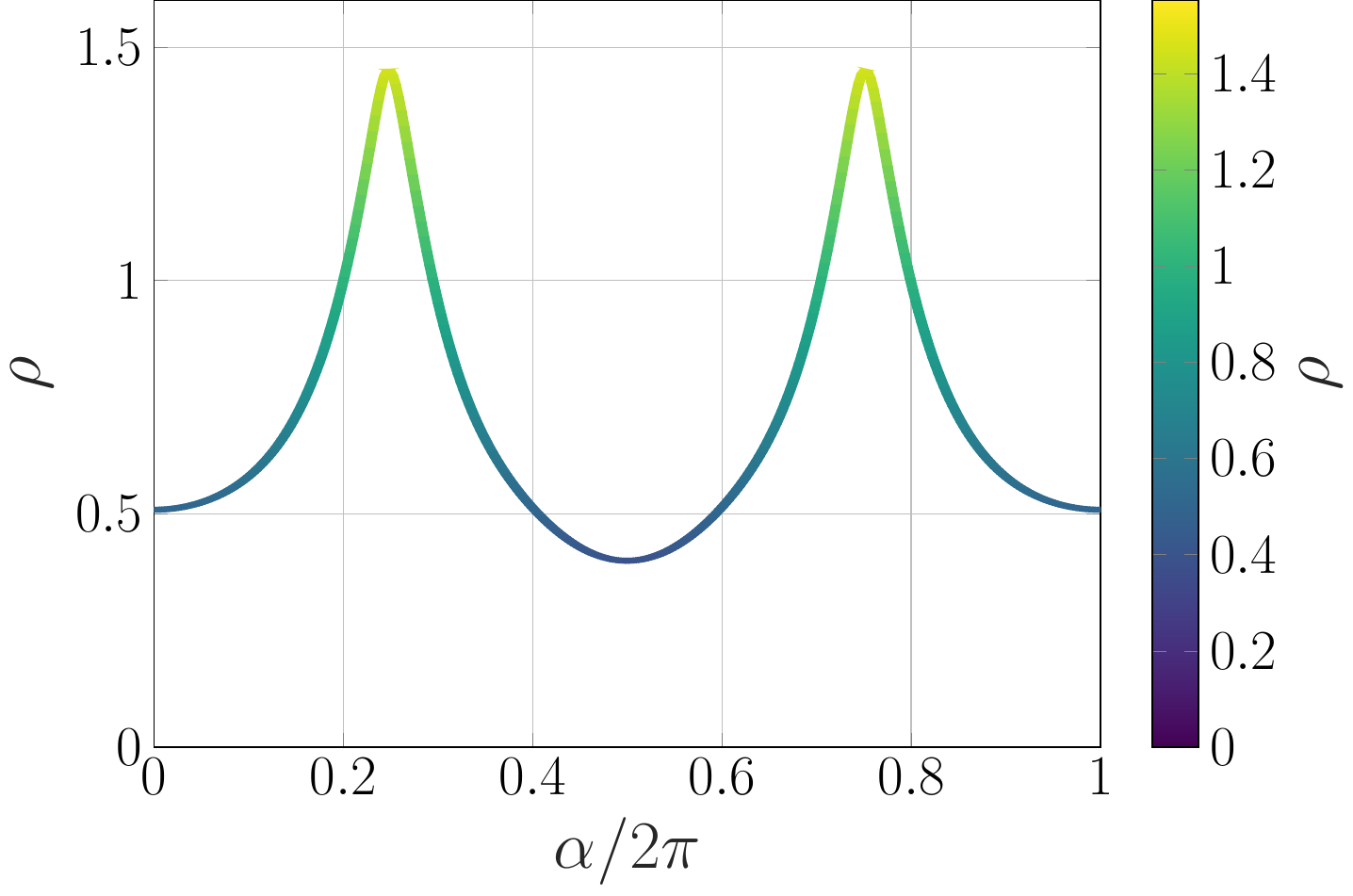}%
	\caption{Surfactant-covered bubbles at time $t=1$. Left: interface position at final time. Cut-off region shown in black box. Right: surfactant concentration at final  time.}
	\label{fig:surfvalid_endtime}
\end{figure}

This case provides a test case for the coupling of interface and surfactant concentration, as well as for both diffusion and convection of the surfactants along the interface. The bubbles are initially circular, covered with a uniform surfactant concentration $\rho_0=1$ and centred around $\pm1.201i$ which corresponds to $\Phi(0) = 0.2875$. In this problem, the elasticity number $E=0.5$, Peclet number $Pe_\Gamma=10$ and the bubbles are placed in an extensional flow with Capillary number $Q=0.5$. The bubbles will deform until time $t=1$, and their change over time can be seen in Figure~\ref{fig:surfvalid_alltime} (left) together with the change of surfactant concentration over time (right). The interface position and surfactant concentration at the final time $t=1$ can be seen in Figure~\ref{fig:surfvalid_endtime}. At time $t=1$ the minimum distance between the bubbles is $0.16$. Similarly to the clean case, only the region where the bubbles are close to each other will be studied, corresponding to the region between $\alpha=\frac{2\pi}{3}$ and $\alpha=\frac{4\pi}{3}$, shown in black in Figure~\ref{fig:surfvalid_endtime} (left).

\begin{figure}[h!]
	\centering
	\includegraphics[width=0.45\textwidth]{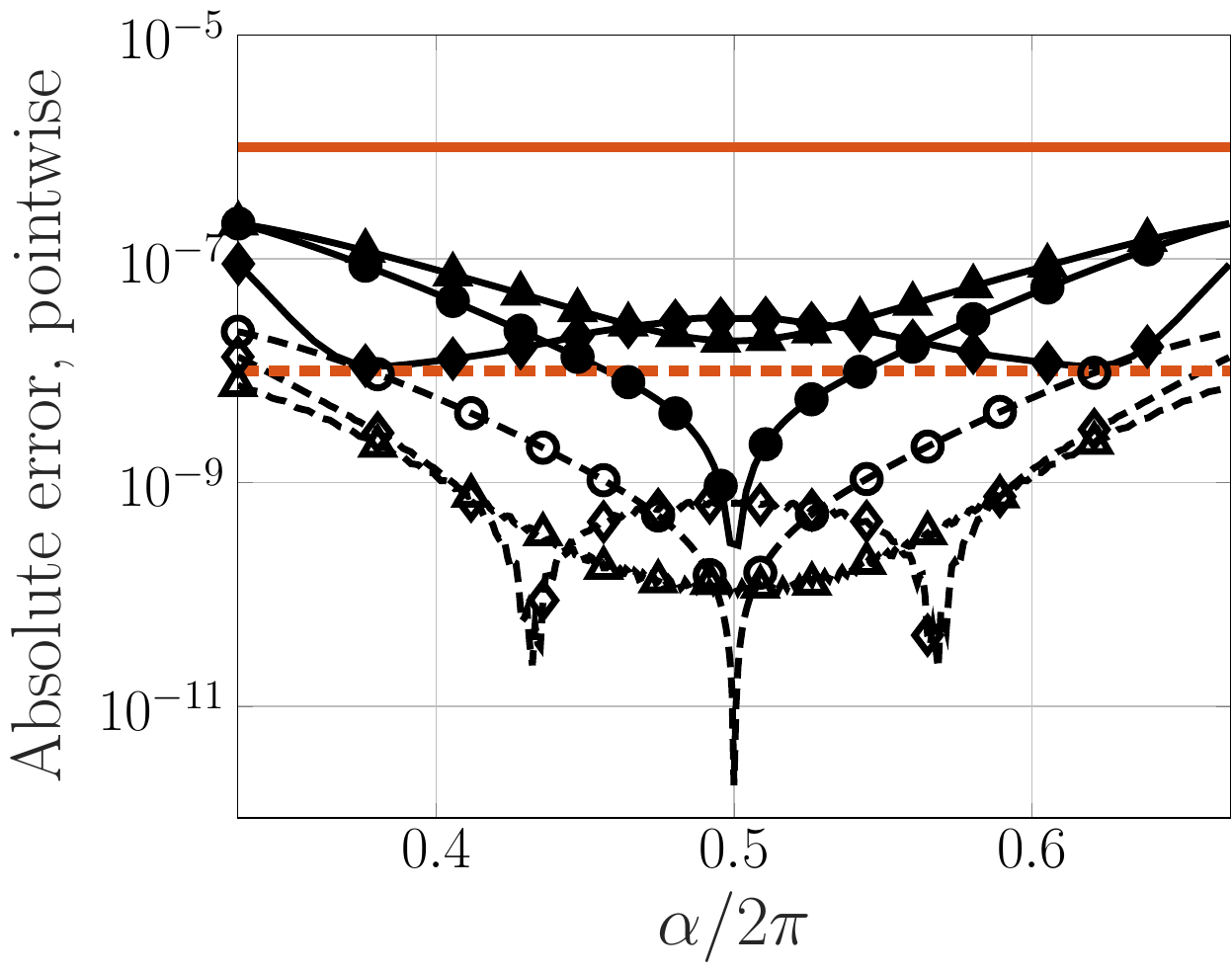}%
	\hspace{1.5mm}
	\includegraphics[width=0.45\textwidth]{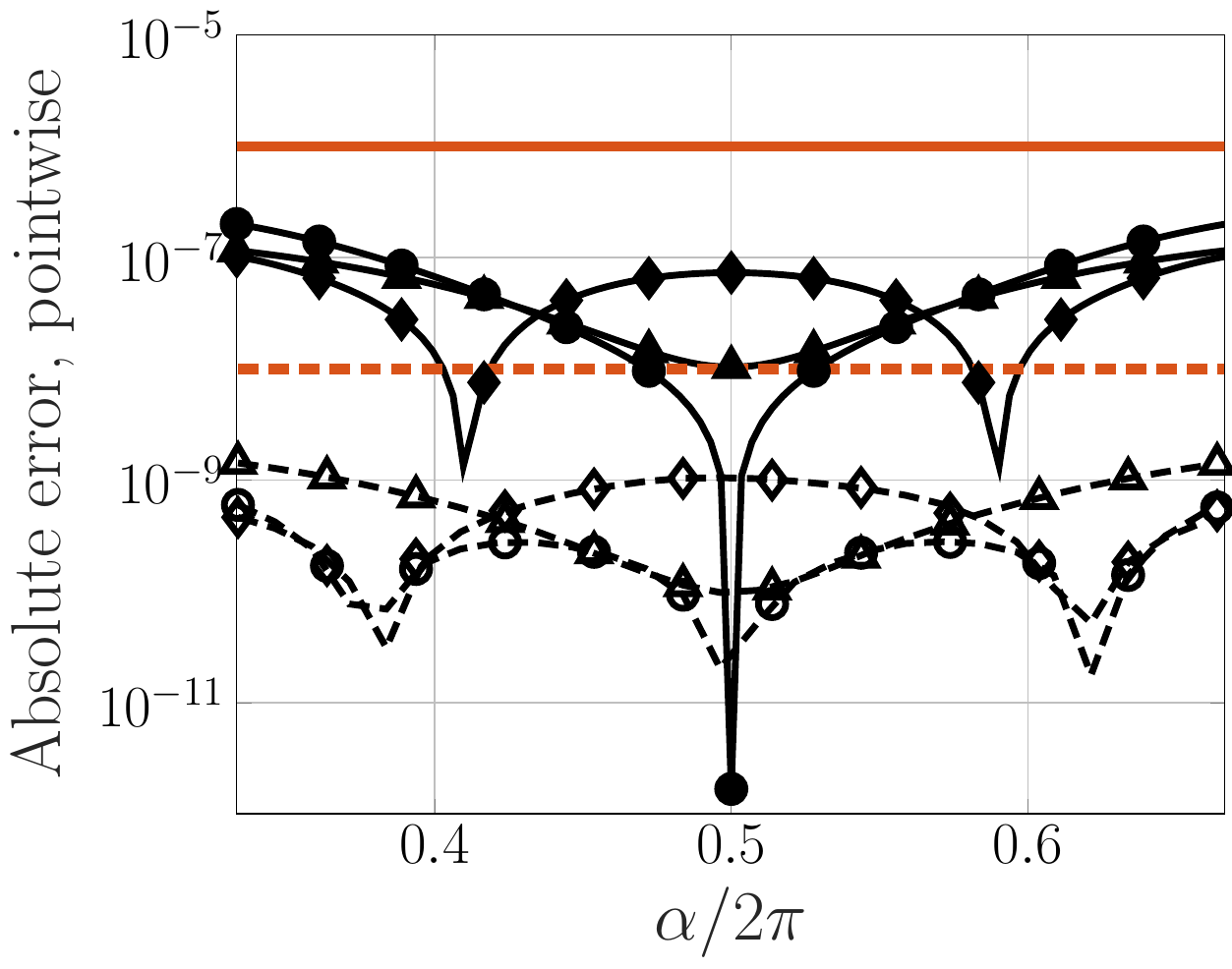}%
	\caption{Point-wise absolute difference in $x$- and $y$-coordinates and surfactant concentration for validation method (left) and BIE method (right) against reference solution. Solid lines/markers show time-step tolerance $10^{-6}$, dashed for tolerance $10^{-8}$. BIE method computed with $576$ and $800$ points per bubble respectively. Validation method computed with $2049$ and $8193$ points respectively. Markers: $\circ$, $\vartriangle$ and $\diamond$ corresponds to $e_x$, $e_y$ and $e_\rho$ respectively.}
	\label{fig:surfvalid_p14o15}
\end{figure}

\begin{figure}[h!]
	\centering
	\includegraphics[width=0.45\textwidth]{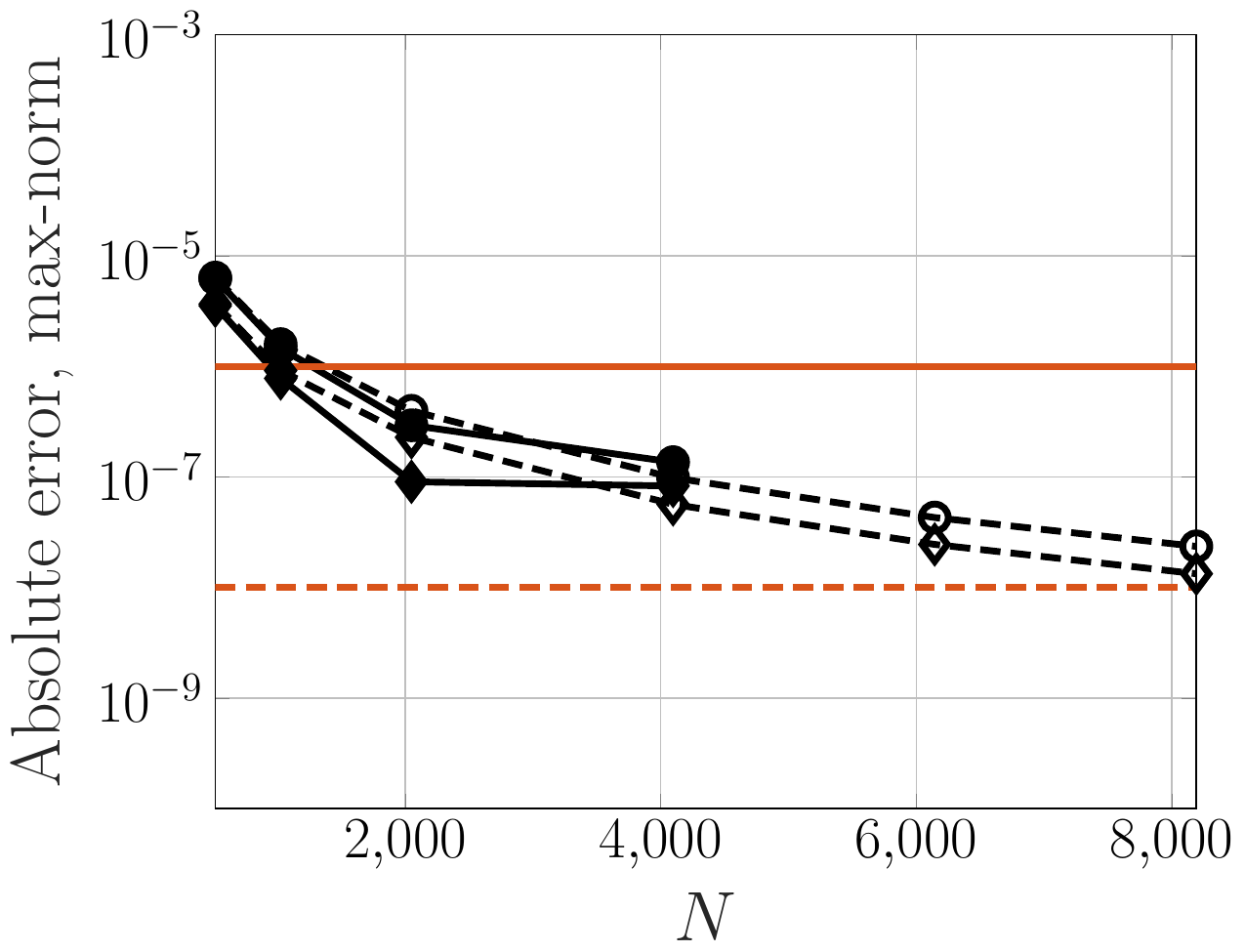}%
	\hspace{1.5mm}
	\includegraphics[width=0.45\textwidth]{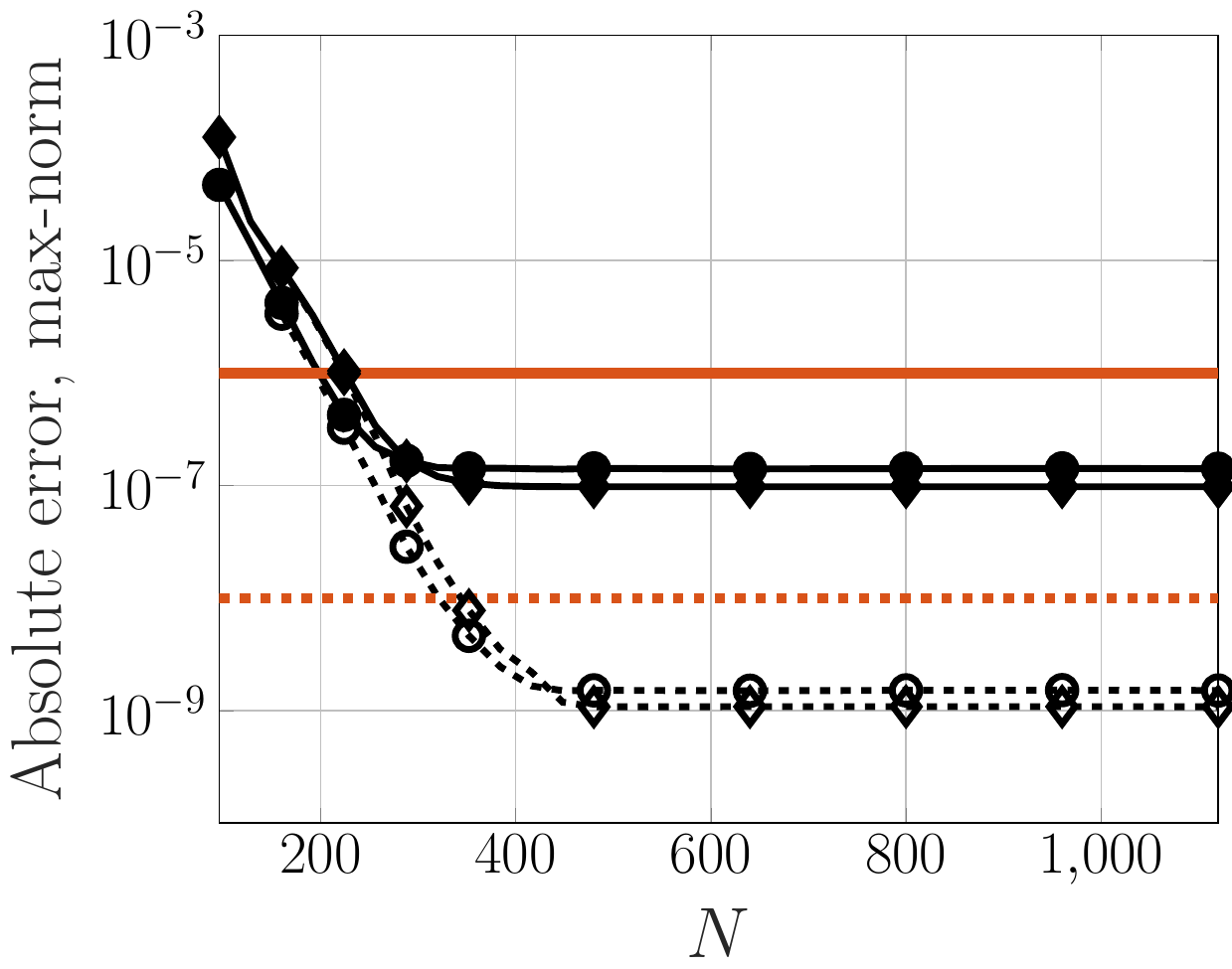}%
	\caption{Error as a function of number of discretisation points for validation method (left) and BIE method (right). Solid lines/markers show results for time-step tolerance $10^{-6}$ and dashed lines tolerance $10^{-8}$. The spatial adaptivity of the BIE method was turned off to keep $N$ constant throughout the simulation. Markers: $\circ$ and $\diamond$ corresponds to $e_z$ and $e_\rho$ respectively.}
	\label{fig:surfvalid_vsN}
\end{figure}

To compare the two methods, a reference solution computed with the BIE method is used. It is computed with a time-step tolerance of $10^{-10}$ and $\Delta s\approx 6.5\cdot10^{-4}$. In Figure~\ref{fig:surfvalid_p14o15} the point-wise absolute error between the two methods and the reference solution is shown, for two time-step tolerances: $tol_1 = 10^{-6}$ (solid lines) and $tol_2 = 10^{-8}$ (dashed lines). The errors considered are those in $x$- and $y$-coordinates, $e_x$ and $e_y$, as well as those in surfactant concentration $e_\rho$. The BIE method was computed with $576$ discretisation points per bubble for $tol_1$ and $800$ for $tol_2$ and all errors stay within time-step tolerance. Note that the spatial adaptivity of the BIE method was turned off for these simulations. In contrast, for the validation method $N_V=2049$ points was used for tolerance $10^{-6}$ and $N_V=9001$ points for $10^{-8}$. Again, for this discretisation the method does not meet the stricter tolerance, due to the spatial accuracy towards the end of the cut-off box. It is clear, however, that given enough spatial resolution the errors would come down also to $tol_2$. In Figure~\ref{fig:surfvalid_vsN} the convergence of the two methods is studied. To reach the set tolerances with the BIE method, $36$ and $50$ panels per bubble are needed for $tol_1$ and $tol_2$ respectively. In comparison, the validation method needs approximately $2049$ points for $tol_1$ and more than $8192$ to reach $tol_2$ due to spatial resolution.

%% file: srcfiles/conclusions2.tex

\section{Conclusions}
\label{sec:conc}
\noindent This paper presents a highly accurate boundary integral method to simulate deforming droplets in Stokes flow. The method can be used for both droplets and inviscid bubbles, both which may be covered by insoluble surfactants. The boundary integral method is coupled with a pseudo-spectral method for surfactant concentration, and together they are spectrally accurate in space.

The errors introduced by the numerical evaluation of the integrals in the boundary integral equation when droplets get close to each other are accurately estimated using contour integrals. These errors are removed using a special quadrature, which enables the method to simulate close drop-drop interactions very accurately.

Given a sufficient spatial resolution, the accuracy of the method is limited only by the set tolerance of the adaptive time-stepping scheme. To the knowledge of the authors, this method is more accurate than other methods to simulate Stokes flow in 2D currently available.

By presenting a set of easily accessible algorithms based on exact and semi-analytical solutions, the hope of the authors is that this will set a standard of validation for any proposed new method for surfactant-laden drops in 2D.

%% file: srcfiles/acknowledgement.tex

\section{Acknowledgement}
\label{sec:acknow}
\noindent The authors gratefully acknowledge the support by the G\"{o}ran Gustafsson Foundation for Research in Natural Sciences and Medicine (S.P. and A.K.T) as well as the support from the National Science Foundation grant DMS-1412789 (M.S.) and by the Knut and Alice Wallenberg foundation under grant no. 2014-0338 (M.S. and A.K.T).

%% file: srcfiles/appendix_est.tex

\subsection{Error estimates}
\label{sec:app_est}
\noindent As discussed in \S\ref{sec:est}, the integrals in equations \eqref{eq:bie_mu} and \eqref{eq:bie_u} both contain terms $(z_j-z_i)^{-1}$, which is problematic when $\|z_i-z_j\|\ll 1$, $i\neq j$. The integrals are then said to become near-singular, which is a numerical problem where large errors are introduced through under-resolving the integrand. To see how the errors grow as an evaluation point $z_0\in\Omega$ is approaching an interface $\Gamma$, the quadrature errors of such near-singular integrals can be approximated \cite{AfKlinteberg2017,klintebergadapt}.

In the following, the error of the near-singular integrals when solving Stokes equations in a domain with non-deforming boundaries, $\Omega$, is considered.  Similarly to the BIE formulation in \S\ref{sec:bieform}, a complex-valued density is introduced and solved for on the boundary $\Gamma$ through
\begin{align}
	\mu(z) + \dfrac{1}{\pi} \int_\Gamma \mu(\tau)\Im\left\{\dfrac{d\tau}{\tau-z}\right\} - \dfrac{1}{\pi} \int_\gamma \overline{\mu(\tau)}\dfrac{\Im\left\{ d\tau\left(\overline{\tau}-\overline{z}\right)\right\}}{\left(\overline{\tau}-\overline{z}\right)^2} = f(z), \; z\in\Gamma,
	\label{eq:est_mu}
\end{align} 
where $f(z)$ is the prescribed boundary condition of $\Gamma$. Setting the right hand side $f(z)$, $z\in\Gamma$, as a sum of Stokeslets for point sources surrounding $\Omega$, the analytical solution $u(z)=f(z)$ for all $z\in\Omega$. For a given set of point sources: $f_1=4\pi+4i\pi$ located at $x_1=1.1+1.3i$, $f_2=\pi-2i\pi$ at $x_2= -1.4-1.3i$ and $f_3=-0.5\pi+3.5i\pi$ at $x_3=1.3-0.75i$, the solution is given in Figure~\ref{fig:est_domu} (right). The domain and the location of the point sources can be seen in Figure~\ref{fig:est_domu} (left). The solution $u(z_0)$ for $z_0\in\Omega$ is computed through
\begin{align}
u(z_0) = -\dfrac{i}{\pi}\int_\Gamma \mu(\tau) \Im\left(\dfrac{d\tau}{\tau-z_0}\right) + \dfrac{i}{\pi}\int_\Gamma \overline{\mu(\tau)}\dfrac{\Im\left\{d\tau(\overline{\tau}-\overline{z_0})\right\}}{(\overline{\tau}-\overline{z_0})^2}.
\label{eq:est_u}
\end{align}
Comparing \eqref{eq:est_u} to \eqref{eq:bie_u}, it is clear that the only difference in the integrals between the two is the imaginary part in the first integral in \eqref{eq:est_u}. This should not affect the near-singular behaviour of the integral. In the following \eqref{eq:est_mu} is solved to high accuracy, and the errors of the solution will come from solely the evaluation of \eqref{eq:est_u}.

Discretising the integrals in \eqref{eq:est_u} with a composite 16-point Gauss-Legendre quadrature and computing $u(z)$ for $z\in\Omega$, the error compared to the analytical solution is shown in Figure~\ref{fig:est_err}. Far away from the boundary $\Gamma$, the error is very small but as $z_0$ approaches $\Gamma$ the error increases. In the region where the errors are large, special treatment is needed for an accurate solution; the method for this is explained in \S\ref{sec:specq}. In order to decide where such special treatment is needed, it is important to know how the errors behave. For Laplace's and Helmholtz equations, estimates of the quadrature error has been derived by \citeauthor{AfKlinteberg2017} \cite{AfKlinteberg2017}, \cite{klintebergadapt}, based on contour integration and calculus of residues. This approach will be followed here.

Regarding the integral expression in \eqref{eq:est_u}, the error at a point $z_0\in\Omega$ can be seen as the sum of the errors when integrating over each panel $\Gamma_i$:
\begin{align*}
	e(z_0) =\sum_{i=1}^{N_{panels}}e_i(z_0).
\end{align*}
The errors will be estimated on a point-panel basis, where the error $e_i(z_0)$ is estimated by $R_i(z_0)$. For convenience only one panel $\Gamma_i$ will be regarded in the following derivation.
\\ \\
In \cite{AfKlinteberg2017,klintebergadapt} estimates are shown for integrals of the type
\begin{align}
	I = \int_{-1}^1 f(\xi)d\xi,
	\label{eq:est_I}
\end{align}
where the meromorphic function $f(\xi)$ is analytic on $[-1,1]$ and has a pole $\xi_0\in\mathbb{C}$ of order $m+1$. To consider the integrals in \eqref{eq:est_u}, the panel $\Gamma_i$ is mapped to the real line $[-1,1]$ by a numerically constructed mapping $\eta(\xi)$. Also, the point $\xi_0\in\mathbb{C}$ that corresponds to each evaluation point $z_0\in\Omega$ needs to be found. This procedure is explained in \cite{klintebergadapt}. 

To compute the estimates, first the integrals are rewritten on the form of \eqref{eq:est_I}. For the first integral in \eqref{eq:est_u}, this is straightforward:
\begin{align}
\begin{split}
J_1(z_0) := \int_{\Gamma_i}\mu(\tau)\Im\left(\dfrac{d\tau}{\tau-z_0}\right) &= \Im\left(\int_{-1}^1\mu(\xi)\dfrac{\eta^\prime(\xi)}{\eta(\xi)-\eta(\xi_0)}d\xi\right) \\ 
&=: \Im\left(\int_{-1}^1 f_1(\xi) d\xi\right).
\end{split}
\label{eq:est_J1}
\end{align}
The second integral is reformulated as follows
\begin{align}
	\int_\Gamma \overline{\mu(\tau)}\dfrac{\Im\left\{ (\overline{\tau}-\overline{z_0})d\tau\right\}}{(\overline{\tau}-\overline{z_0})^2} = \dfrac{i}{2}\overline{\int_{\Gamma}\dfrac{\mu(\tau)\overline{n}_\tau^2d\tau}{\tau-z_0}} + \frac{i}{2}\overline{\int_{\Gamma}\dfrac{\mu(\tau)(\overline{\tau}-\overline{z_0})}{(\tau-z_0)^2}},
	\label{eq:numm_integral2}
\end{align}
where $n_\tau$ is the (outward) normal at point $\tau$. Then,
\begin{align}
\begin{split}
J_2(z_0) := \int_\Gamma \dfrac{\mu(\tau)\overline{n}_\tau^2d\tau}{\tau-z_0} &=  \int_{-1}^1 \dfrac{\mu(\xi)}{\eta(\xi)-\eta(\xi_0)}\left(\overline{\dfrac{i\eta^\prime(\xi)}{|\eta^\prime(\xi)|}}\right)^2\eta^\prime(\xi)d\xi \\
&=: \int_{-1}^1 f_2(\xi)d\xi,
\end{split}
\label{eq:est_J2}
\end{align}
and 
\begin{align}
\begin{split}
J_3(z_0) := \int_\Gamma \dfrac{\mu(\tau)(\overline{\tau}-\overline{z_0})}{(\tau-z_0)^2} &= \int_{-1}^1 \dfrac{\mu(\xi)(\overline{\eta(\xi)}-\overline{\eta(\xi_0)})}{(\eta(\xi)-\eta(\xi_0))^2}\eta^\prime(\xi)d\xi \\
&=: \int_{-1}^1 f_3(\xi)d\xi.
\end{split}
\label{eq:est_J3}
\end{align}
The whole integral over $\Gamma_i$ can be written as
\begin{align}
I_i(z_0) = -\dfrac{i}{\pi}J_1(z_0) + \dfrac{1}{2\pi}\overline{J_2(z_0)} + \dfrac{1}{2\pi}\overline{J_3(z_0)}.
\end{align}
For each integral, error can be estimated as $R^1[\xi_0]$, $R^2[\xi_0]$ and $R^3[\xi_0]$ respectively. In total, the error of the sum of the integrals over the panel $\Gamma_i$ can then be estimated as the absolute value of 
\begin{align}
	R_i(z_0) = \left|-\dfrac{i}{\pi}R[f_1]+\dfrac{1}{2\pi}\overline{R[f_2]} + \dfrac{1}{2\pi}\overline{R[f_3]}\right|.
	\label{eq:est_Rtot}
\end{align}
In the following, how to obtain each estimate $R[f_1]$, $R[f_2]$ and $R[f_3]$ will be explained.

\citeauthor{klintebergadapt} showed that for the integral \eqref{eq:est_I} the quadrature error $R_n[f]$ can be approximated by
\begin{align}
	R_n[f]\approx -\dfrac{1}{m!}\lim_{\xi\rightarrow \xi_0} \dfrac{d^m}{d\xi^m}\left((\xi-\xi_0)^{m+1}f(\xi)k_n(\xi)\right).
	\label{eq:est_Rn}
\end{align}
Here, $k_n(\xi)$ is the characteristic remainder function, approximated as
\begin{align}
	k_n(\xi) \approx \dfrac{2\pi}{(\xi+\sqrt{\xi^2-1})^{2n+1}},
	\label{eq:est_kn}
\end{align}
where $n$ is the order of the Gauss-Legendre quadrature rule used to approximate $I$, in this paper $n=16$. They also showed that the derivatives of $k_n(\xi)$ can be approximated as 
\begin{align}
	k_n^{(m)}(\xi) \approx k_n(\xi) \left(-\dfrac{2n+1}{\sqrt{\xi^2-1}}\right)^m,
	\label{eq:est_knder}
\end{align}
for small $m$, and that
\begin{align*}
	\lim_{\xi\rightarrow \xi_0} \dfrac{(\xi-\xi_0)^{m+1}}{(\eta(\xi)-\eta(\xi_0))^{m+1}} = \dfrac{1}{(\eta^\prime(\xi_0))^{m+1}}.
\end{align*}

For $J_1$, $f_1$ is a meromorphic function with a pole at $\xi_0$ of order one, and the estimate is defined by \eqref{eq:est_Rn} as 
\begin{align*}
	R[f_1] &= \Im\left( \lim_{\xi\rightarrow \xi_0}(\xi-\xi_0)f_1(\xi)k_n(\xi)\right) \\
	&\approx \Im\left( k_n(\xi_0)\mu(\xi_0) \right),
\end{align*}
where \eqref{eq:est_knder} is used and $k_n(\xi_0)$ can be approximated as in \eqref{eq:est_kn}. Instead of finding the value of $\mu(\xi)$ at the parameter $\xi_0$, it is sufficient to consider the max-norm of $\mu$ on the panel $\Gamma_i$: $\|\mu^i\|_\infty$. Thus the error estimate for $J_1(z_0)$ becomes
\begin{align}
	R[f_1] \approx \Im \left( k_n(\xi_0)\|\mu^i\|_\infty \right).
	\label{eq:est_R1}
\end{align}

Similarly, for the second integral $J_2$ the function $f_2$ is meromorphic with a pole of order one at $\xi_0$. The estimate reads
\begin{align}
	R[f_2] &= \lim_{\xi\rightarrow \xi_0} (\xi-\xi_0)f_2(\xi)k_n(\xi) \nonumber \\
	& \approx k_n(\xi_0)\mu(\xi_0)\overline{n(\xi_0)}^2  \nonumber \\
	& \approx k_n(\xi_0)\|\mu^i\|_\infty\overline{n(\xi_0)}^2.
	\label{eq:est_R2}
\end{align}
Again, \eqref{eq:est_knder} is used. Note that also the value of the normal $n$ needs to be found at $\xi_0$ through the mapping $\eta(\xi)$.

Finally, regard the integral $J_3$. Here, $f_3$ is a meromorphic function with a pole of order two at $\xi_0$. Following \eqref{eq:est_Rn}, the estimate becomes
\begin{align*}
	R[f_3] = \lim_{\xi\rightarrow \xi_0} \dfrac{d}{d\xi}\left((\xi-\xi_0)^2f_3(\xi)k_n(\xi)\right).
\end{align*}
Assuming a good approximation can be obtain by considering only the term with the derivative of $k_n$, see \cite{klintebergadapt}, this can be simplified into
\begin{align}
	R[f_3] &\approx k^\prime_n(\xi_0)\dfrac{1}{\eta^\prime(\xi_0)}\left(\overline{\eta(\overline{\xi_0})}-\overline{\eta(\xi_0)}\right)\mu(\xi_0) \nonumber \\
		&\approx k^\prime_n(\xi_0)\dfrac{1}{\eta^\prime(\xi_0)}\left(\overline{\eta(\overline{\xi_0})}-\overline{\eta(\xi_0)}\right)\|\mu^\|_\infty.
	\label{eq:est_R3)}
\end{align}

Together, the complete error estimate for the integrals in \eqref{eq:est_u} over a panel $\Gamma_i$ can be written as \eqref{eq:est_Rtot}. This error etimate is shown in Figure~\ref{fig:est_est}, where the estimate of errors of order $10^{-p}$ is marked with black level curves for $p=-14,\hdots,0$ with increment $2$.

%% file: srcfiles/appendix_spec.tex

\subsection{Special quadrature}
\label{sec:app_spec}
\noindent An overview of the special quadrature follows here, for more details see \cite{Ojala2015}.

The special quadrature is employed on integrals of two types:
\begin{align*}
	I_1(z) &= \int_{\Gamma_i} \dfrac{f(\tau)d\tau}{\tau-z},  \\
	I_2(z) &= \int_{\Gamma_i} \dfrac{f(\tau)d\tau}{(\tau-z)^2}. 
\end{align*} 
The function $f(\tau)$ is interpolated on the panel $\Gamma_i$, as $f(\tau) \approx \sum_{j=0}^{15} c_j\tau^j$, where the coefficients $c_j$ are computed through solving a Vandermonde system.  For stability, the panel is transformed to have endpoints at $-1$ and $1$. Using the above interpolation, the integrals $I_1$ and $I_2$ can be rewritten as 
\begin{align} 
	I_1(z) \approx& \sum_{j=0}^{15} c_j \int_{-1}^1 \dfrac{\tau^jd\tau}{\tau-z_0} = \sum_{j=0}^{15} c_j p_j, \\
	I_2(z) \approx& \, \alpha \sum_{j=0}^{15}c_j\int_{-1}^1\dfrac{\tau^jd\tau}{(\tau-z_0)^2} = \alpha \sum_{j=0}^{15} c_j q_j.
\end{align}
Here $z_0$ is the target point $z$ under the same transformation as the one applied to the panel and $\alpha =  \frac{2}{\tau_e - \tau_s}$, where $\tau_s$, $\tau_e$ are the start and end points of the untransformed panel respectively. The numbers $p_j$ and $q_j$ can computed through recursion, where
\begin{align}
\begin{cases}
	p_0 &= \int_{-1}^1 \dfrac{d\tau}{\tau-z_0} = \log(1-z_0) - \log(-1-z_0), \\ 
	p_j &= z_0p_{j-1} + \dfrac{1-(-1)^j}{j}, \; j=1,\hdots, 15,
\end{cases}
\end{align}
and 
\begin{align}
\begin{cases}
	q_0 &= \int_{-1}^{1}\dfrac{d\tau}{(\tau-z_0)^2} = -\dfrac{1}{1+z_0} - \dfrac{1}{1-z_0}, \\ 
	q_j &= z_0q_{j-1} + p_j, \; j=1,\hdots, 15.
\end{cases}
\end{align}
Note that if $z_0$ is within the contour created by the real axis from $-1$ to $1$ and the transformed panel $\Gamma_i$, a residue of $2\pi i$ must be added or subtracted from $p_0$ depending on if $z_0$ has a positive or negative imaginary part respectively.
\\ \\
To apply the special quadrature on the integrals in \eqref{eq:bie_mu} and \eqref{eq:bie_u}, they need to be rewritten into the forms of $I_1$ and $I_2$. The second integral can be reformulated as in \eqref{eq:numm_integral2}, i.e.
\begin{align*}
	\int_\Gamma \overline{\mu(\tau)}\dfrac{\Im\left\{ (\overline{\tau}-\overline{z_0})d\tau\right\}}{(\overline{\tau}-\overline{z_0})^2} = \dfrac{i}{2}\overline{\int_{\Gamma}\dfrac{\mu(\tau)\overline{n}_\tau^2d\tau}{\tau-z_0}} + \frac{i}{2}\overline{\int_{\Gamma}\dfrac{\mu(\tau)(\overline{\tau}-\overline{z_0})}{(\tau-z_0)^2}}.
\end{align*}
All integrals are then on the required form, with $f_1(\tau) = \mu(\tau)\tau^\prime$, $f_2(\tau) = \mu(\tau)\bar{n}_\tau^2\tau^\prime$ and $f_3(\tau) = \mu(\tau)(\overline{\tau}-\overline{z})\tau^\prime$ respectively. The improved error computed with $50$ Gauss-Legendre panels can be seen in Figure~\ref{fig:numm_specquad}. Using the special quadrature it is possible to evaluate $u(z)$ for points close to the boundary, with errors below $10^{-10}$.